\documentclass[12pt,oneside,reqno]{amsart}

\usepackage[margin=1in]{geometry}

\usepackage{geometry}

\usepackage{multirow}
\usepackage{shortcuts}
\usepackage{stmaryrd}
\usepackage{lipsum}
\usepackage{graphicx}
\usepackage{subfigure}
\usepackage{setspace}
\usepackage{indentfirst}

\graphicspath{{./figures/}} 
\usepackage{hyperref}
\hypersetup{
    colorlinks,
    citecolor=black,
    filecolor=black,
    linkcolor=black,
    urlcolor=black
}

\newcommand{\cE}{\mathcal{E}}

\newcommand{\quand}{\quad \text{and} \quad}

\newcommand{\dd}{\mathrm{d}}

\newcommand{\ui}{\underline{i}}
\newcommand{\ipi}[2]{\ip{#1}{#2}_{I_i}}

\theoremstyle{definition}

\usepackage{etoolbox} 
\usepackage{adjustbox} 

\BeforeBeginEnvironment{tabular}{\begin{adjustbox}{scale=1}} 
	\AfterEndEnvironment{tabular}{\end{adjustbox}}

\usepackage{chngcntr}

\counterwithin*{equation}{section}

\title[Reducing polynomial degree for RKDG]{Reducing polynomial degree by one for inner-stage operators affects neither stability type nor accuracy order of the Runge--Kutta discontinuous Galerkin method}

\author[Z.~Sun]{Zheng Sun}
\address{Department of Mathematics,
	The University of Alabama,
	Box 870350,
	Tuscaloosa, AL 35487, USA}
\email{zsun30@ua.edu}

\newcommand{\nm}[1]{%
	\relax\ifinner
	\| #1 \|
	\else
	\left\| #1 \right\|
	\fi
}

\newcommand{\RRR}{\widetilde{R}}
\newcommand{\Pip}{\Pi_{\perp}}
\newcommand{\Pig}{\Pi_{G}}
\newcommand{\Pis}{\Pi_{\star}}
\newcommand{\Pik}{\Pi}
\newcommand{\Piz}{\Pi_{k-1}}
\newcommand{\PiPi}{\overline{\Pi}}
\newcommand{\LL}{\widetilde{L}}

\newcommand{\Mod}[1]{\ (\mathrm{mod}\ #1)}
\begin{document}


\begin{abstract}
	The Runge--Kutta (RK) discontinuous Galerkin (DG) method is a mainstream numerical algorithm  for solving hyperbolic equations. In this paper, we use the linear advection equation in one and two dimensions as a model problem to prove the following results: For an arbitrarily high-order RKDG scheme in Butcher form, as long as we use the $P^k$ approximation in the final stage, even if we drop the $k$th-order polynomial modes and use the $P^{k-1}$ approximation for the DG operators at all inner RK stages, the resulting numerical method still maintains the same type of stability and convergence rate as those of the original RKDG method. Numerical examples are provided to validate the analysis. The numerical method analyzed in this paper is a special case of the Class A RKDG method with stage-dependent polynomial spaces proposed in \cite{chen2024runge}. Our analysis provides theoretical justifications for employing cost-effective and low-order spatial discretization at specific RK stages for developing more efficient DG schemes without affecting stability type and accuracy order of the original method.
\end{abstract}

\maketitle 

\section{Introduction}

In this paper, we perform stability analysis and optimal error estimates of a class of Runge--Kutta (RK) discontinuous Galerkin (DG) schemes for solving a linear hyperbolic equation. In these schemes, the DG operators at all inner RK stages adopt polynomial approximations of one degree lower than those of the original RKDG schemes. The numerical method analyzed in this paper belongs to Class A of the RKDG method with stage-dependent polynomial spaces (sd-RKDG method), which we proposed in \cite{chen2024runge}. Hence, we will refer to it as the sdA-RKDG method. Its general formulation is given in \eqref{eq:sdA-rkdg-Butcher} and \eqref{eq:sdA-RKDG}.

\subsection{RKDG method and its variant.} The DG method is a class of finite element methods that utilize discontinuous piecewise polynomial spaces. It was initially introduced by Reed and Hill in the 1970s for solving steady-state transport equations \cite{reed1973triangular}. Later, in a series of works by Cockburn et al. \cite{rkdg1,rkdg2,rkdg3,rkdg4,rkdg5}, the DG method was integrated with the Runge--Kutta (RK) time discretization method for solving time-dependent hyperbolic conservation laws. The resulting fully discrete method is referred to as the RKDG method. The RKDG method features many favorable properties, such as good stability, high-order accuracy, flexibility for handling complex geometries, and capability for sharp capturing of shock waves. It has become one of the mainstream numerical methods for computational fluid dynamics \cite{cockburn2001runge,cockburn2012discontinuous}. 

The construction of the RKDG scheme follows the classical method-of-lines framework. For example, consider the initial value problem of the linear advection equation
\begin{equation}\label{eq:adv}
	u_t + \pmb{\beta}\cdot \nabla u = 0, \quad u = u(\pmb{x},t), \quad \pmb{x} \in \Omega \subset \mathbb{R}^d, \quad d = 1,2, \quad 0 < t \leq  T,
\end{equation}
where $\pmb{\beta}\in \mathbb{R}^d$ is a constant or a constant vector. To avoid unnecessary technicalities, we will consider a periodic rectangular domain with Cartesian meshes. If one wishes to construct a second-order RKDG scheme, the first step is to apply the $P^1$-DG spatial discretization to the linear advection equation \eqref{eq:adv}, which yields a semi-discrete scheme in the form of an ordinary differential equation (ODE) system
\begin{equation}
	\frac{\dd}{\dd t} u_h = L_h u_h,
\end{equation}
where $u_h$ is the semi-discrete solution, and $L_h$ is the DG operator. To obtain the fully discrete scheme, a second-order RK method is used for time discretization. For example, using the explicit midpoint rule with a time step size of $\tau$, we have
\begin{subequations}\label{eq:rkdg2}
	\begin{align}
		u_h^{(2)} =& u_h^{n} +\frac{\tau}{2} L_h u_h^n,\label{eq:rkdg2-1}\\
		u_h^{n+1} =& u_h^{n} + {\tau} L_h u_h^{(2)}.
	\end{align}
\end{subequations}
This formulation results in a stable second-order RKDG scheme.

However, in our recent work \cite{chen2024runge}, we observe that even if one uses the first-order $P^0$ approximation in the first stage \eqref{eq:rkdg2-1}, the resulting numerical scheme is still stable and admits optimal accuracy. Let $\LL_h$ be the $P^0$ approximation of $L_h$. In other words, it is the $P^1$-DG operator without the evaluation of linear modes. See \eqref{eq:def-LL_h} for the rigorous setup. Then the scheme
\begin{subequations}\label{eq:sdArkdg2}
	\begin{align}
		u_h^{(2)} =& u_h^{n} +\frac{\tau}{2} \LL_h u_h^n,\\
		u_h^{n+1} =& u_h^{n} + {\tau} L_h u_h^{(2)}
	\end{align}
\end{subequations}
is still stable under the usual Courant–Friedrichs–Lewy (CFL) condition and maintains second-order accuracy. \eqref{eq:sdArkdg2} is an example of the sdA-RKDG scheme \eqref{eq:sdA-rkdg-Butcher}.

In the implementation of the RKDG scheme \eqref{eq:rkdg2}, one needs to evaluate the second-order $P^1$-DG operator $L_h$ twice, while in the implementation of the sdA-RKDG scheme \eqref{eq:sdArkdg2}, one only needs to evaluate a first-order $P^0$-DG operator once and a second-order $P^1$-DG operator once. The ratio of the total numbers of polynomial coefficients one needs to evaluate in \eqref{eq:sdArkdg2} and \eqref{eq:rkdg2} is $d+2:2d+2$, which is $75\%$ in one dimension (1D, also for one-dimensional) and about $66.7\%$ in two dimensions (2D, also for two-dimensional). As one can see, in \eqref{eq:sdArkdg2}, by breaking the method-of-lines structure and allowing DG operators with reduced polynomial order, we can obtain a more efficient numerical scheme without sacrificing stability and optimal order of accuracy.

\subsection{A truncation error argument for optimal accuracy.}\label{sec:trunc} The stability and accuracy of \eqref{eq:sdArkdg2} were analyzed through the Fourier approach in \cite{chen2024runge}. For an intuitive explanation on why low-order approximations of inner stages do not cause order reduction, the following local truncation error analysis for the finite difference schemes may offer some insights.\footnote{However, this argument may not apply to the DG spatial discretization. See Section \ref{sec:lim} for further discussion.} Indeed, now consider $L_h$ and $\LL_h$ both to be finite difference operators. One can regard $\LL_h$ as a first-order perturbation of the second-order operator $L_h$, namely, $\LL_h u = L_h u + \mathcal{O}(h)$, where $u$ is a smooth function, $h$ is the spatial mesh size, and $\mathcal{O}(h)$ corresponds to the pointwise error. Replacing $u_h^n(\pmb{x})$ by $u^n(\pmb{x}) = u(\pmb{x},t^n)$ in \eqref{eq:sdArkdg2} and substituting in this perturbation relationship, one obtains the local truncation error for \eqref{eq:sdArkdg2}
\begin{equation}\label{eq:intro-taylor-sdArkdg2}
	\begin{aligned}
	\mathrm{LTE} =& u^{n+1} - \left(I + \tau L_h + \frac{\tau^2}{2} L_h\LL_h\right)u^n
	= u^{n+1} - \left(I + \tau L_h + \frac{\tau^2}{2} L_h^2\right)u^n + \mathcal{O}(\tau^2 h). 
	\end{aligned}
\end{equation}
As can be seen, the perturbation in the first stage will be multiplied by $\tau^2$ in \eqref{eq:intro-taylor-sdArkdg2}, resulting in an additional error term $\mathcal{O}(\tau^2h)$. Without this term, \eqref{eq:intro-taylor-sdArkdg2} retrieves the local truncation error of the scheme \eqref{eq:rkdg2}. Therefore, the perturbation from $L_h$ to $\LL_h$ will contribute only to a third-order local truncation error $\mathcal{O}(\tau^2h)$ under the usual time step constraint $\tau \leq Ch$, which does not affect the global-in-time second-order convergence.

Indeed, similar analysis of the low-order (or low-precision) perturbation of RK stages has been systematically studied in the context of ODEs by Grant in \cite{grant2022perturbed}. See also \cite{burnett2021performance,burnett2023stability}. As a special case of Grant's work, we have the following observation: For an RK method in Butcher form of $r$th order, if one perturbs all inner-stage operators by $\mathcal{O}(\varepsilon)$ and leaves the final stage untouched, the resulting perturbed RK method admits the local truncation error $\mathcal{O}(\tau^{r+1})+\mathcal{O}(\varepsilon \tau^2)$. Thanks to the extra $\tau$ and the fact $\tau \sim h$, one can take the perturbation as large as $\varepsilon\sim h^k$ at inner stages and still achieve optimal accuracy in the global error $\mathcal{O}(\tau^r + h^{k+1})$. In other words, for high-order RK schemes in Butcher form, reducing the accuracy of the inner-stage operators by one may not affect the overall convergence rate. In \cite{chen2023runge}, we used this idea to develop an RKDG method with compact stencils. In \cite{chen2024runge}, we further explored this idea and developed the sd-RKDG method, and Grant's theory may relate to the performance of the Class A schemes of the sd-RKDG method.

\subsection{Limitations of the existing analysis.}\label{sec:lim} Although the Fourier analysis in \cite{chen2024runge} and the local truncation error analysis in \cite{grant2022perturbed} shed light on the stability and accuracy of the sdA-RKDG schemes, these analyses have essential limitations and should not be considered as a complete theoretical justification for the sdA-RKDG method. The Fourier analysis only applies to linear problems with constant coefficients on periodic domains. Moreover, this approach requires to analyze matrix eigenvalues, which are difficult to obtain analytically for high-order schemes or in multiple dimensions. The local truncation error analysis of the perturbed RK method for ODEs in \cite{grant2022perturbed} relies on the assumption that the ODE operator has bounded derivatives. However, in the context of the fully discrete schemes for partial differential equations (PDEs), the operators $L_h$ and $\LL_h$ scale inversely proportional to the mesh size $h$, which violates the essential assumptions in numerical ODE analysis as the limit $h\to 0$ is approached. Although this limitation can be circumvented to cover finite difference methods by analyzing the fully discrete local truncation error (see Section 1.2), the argument cannot be applied to DG methods. Indeed, in the design of sd-RKDG schemes \cite{chen2024runge}, we have encountered cases where the local truncation error analysis mismatches the actual convergence rate for the sd-RKDG schemes. It may indicate that the ODE approach in \cite{grant2022perturbed} may not fully explain the numerical performance of the DG solvers for PDEs. In fact, for the error estimates in this paper, a very different argument will be used for the proof.

\subsection{Literature review.} The classical approach for analyzing finite element schemes is through an energy-type argument, in which one exploits the weak formulation of the numerical scheme and derives norm estimates by taking appropriate test functions \cite{thomee2007galerkin,arnold2002unified,burman2010explicit}. In terms of the DG method for hyperbolic problems, although there are many earlier works on the analysis of the DG schemes for the steady-state transport equations \cite{lesaint1974finite,johnson1986analysis,peterson1991note,brezzi2004discontinuous} and semi-discrete DG schemes for the time-dependent hyperbolic conservation laws \cite{jiang1994cell,shu2009discontinuous}, a systematic stability and error analysis was not readily available until recently. In \cite{zhang2004error,zhang2010stability}, Zhang and Shu analyzed the stability and error estimates of the second- and third-order RKDG schemes. These works also relate to \cite{burman2010explicit} on the study of the stabilized finite element methods by Burman et al. and also to \cite{tadmor2002semidiscrete} on the stability analysis of the third-order RK method for generic linear semi-negative problems by Tadmor. The analysis of RKDG schemes of higher order was hindered by the stability analysis of the fourth-order RK scheme for some time, which was then resolved in \cite{sun2017rk4}. With the experience on the fourth-order RK method \cite{sun2017rk4,ranocha2018l_2}, the stability analysis for both the RK method for generic linear semi-negative problems and for the RKDG method was extended to very high order \cite{sun2019strong,xu2019L2,xu2020superconvergence,sun2022energy}. The stability theory also facilitated the error estimates of the RKDG schemes. See \cite{xu2020superconvergence,xu2020error,sun2020error,ai20222,xu2024sdf}.

\subsection{Challenges, main ideas, and main results.} In this paper, we perform stability analysis and optimal error estimates of the sdA-RKDG method for the linear advection equation in both 1D and 2D. The analysis is conducted for schemes of arbitrarily high order.  

The major challenge in the analysis is that the sdA-RKDG method goes beyond the classical method-of-lines framework and involves two different spatial operators within each time step. While in most previous works, the fully discrete analysis typically only involves the same spatial operator across all stages. From this perspective, the challenge in the analysis of the sdA-RKDG method shares similarities with the analysis of the Lax--Wendroff DG method \cite{sun2017stability} and the RKDG method with stage-dependent numerical fluxes \cite{xu2024sdf}. However, they are also very different in the sense that \cite{sun2017stability,xu2024sdf} mainly involve different flux terms, while the sdA-RKDG method involves operators associated with different polynomial spaces. This difference presents significant challenges: For the sdA-RKDG method, even the solution $u_h^n$ itself may not be in the test function space that defines the low-order operator $\LL_h$ in  \eqref{eq:def-LL_h}.

The main idea to overcome these difficulties is the following. We will work on the strong form, instead of the weak bilinear form, of the sdA-RKDG scheme. We note that $\LL_h = (I-\Pip)L_h$, where $(I-\Pip)$ is the $L^2$ projection to the $P^{k-1}$ space (and $\Pip$ is the projection to the orthogonal complement of the $P^{k-1}$ space in the $P^k$ space). Therefore, the $\LL_h$ operator can be embedded into the usual $P^k$-DG space, which avoids the complication of working with different spaces. For the stability analysis, we reformulate the sdA-RKDG scheme as the original RKDG scheme plus some perturbation terms consisting of compositions of $L_h$ and $\Pip$. The main results for the stability analysis are summarized below.
\begin{itemize}
	\item In each time step, the $L^2$ energy of the sdA-RKDG solution can be bounded by that of the original RKDG solution along with some diminishing jump terms  (Lemma \ref{lem:keykey}). Through this energy estimate, one can deduce that the sdA-RKDG method retains the same type of stability as the original RKDG method (Theorem \ref{thm:key}).
\end{itemize}
For the error estimates, we derive the error equation in a compact form. We design a class of special approximation operators to help control and eliminate the remainder terms in the error equation. These approximation operators serve a role similar to that of interpolation or projection operators in classical finite element analysis, but they are neither interpolation nor projection in general. Their construction is also very different, which incorporates the time step size $\tau$. The main results for the error analysis are summarized below.

\begin{itemize}
	\item We construct a class of approximation operators that are well-defined and admit the optimal approximation estimate under the usual CFL condition (Lemma \ref{lem:newproj-general}). With these new approximation operators, we prove the optimal error estimates of the sdA-RKDG schemes (Theorem \ref{thm:err-g}).
\end{itemize}

We want to remark that the error estimates in this paper are written in a different language from that in previous works \cite{zhang2004error,zhang2006error,xu2020superconvergence,xu2020error}, where the analysis is conducted using the multi-stage weak form of the RKDG schemes. To leverage the stability estimates in Section \ref{sec:stab}, the error estimates in this paper are performed under the compact strong form of the sdA-RKDG method with all stages combined. The language in this paper is more similar to that in the classical work by Brenner et al.  \cite{brenner1982single}. However, despite the differences in notations at first glance, the two approaches are indeed closely interrelated. See Appendix \ref{sec:multistage}. 

\subsection{Significance and novelty.} In this work, we utilize the RKDG method for the linear advection equation as an illustrative example. Through systematic and rigorous analysis, we demonstrate the feasibility of employing cost-effective and coarse spatial discretization at specific RK stages for  developing a numerical scheme that is highly efficient while maintaining the stability and accuracy of the original method. Our findings provide valuable insights, laying the groundwork for future research of fully discrete RK solvers for time-dependent PDEs involving low-order perturbed spatial operators. The ideas and techniques presented in this paper may extend beyond DG spatial discretization, encompassing other methods such as the continuous finite element method and the spectral method.

Below, we summarize our major theoretical novelties and contributions.
\begin{itemize}
	\item This paper exemplifies a novel framework for conducting stability and error analysis for numerical PDE schemes beyond the method of lines.
	\item In the stability analysis, we establish energy estimates bounding the $L^2$ norm of the high-order DG modes by jump terms. See Lemmas \ref{lem:PiL}, \ref{lem:nmLPL}, and \ref{lem:Luiwv}. These estimates indicate that filtering out high-order DG modes may have effects comparable to imposing jump penalties. The modal filter and jump penalty are two very commonly used approaches for stabilizing DG schemes. Our estimates reveal the inherent connections between the two stabilization approaches and demonstrate the stability effects in the fully discrete setting.
	\item In the error estimates, by adopting the compact strong form in the analysis, we gain further insights into the error accumulation through the RK iteration.
	\begin{enumerate}
		\item[$\circ$] Firstly, through the error equation \eqref{eq:err-xi-g} and its norm estimates \eqref{eq:xi-g}, one can directly see how the norm increment or decrement of the RK scheme can affect the growth of the projected error $\xi_h^n$, which bridges the error estimates of the RKDG scheme with recent works on the norm estimates of the RK method for homogeneous linear semi-negative problems \cite{tadmor2002semidiscrete,sun2017rk4,ranocha2018l_2,sun2019strong,sun2022energy}.
		\item[$\circ$] Secondly, by writing the error equation in the compact strong form, one can conveniently identify the remainder terms after a complete time step. It facilitates the design of novel approximation operators with information both in space and time for optimal error estimates.
	\end{enumerate}
\end{itemize}

\subsection{Organization of the paper.} As we are performing a systematic and unified analysis for sdA-RKDG schemes of arbitrarily high order, the paper unavoidably becomes lengthy and technical. Below, we explain the organization and provide tips on reading the paper.
\begin{itemize}
	\item Since there are no significant differences in the 1D and 2D analyses, we use unified notations for both cases and conduct the analysis simultaneously. However, to avoid being distracted by the complications in the 2D notations, it would be helpful for readers to concentrate on the 1D case when verifying propositions and lemmas.
	\item While the scope of the paper focuses on generic high-order sdA-RKDG schemes, to enhance readability, we begin each section on stability and error analysis with the second-order case \eqref{eq:sdArkdg2}. The lemmas introduced in the analysis of the second-order scheme will serve as the base case for extension to higher-order schemes. To avoid being confused by the general high-order cases, we suggest the readers to skip Sections \ref{sec:stab-g-thm}, \ref{sec:err-g}, \ref{sec:lem:Luiwv}, and \ref{newproj-general} in the first-round reading.
	\item The key to the stability analysis is to acquire estimates in Lemmas \ref{lem:PiL} and \ref{lem:Luiwv}, and the key to the optimal error estimates is to construct the approximation operators in Lemmas \ref{lem:newproj} and \ref{lem:newproj-general}. The proofs of these lemmas are technical. To avoid disrupting the flow of the paper, we postpone these proofs to Section \ref{sec:proofs} near the end of the paper. 
\end{itemize}

The rest of the paper is organized as follows. In Section \ref{sec:not}, we introduce the notations and preliminaries. Sections \ref{sec:stab} and \ref{sec:err} provide stability analysis and optimal error estimates, respectively. In each of these sections, we begin with the second-order case before presenting results for higher-order cases. Section \ref{sec:num} presents numerical examples to verify our theoretical analysis, while Section \ref{sec:proofs} contains proofs of several technical lemmas. Finally, conclusions are given in Section \ref{sec:conclusion}.

\section{Notations and preliminaries}\label{sec:not}

\subsection{DG spatial discretization}


\subsubsection{1D case} 

\;\smallskip

Consider the 1D linear advection equation
\begin{equation}
	u_t + \beta u_x = 0, \quad u = u(x,t), \quad x\in \Omega \subset \mathbb{R}, \quad 0<t\leq T.
\end{equation}
We assume the periodic boundary condition, and $\beta > 0$ is a positive constant. 

Let $\Omega_h = \{I_i\}_{i = 1}^N$, with $I_i = (x_{i-1/2}, x_{i+1/2})$ be a quasi-uniform mesh partition of $\Omega$. We denote by $h_i = x_{i+1/2} - x_{i-1/2}$ and $h = \max_{i = 1,\cdots, N} h_i$. The finite element space of the DG method is chosen as 
\begin{equation}
	{V}_h^\ell = \{v\in L^2(\Omega): v|_{I_i} \in {P}^\ell(I_i)\quad \forall i = 1,\cdots,N\},
\end{equation}
where ${P}^\ell(I_i)$ is the linear space for polynomials on $I_i$ of degree less than or equal to $\ell$. 

Let $\ip{w}{v} = \sum_{i=1}^N\ipi{w}{v} = \sum_{i=1}^N \int_{I_i} w v \dd x$ be the $L^2$ inner product and $\nm{v} = \sqrt{\ip{v}{v}}$ be the $L^2$ norm on $\Omega_h$. We denote by $L_h$ the DG approximation of the operator $L = -\beta \partial_x$ with the upwind flux, where $L_h: {V}_h^k \to {V}_h^k$ is defined such that
\begin{equation}\label{eq:Lh1d}
	\ip{L_h w_h}{v_h} = \beta\ip{w_h}{(v_h)_x} -\beta \sum_{i=1}^N \left((w_h^- v_h^-)_{i+\hf} - (w_h^- v_h^+)_{i-\hf}\right) \qquad \forall v_h \in {V}_h^k. 
\end{equation}
Here $(v_h^\pm)_{i+1/2} = \lim_{\varepsilon\to 0^\pm} v_h(x_{i+1/2}+\varepsilon)$ is the limit of $v_h$ from the left or the right of $x_{i+1/2}$.

Moreover, we introduce the following notations associated with jumps and traces.
\begin{align}
	[v_h] = v_h^+- v_h^-, \quad \jpip{w_h}{v_h} = \beta \sum_{i=1}^N ([w_h][v_h])_{i+1/2}, \quad \jpnm{v_h}= \jpip{v_h}{v_h}^{1/2}
\end{align} 
and 
\begin{equation}
	\nmG{v_h} = \left(\beta \sum_{i=1}^N \left((v_h^-)_{i+\hf}^2+(v_h^+)_{i+\hf}^2\right)\right)^{1/2}.
\end{equation}
Here the subscript $\Gamma$ corresponds to the mesh skeleton $\Gamma_h = \{\partial I_i:i = 1,\cdots, N\}$.

\subsubsection{2D case}

\;\smallskip

Consider the 2D linear advection equation
\begin{equation}
	u_t + \beta^x u_x + \beta^y u_y = 0, \quad u = u(x,y,t), \quad (x,y) \in \Omega \subset \mathbb{R}^2, \quad 0<t\leq T.
\end{equation}
To avoid unnecessary technicalities, we consider a rectangular domain with periodic boundary conditions.  $\beta^x>0$ and $\beta^y>0$ are positive constants. 

Let $\Omega_h = \{K_{ij}\}$ be a partition of the domain, where $K_{ij} = (x_{i-1/2},x_{i+1/2})\times (y_{j-1/2},y_{j+1/2})$, $i = 1,\cdots, N_x$ and $j = 1,\cdots, N_y$. We denote by $h_x^i = x_{i+1/2} - x_{i-1/2}$ and $h_y^j = y_{j+1/2} - y_{j-1/2}$.  Then we set
\begin{equation}
{V}_h^\ell = \{v\in L^2(\Omega): v|_{K_{ij}} \in {P}^\ell(K_{ij})\quad \forall i,j\}
\end{equation}
where ${P}^\ell(K_{ij})$ is the linear space for polynomials on $K_{ij}$ of degree less than or equal to $\ell$. 

Let $\ip{w}{v} = \sum_{i=1}^{N_x}\sum_{j=1}^{N_y}\ip{w}{v}_{K_{ij}} = \sum_{i=1}^{N_x}\sum_{j=1}^{N_y} \int_{K_{ij}} w v \dd x \dd y$ be the $L^2$ inner product and $\nm{v} = \sqrt{\ip{v}{v}}$ be the $L^2$ norm on $\Omega_h$. We denote by $L_h$ the DG approximation of the operator $L = -\beta^x \partial_x-\beta^y \partial_y$, where $L_h: V_h^k \to V_h^k$ is defined such that for all $v_h \in V_h^k$, 
\begin{equation}\label{eq:Lh2d}
\begin{aligned}
&\ip{L_h w_h}{v_h}  = \ip{w_h}{\beta^x(v_h)_x+\beta^y(v_h)_y}\\
&-\sum_{i=1}^{N_x}\sum_{j=1}^{N_y}\int_{x_{i-\hf}}^{x_{i+\hf}} \beta^y \left(w_h\left(x,y_{j+\hf}^-\right) v_h\left(x,y_{j+\hf}^-\right)-w_h\left(x,y_{j-\hf}^-\right) v_h\left(x,y_{j-\hf}^+\right)\right)\dd x\\
&-\sum_{i=1}^{N_x}\sum_{j=1}^{N_y} \int_{y_{j-\hf}}^{y_{j+\hf}}  \beta^x \left(w_h\left(x_{i+\hf}^-,y\right) v_h\left(x_{i+\hf}^-,y\right)-w_h\left(x_{i-\hf}^-,y\right) v_h\left(x_{i-\hf}^+,y\right)\right)\dd y.
\end{aligned}
\end{equation}
Here the superscript $-$ ($+$) indicates the limit from the left or below (right or above). 

Moreover, we introduce the following notations associated with jumps and traces.
\begin{align}
	[v_h]_{i,j+\hf} =& v_h\left(x,y_{j+\hf}^+\right) -v_h\left(x,y_{j+\hf}^-\right), \\
	[v_h]_{i+\hf,j} =& v_h\left(x_{i+\hf}^+,y\right) -v_h\left(x_{i+\hf}^-,y\right),\\
	\jpip{w_h}{v_h}
	=& \sum_{i=1}^{N_x}\sum_{j=1}^{N_y}\left(\int_{x_{i-\hf}}^{x_{i+\hf}} \beta^y\left([w_h][v_h]\right)_{i,j+\hf} \dd x + \int_{y_{j-\hf}}^{y_{j+\hf}}  \beta^x 
	\left([w_h][v_h]\right)_{i+\hf,j}
	\dd y\right),\\
	\jpnm{v_h} =& \jpip{v_h}{v_h}^{\frac{1}{2}},
	\end{align}
\begin{equation}
	\begin{aligned}
\nmG{v_h}
=& \left(\sum_{i=1}^{N_x}\sum_{j=1}^{N_y}\beta^y\int_{x_{i-\hf}}^{x_{i+\hf}} \left(v_h\left(x,y_{j+\hf}^-\right)\right)^2 + \left(v_h\left(x,y_{j+\hf}^+\right)\right)^2\dd x \right. \\
&\left.+ \sum_{i=1}^{N_x}\sum_{j=1}^{N_y} \beta^x\int_{y_{j-\hf}}^{y_{j+\hf}}   \left(v_h\left(x_{i+\hf}^-,y\right)\right)^2 + \left(v_h\left(x_{i+\hf}^+,y\right)\right)^2
\dd y\right)^\hf.
\end{aligned}
\end{equation}
Here the subscript $\Gamma$ corresponds to the mesh skeleton $\Gamma_h = \{\partial K_{ij}: i = 1,\cdots, N_x,j = 1,\cdots, N_y\}$.
\;

\subsubsection{Notations for both 1D and 2D cases.}

\;\smallskip

Throughout the paper, we use $C$, possibly with subscripts, for a generic constant independent of the time step size $\tau$ and the spatial mesh size $h$. 

Consider a mesh partition of the domain $\Omega_h = \{K\}$. The DG space is set as
\begin{equation}
	V_h^\ell = \{v\in L^2(\Omega): v|_K \in P^{\ell}(K)\quad \forall K\in \Omega_h\}.
\end{equation}
The $P^k$-DG operator $L_h: V_h^k \to V_h^k$ is defined in \eqref{eq:Lh1d} and \eqref{eq:Lh2d}. Note that $L_h$ is a semi-negative operator in the sense that $\ip{(L_h+L_h^\top)v_h}{v_h}\leq 0$ for all $v_h \in V_h^k$. This fact is implied by \eqref{eq:jp-defL2} in Proposition \ref{prop:sn} below. For details, we refer to \cite{zhang2010stability,cheng2017application,xu2019L2}.

\begin{PROP}[Negative semi-definiteness of the DG operator]\label{prop:sn} In both 1D and 2D,
\begin{subequations}
	\begin{align}
	\jpip{w_h}{v_h} = -\ip{(L_h + L_h^\top)w_h}{v_h} \quad \forall w_h, v_h \in V_h^k\label{eq:jp-defL},\\		
	\jpnm{v_h}^2 = -\ip{(L_h+L_h^\top)v_h}{v_h}\geq 0 \quad \forall v_h \in V_h^k\label{eq:jp-defL2}.
\end{align}
\end{subequations}
\end{PROP}

The semi-discrete $P^k$-DG scheme can be written as a semi-negative ODE system
\begin{equation}\label{eq:ODE}
\frac{\dd}{\dd t} u_h = L_h u_h. 	
\end{equation}

In the sdA-RKDG scheme, we also need the $P^{k-1}$-DG operator, which is defined as $\LL_h: V_h^k\to V_h^{k-1}$ such that (please pay attention to the test function space)
\begin{equation}\label{eq:def-LL_h}
\ip{\LL_h w_h}{v_h} = \ip{L_h w_h}{v_h} \quad \forall v_h \in {V}_h^{k-1}. 
\end{equation}
Let $V_h^\perp$ be the orthogonal complement of $V_h^{k-1}$ in $V_h^k$ under the $L^2$ inner product. We introduce the operators
\begin{equation}
	\Pi: L^2(\Omega) \to V_h^k, \quad \Piz: L^2(\Omega) \to V_h^{k-1}, \quand \Pip: L^2(\Omega) \to V_h^\perp
\end{equation} for the $L^2$ projection to $V_h^k$, $V_h^{k-1}$, and $V_h^\perp$, respectively. It can be seen that the $P^{k-1}$-DG operator $\LL_h:V_h^k\to V_h^{k-1}$ satisfies 
\begin{equation}
	\LL_h = \Piz L_h = (I-\Pip)L_h.
\end{equation}
Since $V_h^{k-1}\subset V_h^k$, $\LL_h$ is also an operator on $V_h^k$, namely, $\LL_h : V_h^{k}\to V_h^{k}$. 

With the inverse estimate, it can be seen that
\begin{equation}\label{eq:inv1}
\jpnm{v_h}\leq C\nmG{v_h} \leq Ch^{-\hf}\nm{v_h}\quad \forall v\in V_h^k.
\end{equation}
Then using the definition of $L_h$, we have $\nm{L_h v_h}\leq Ch^{-1}\nm{v_h}$ for all $v_h \in V_h^k$, or equivalently, $\nm{L_h}\leq Ch^{-1}$. Let us denote by $\lambda =  \tau/h$. Hence we have $\tau \nm{\LL_h} = \tau \nm{\Piz L_h} \leq \tau \nm{L_h}\leq C\lambda$. Without loss of generality, we assume $\lambda \leq 1$ throughout the paper.

\subsection{RKDG and sdA-RKDG schemes}

\;\smallskip

The original RKDG scheme in Butcher form is written as
\begin{subequations}\label{eq:rkdg-Butcher}
\begin{align}
u_h^{(i)} =& u_h^n + \tau \sum_{j = 1}^{i-1} a_{ij}L_h u_h^{(j)},\qquad i = 1,2,\cdots, s,\label{eq:rkdg-Butcher1}\\
u_h^{n+1} =& u_h^n + \tau \sum_{i = 1}^{s} b_{i}L_h u_h^{(i)}.
\end{align}
\end{subequations}
We will only consider explicit RK schemes in this paper and assume $a_{ij} = 0$ if $j\geq i$. In this case, for an $r$th-order RK method, the RK operator for the linear autonomous problem \eqref{eq:ODE} should match its stability polynomial, and can hence be written in a compact form as
\begin{equation}\label{eq:RKDG}
u_h^{n+1} = R_h u_h^n :=\sum_{i = 0}^s \alpha_i (\tau L_h)^i u_h^n = \left(I + \tau L_h \sum_{i = 1}^s \alpha_i (\tau L_h)^{i-1}\right) u_h^n,
\end{equation}
where $\alpha_i = (i!)^{-1} \text{ for } i = 0,1,\cdots, r$.

In the sdA-RKDG scheme, we replace the $P^k$-DG operator at inner stages by the $P^{k-1}$-DG operator $\LL_h$. Then the scheme can be written as
\begin{subequations}\label{eq:sdA-rkdg-Butcher}
\begin{align}
u_h^{(i)} =& u_h^n + \tau \sum_{j = 1}^{i-1} a_{ij}\LL_h u_h^{(j)},\qquad i = 1,2,\cdots, s,\label{eq:sdA-rkdg-Butcher1}\\
u_h^{n+1} =& u_h^n + \tau \sum_{i = 1}^{s} b_{i}L_h u_h^{(i)}.
\end{align}
\end{subequations}
It yields the compact form
\begin{equation}\label{eq:sdA-RKDG}
u_h^{n+1} = \RRR_h u_h^n:=\left(I + \tau L_h \sum_{i = 1}^s \alpha_i (\tau \LL_h)^{i-1}\right) u_h^n,
\end{equation}
where $\alpha_i = (i!)^{-1} \text{ for } i = 0,1,\cdots, r$. In other words, we only use the $P^k$-DG operator $L_h$ on the very left, and we use the $P^{k-1}$-DG operator $\LL_h$ at all other places.

When there is no confusion, we will refer to \eqref{eq:RKDG} as an RK$r$DG$k$ scheme, and \eqref{eq:sdA-RKDG} as an sdA-RK$r$DG$k$ scheme. 

\begin{REM}[Efficiency]
Compared to the original RKDG method, the sdA-RKDG method requires fewer computational degrees of freedom for the inner RK stages, which reduces the number of floating-point operations needed in each time step. The specific efficiency gain depends on the problem's settings, implementation, and solver compilation. For fluid problems, DG schemes typically involve multiple quadrature evaluations per time step. The sdA-RKDG method can lower computational costs by requiring fewer quadrature evaluations. Its actual savings in CPU time will depend on the cost of nodal evaluation of the flux functions and the implementation of the quadrature rule. For kinetic problems, the convection part is usually linear and could be in high dimensions. In some of these applications, the \( L_h \) and \( \LL_h \) matrices are assembled offline, with matrix multiplications performed online during time-stepping. In such cases, the better sparsity of $\LL_h$ could enhance solver efficiency. For further discussions on efficiency, we refer to \cite[Remark 2.3]{chen2024runge}.
\end{REM}

\begin{REM}[RK method in Butcher form]
The analysis in this paper only considers the RK method in Butcher form and may not be directly applicable to the RK method in Shu--Osher form (also known as the strong-stability-preserving form \cite{gottlieb2001strong,gottlieb2011strong}), unless a transformation is performed to rewrite a Shu--Osher form RK method back into the Butcher form. This is because our analysis relies on the following crucial fact: The operator on the very left of \eqref{eq:sdA-RKDG} is the full DG operator $L_h$, and does not involve the reduced operator $\LL_h$. Altering the inner stages of an RKDG scheme in Shu--Osher form could result in having $\LL_h$ present in the final stage, thereby disrupting the structure we have discussed. Such a method could yield a Class B sd-RKDG method, as described in \cite{chen2024runge}, which maintains optimal convergence for problems without sonic points but could suffer from accuracy degeneracy when sonic points occur.

We note that the RK method in Butcher form is also widely used for numerical solutions for hyperbolic conservation laws. See, for example, \cite{pazner2017stage,ranocha2020relaxation,ern2022invariant}.
\end{REM}
\section{Stability analysis}\label{sec:stab}

\subsection{Stability analysis: the second-order case}

\subsubsection{Main result}

\;\smallskip

In this section, we prove the stability of the sdA-RK2DG1 scheme \eqref{eq:sdArkdg2}.
\begin{THM}[Stability of sdA-RK2DG1]\label{thm:main-stab2}
The sdA-RK2DG1 scheme \eqref{eq:sdArkdg2} is monotonically stable, namely,  $
	\nm{u_h^{n+1}} \leq 
	\nm{u_h^{n}}$ for $\lambda \ll 1$. 
\end{THM}

\subsubsection{Preliminaries: stability of the RK2DG1 scheme}\label{sec:prelim:RK2DG1}

\;\smallskip

Before proceeding to prove Theorem \ref{thm:main-stab2}, we first recall the stability analysis of the standard RK2DG1 scheme \eqref{eq:rkdg2}, which can be written in the compact form as
\begin{equation}\label{eq:rkdg2-compact}
	u_h^{n+1} = \left(I + \tau L_h + \frac{\tau^2}{2}L_h^2\right)u_h^n:= R_{h,2} u_h^n.
\end{equation} 

The energy estimate of the RK2 scheme can be found in many previous works. See, for example, \cite{zhang2010stability,sun2017stability,sun2019strong,xu2019L2}.
\begin{PROP}[Energy identity of RK2]\label{prop:stabR2}
	\begin{align}\label{eq:stabR2}
		\nm{R_{h,2} w_h}^2  = \nm{w_h}^2 + \frac{\tau^4}{4}\nm{L_h^2 w_h}^2 - \tau \jpnm{w_h}^2 - {\tau^2}\jpip{w_h}{L_h w_h} - \frac{\tau^3}{2}\jpnm{L_h w_h}^2.
	\end{align}
\end{PROP}
In general, the right hand side of \eqref{eq:stabR2} may not be bounded by $\nm{w_h}^2$ appropriately. Indeed, it is known that the RK$2$DG$k$ method is unstable under the usual CFL condition if $k\geq 2$ \cite{cockburn2001runge}. But when coupled with $P^1$ elements, the RK2DG1 method is stable due to the following norm estimate, whose proof can be found in \cite{sun2017stability,xu2019L2,xu2024note}. 
\begin{PROP}[Bounding high-order DG derivatives hy jumps]\label{prop:estnmjp}
	With $P^1$ elements, 
	\begin{equation}\label{eq:estnmjp}
		{\tau^4}\nm{L_h^2 w_h}^2 \leq C\left(\tau \lambda^3\jpnm{w_h}^2 + \tau^3\lambda \jpnm{L_h w_h}^2\right).
	\end{equation}
\end{PROP}
Simple algebra shows that $\tau^2\jpip{w_h}{L_hw_h} \leq 
3/4\cdot \tau\jpnm{w_h}^2 + 1/3\cdot \tau^3\jpnm{L_hw_h}^2$, and hence this cross term can be absorbed by the negative jump terms in \eqref{eq:stabR2}. Thus, Propositions \ref{prop:stabR2} and \ref{prop:estnmjp} together imply the following stability results for the RK2DG1 scheme \eqref{eq:rkdg2-compact}. 
\begin{PROP}[Monotonicity stability of RK2DG1]\label{prop:stabR2+} Suppose $\lambda \ll 1$. Then there exists a positive constant $\varepsilon>0$ such that  
	\begin{align}
		\nm{R_{h,2} w_h}^2  \leq  \nm{w_h}^2 - \varepsilon \left(\tau \jpnm{w_h}^2 + {\tau^3}\jpnm{L_h w_h}^2\right).
	\end{align}
\end{PROP}

\subsubsection{Proof of Theorem \ref{thm:main-stab2}}\label{sec:proof-RK2DG1}

\;\smallskip

Now we prove Theorem \ref{thm:main-stab2}. The main idea for the proof is to write the sdA-RK2DG1 scheme \eqref{eq:sdArkdg2} as a perturbation of the standard RK2DG1 scheme \eqref{eq:rkdg2}. Recall that $\LL_h = (I-\Pip)L_h$. Hence we have
\begin{align}\label{eq:sds21line}
	u_h^{n+1} =& u_h^n + \tau L_hu_h^{n} + \frac{\tau^2}{2}L_h\LL_h u_h^n = R_{h,2} u_h^n - \frac{\tau^2}{2} L_h \Pip L_h u_h^n.
\end{align}
Note that each term is in $V_h^k$. We can take the square on both sides of \eqref{eq:sds21line} to get
\begin{align}\label{eq:stab2exp}
	\nm{u_h^{n+1}}^2 =&\nm{R_{h,2} u_h^n}^2 + \frac{\tau^4}{4} \nm{L_h \Pip L_h u_h^n}^2 - \tau^2\ip{L_h\Pip L_h u_h^n}{R_{h,2}u_h^n}.
\end{align}
The estimate of $\nm{R_{h,2} w_h}^2$ has already been shown in Proposition \ref{prop:stabR2+}. It remains to estimate $\nm{L_h \Pip L_h w_h}^2$ and terms of the form $\ip{L_h\Pip L_h w_h}{v_h}$. 

To estimate such terms, we need to introduce Lemma \ref{lem:PiL}, which basically states that the highest-order polynomial modes produced by the DG operator can be bounded by the jump term. Its proof is postponed to Section \ref{sec:proof-lem:PiL}. 
\begin{LEM}[Bounding highest-order DG modes by jumps] \label{lem:PiL}	With $P^k$ elements,
	\begin{subequations}
		\begin{align}
		\ip{\Pip L_h w_h}{v_h} \leq C h^{-\hf}\jpnm{w_h}\nm{v_h},\label{eq:key1}\\
		\ip{\Pip L_h^\top w_h}{v_h} \leq C  h^{-\hf}\jpnm{w_h}\nm{v_h}.\label{eq:key2}
	\end{align}
	\end{subequations}
\end{LEM}
With Lemma \ref{lem:PiL}, we can estimate $|\ip{L_h\Pip L_h w_h}{v_h}|$ and $\nm{L_h \Pip L_h w_h}$ in Lemmas \ref{lem:LPL} and \ref{lem:nmLPL}, respectively. 
\begin{LEM}\label{lem:LPL}With $P^k$ elements,
	\begin{align}
		|\ip{L_h\Pip L_h w_h}{v_h}| \leq Ch^{-1}\jpnm{w_h}\jpnm{v_h}.\label{eq:key}
	\end{align}
\end{LEM}
\begin{proof}
	Using the definition of the adjoint operator and the property of the $L^2$ projection, one can get
	\begin{equation}\label{eq:proof:LPL}
		\begin{aligned}
			\ip{L_h\Pip L_h w_h}{v_h} =& 	\ip{\Pip L_h w_h}{L_h^\top v_h}  = 	\ip{\Pip L_h w_h}{\Pip L_h^\top v_h}.
		\end{aligned}
	\end{equation}
	From Lemma \ref{lem:PiL}, we have $\nm{\Pip L_h w_h} \leq C h^{-\hf} \jpnm{w_h}$ and $\nm{\Pip L_h^\top v_h} \leq C h^{-\hf} \jpnm{v_h}$. Therefore, applying the Cauchy--Schwarz inequality to the right side of \eqref{eq:proof:LPL}, we have 
	\begin{equation}
	\begin{aligned}
		|\ip{L_h\Pip L_h w_h}{v_h}|
		\leq Ch^{-\hf}\jpnm{w_h}\cdot Ch^{-\hf}\jpnm{v_h} \leq  Ch^{-1}\jpnm{w_h}\jpnm{v_h}.
	\end{aligned}
	\end{equation}
\end{proof}
\begin{LEM}\label{lem:nmLPL}With $P^k$ elements, 
	\begin{equation}\label{eq:nmLPL}
		\nm{L_h\Pip L_h w_h} \leq  Ch^{-\frac{3}{2}}\jpnm{w_h}.
	\end{equation}
\end{LEM}
\begin{proof}
	Note that \eqref{eq:key1} of Lemma \ref{lem:PiL} implies
	$\nm{\Pip L_h w_h} \leq C h^{-\hf}\jpnm{w_h}$. 	Recall that $\nm{L_h v_h}\leq Ch^{-1} \nm{v_h}$. Thus $\nm{L_h\Pip L_h w_h} \leq Ch^{-1}\nm{\Pip L_h w_h}\leq C h^{-3/2}\jpnm{w_h}$.
\end{proof}

Now we continue to prove Theorem \ref{thm:main-stab2}. Applying Lemma \ref{lem:nmLPL}, we have
\begin{equation}
	\frac{\tau^4}{4}\nm{L_h\Pip L_h u_h^n}^2 \leq C \lambda^3 \tau \jpnm{u_h^n}^2.\label{eq:st1}
\end{equation}
Using Lemma \ref{lem:LPL}, the triangle inequality, and Young's inequality $ab\leq C(a^2+b^2)$, we have
\begin{equation}
	\begin{aligned}
	\tau^2|\ip{L_h\Pip L_h u_h^n}{R_{h,2}u_h^n}| \leq & C\lambda \tau \jpnm{u_h^n} \jpnm{R_{h,2}u_h^n}\\
	\leq& C\lambda\tau \jpnm{u_h^n}\left( \jpnm{u_h^n} + \tau \jpnm{L_h u_h^n} + \frac{\tau^2}{2} \jpnm{L_h^2 u_h^n}\right)\\
	\leq& C\lambda \tau \left(\jpnm{u_h^n}^2 + \tau^2 \jpnm{L_hu_h^n}^2 + \tau^4 \jpnm{L_h^2u_h^n}^2\right).
	\label{eq:st2}
\end{aligned}
\end{equation}
Substituting \eqref{eq:st1} and \eqref{eq:st2} into \eqref{eq:stab2exp}, it gives
\begin{align}\label{eq:stab-rkdg2-jp}
	\nm{u_h^{n+1}}^2 = \nm{\RRR_{h,2} u_h^n}^2\leq \nm{R_{h,2} u_h^n}^2 +C\lambda\left(\tau \jpnm{u_h^n}^2 +\tau^3 \jpnm{L_h u_h^n}^2 + \tau^5 \jpnm{L_h^2 u_h^n}^2\right).
\end{align}
	Recall \eqref{eq:inv1}. We have $\tau \jpnm{v_h}^2 \leq C\lambda \nm{v_h}^2$. Thus, with Proposition \ref{prop:estnmjp}, we get
\begin{equation}\label{eq:estL2u}
	\tau^5 \jpnm{L_h^2u_h^n}^2 \leq C \lambda \tau^4  \nm{L_h^2 u_h^n}^2 \leq C \left(\lambda^4 \tau \jpnm{u_h^n}^2 + \lambda^2 \tau^3 \jpnm{L_h u_h^n}^2\right).
\end{equation}
Substituting \eqref{eq:estL2u} into \eqref{eq:stab-rkdg2-jp}, one can get \begin{equation}
	\nm{u_h^{n+1}}^2 \leq \nm{R_{h,2} u_h^n}^2 +C\lambda\left(\tau \jpnm{u_h^n}^2 +C\lambda\tau^3 \jpnm{L_h u_h^n}^2\right).
\end{equation} 
Applying the stability estimate of $\nm{R_{h,2}u_h^n}^2$ in Proposition \ref{prop:stabR2+}, it yields
\begin{equation}
	\nm{u_h^{n+1}}^2 \leq \nm{u_h^n}^2 +(C\lambda-\varepsilon)\tau \jpnm{u_h^n}^2 +(C\lambda-\varepsilon)\tau^3 \jpnm{L_h u_h^n}^2.
\end{equation}
Recall that $\varepsilon$ is a fixed positive constant. Therefore, when $\lambda\ll 1$, we have $C\lambda-\varepsilon<0$, which implies $\nm{u_h^{n+1}}^2 \leq \nm{u_h^{n}}^2$. This completes the proof of Theorem \ref{thm:main-stab2}.

\subsection{Stability analysis: general cases}

\subsubsection{Main results}

\;\smallskip

The main result of our energy estimate is given in Lemma \ref{lem:keykey}, which states that the $L^2$ energy of the sdA-RKDG solution can be bounded by that of the original RKDG solution along with some diminishing jump terms as $\lambda$ approaches zero. The proof of Lemma \ref{lem:keykey} is given in Section \ref{sec:lem:keykey}.
\begin{LEM}[Generalization of \eqref{eq:stab-rkdg2-jp}]\label{lem:keykey} For $\RRR_h$ defined in \eqref{eq:sdA-RKDG}, and more generally, for $\RRR_h$ defined in \eqref{eq:sds-general}, we have
	$\nm{\RRR_h w_h}^2 \leq \nm{R_h w_h}^2 + C\lambda \sum_{i = 0}^{2s-2}\tau^{2i+1}\jpnm{L_h^i w_h}^2$. 
\end{LEM}

With Lemma \ref{lem:keykey}, we can prove that the sdA-RKDG schemes share the same stability property as the original RKDG schemes, as explained in Theorem \ref{thm:key}. A sketch proof of the theorem is given in Section \ref{sec:stab-g-thm}.
\begin{THM}[Stability of sdA-RKDG schemes]\label{thm:key}
	Stability results in Theorem \ref{thm:stabr} and Theorem \ref{thm:stabr-low} remain valid after replacing $R_h= \sum_{i = 0}^s \alpha_i (\tau L_h)^{i}$ by $\RRR_h = I + \tau L_h \sum_{i = 1}^s \alpha_i (\tau \LL_h)^{i-1}$.
\end{THM}

\subsubsection{Preliminaries: stability of generic RKDG schemes}

\;\smallskip

The RK schemes admit the following energy expansion \cite{sun2019strong,sun2022energy}. 
\begin{LEM}[Energy identity of RK]\label{lem:energy}
	Let $R_h w_h  = \sum_{i = 0}^s\alpha_i (\tau L_h)^i w_h$. Then
	\begin{equation}\label{eq:rk-energy}
		\nm{R_h w_h}^2 = \sum_{i=0}^s \beta_i \tau^{2i}\nm{L_h^i w_h}^2 +  \sum_{i,j = 0}^{s-1}\gamma_{ij} \tau^{i+j+1}\jpip{L_h^i w_h}{L_h^jw_h}.
	\end{equation}
	Here 
	\begin{align}
	\beta_i = \sum_{\ell = \max\{0,2i-s\}}^{\min\{2i,s\}}\alpha_{\ell}\alpha_{2i-\ell}(-1)^{i-\ell}\quand
	\gamma_{ij} = \sum_{\ell=\max\{0,i+j+1-s\}}^{\min\{i,j\}}(-1)^{\min\{i,j\}+1-\ell}\alpha_{\ell}\alpha_{i+j+1-\ell}.
	\end{align}
\end{LEM}

In general, the stability of the RK scheme is determined by the sign of the first nonzero $\beta_i$ and the size of the largest negative-definite leading principal minor of $(\gamma_{ij})$. The general stability criteria with arbitrary coefficients $\alpha_i$ can be found in \cite{sun2019strong,xu2020superconvergence,sun2022energy}. Below, we list the stability of $r$-stage $r$-th order RK scheme $R_h = \sum_{i = 0}^r(i!)^{-1} (\tau L_h)^i$. 

\begin{DEFN}[Different types of stability]\;
	\begin{enumerate}
		\item The method is weakly$(\gamma)$ stable if 
		$\nm{R_h w_h} \leq (1+ C\lambda^{\gamma})\nm{w_h}$. This implies $\nm{u_h^n}\leq e^{Ct^n}\nm{u_h^0}$ if $\lambda \leq \lambda_0 h^{1/(\gamma-1)}$ for some constant $\lambda_0$.  
		\item The method is strongly($m_\star$) stable if $\nm{R_h^{m} w_h} \leq \nm{w_h}$ for $m\geq m_\star$. This implies $\nm{u_h^{n+m}}\leq \nm{u_h^n}$ if $\lambda \leq \lambda_0$ for some constant $\lambda_0$.  
		\item The method is monotonically stable if it is strongly(1) stable. 
	\end{enumerate}
\end{DEFN}
\begin{THM}[Stability of the $r$-stage RK$r$DG$k$]\label{thm:stabr}
	Suppose $L_h+L_h^\top \leq 0$ is negative semi-definite. Let $R_h = \sum_{i = 0}^r(i!)^{-1} (\tau L_h)^i$. For all $k\geq 0$, we have the following stability results.
	\begin{enumerate}
		\item If $r \equiv 1 \Mod{4}$, then the method is weakly($r+1$) stable.  
		\item If $r \equiv 2 \Mod{4}$, then the method is weakly($r+2$) stable.
		\item If $r \equiv 3 \Mod{4}$, then the method is monotonically stable.
		\item If $r \equiv 0 \Mod{4}$, then the method is strongly($2$) stable.
	\end{enumerate}
\end{THM}
Additionally, similar to Proposition \ref{prop:estnmjp} for the RK2DG1 case, by noting that $\nm{L_h^i w_h}$ can be bounded by jump terms $\jpnm{L_h^jw_h}$ when the polynomial degree $k$ is small compared to $i$, the following improved stability results can be acquired \cite{xu2019L2,xu2020superconvergence}. 
\begin{THM}[Stability of RK$r$DG$k$ with low-order polynomials]\label{thm:stabr-low}
	Consider the RK$r$DG$k$ method $R_h=\sum_{i = 1}^r(i!)^{-1} (\tau L_h)^i$. We have the following stability results.
	\begin{enumerate}
		\item The method is monotonically stable if $k\leq \lfloor r/2\rfloor - 1$. 
		\item The method is strongly($m_\star$) stable for some $m_\star\geq 1$ if $k\leq \lfloor r/2\rfloor$. If $m_\star = 1$, then the method is monotonically stable.
	\end{enumerate} 
\end{THM}

\begin{REM}[``Strong stability"]
	In \cite{sun2019strong,sun2022energy}, we referred to strong$(1)$ stability or monotonicity stability as ``strong stability," while in \cite{xu2019L2,xu2020superconvergence}, the authors referred to the stability implied by strong$(m_\star)$ stability as the ``strong stability." 	Here we avoid using the term ``{strong stability}" to prevent ambiguity, as it has been referred to differently in previous works. 
\end{REM}

\begin{REM}[CFL numbers]
Although the sdA-RKDG method preserves the same type of stability as the original RKDG method, the corresponding CFL numbers may differ slightly. We used the Fourier approach as described in \cite{chen2024runge} to compute the CFL numbers of the sdA-RK$r$DG$k$ schemes with $r = k+1$ up to eighth order. The results are documented in Table \ref{tab:CFL}. For the cases examined, the CFL numbers of the sdA-RKDG schemes are slightly smaller than those of the corresponding RKDG schemes. However, the difference does not appear to be significant, nor does it result in a substantial reduction in applicable time step sizes, at least for these commonly used cases. We also want to remark that, despite the slightly smaller CFL numbers of the sdA-RKDG schemes in Table \ref{tab:CFL}, several Class B sd-RKDG schemes in \cite{chen2024runge} achieve larger CFL numbers compared to the corresponding RKDG schemes.

	\end{REM}
	\begin{table}[h!]
		\centering
		\begin{tabular}{cccccccc}
			\hline
			$r$&$2$&$3$&$4$&$5$&$6$&$7$&$8$\\
			\hline
			sdA-RKDG&$0.333$&$0.191$&$0.127$&$0.104$&$0.085$&$0.076$&$0.064$\\
			RKDG&$0.333$&$0.209$&$0.145$&$0.115$&$0.093$&$0.080$&$0.070$\\
			\hline
			Ratio& 100\%&91.4\%&87.6\%&90.4\%&91.4\%&95.0\%&91.4\%\\
			\hline
		\end{tabular}
		\caption{CFL numbers of the sdA-RK$r$DG$k$ schemes with $r = k+1$.}\label{tab:CFL}
	\end{table}  

\subsubsection{Proof of Lemma \ref{lem:keykey}}\label{sec:lem:keykey}

\;\smallskip

To explain the motivation, we start by taking a glance at the third-order sdA-RKDG scheme $u_h^{n+1} = \RRR_{h,3} u_h^n$, where
\begin{equation}\label{eq:RRRh3}
	\RRR_{h,3}w_h = w_h + \tau L_h w_h + \frac{\tau^2}{2}L_h\LL_h w_h + \frac{\tau^3}{6}L_h\LL_h^2w_h.
\end{equation}
With $\LL_h = (I-\Pip)L_h$, \eqref{eq:RRRh3} can be reformulated as 
\begin{equation}\label{eq:sdA-RKDG3}
	\begin{aligned}
		\RRR_{h,3}w_h =& \left(I + \tau L_h + \frac{\tau^2}{2}L_h^2 + \frac{\tau^3}{6}L_h^3\right)w_h -\frac{\tau^2}{2} L_h \Pip L_h w_h - \frac{\tau^3}{6}L_h\left(L_h^2-\LL_h^2\right)w_h\\
		=&R_{h,3} w_h -\frac{\tau^2}{2} L_h \Pip L_h w_h - \frac{\tau^3}{6}\left(L_h^2\Pip L_h +L_h\Pip L_h^2 -L_h(\Pip L_h)^2\right)w_h.
	\end{aligned}
\end{equation}
From \eqref{eq:sdA-RKDG3}, it can be seen that in general, we need to deal with operators in the form of a noncommunicative product of $L_h^i$ and $\Pip$. To this end, we need to introduce some notations. 
We denote by $\ui = (i_1,i_2,\cdots, i_{[\ui]})$ a vector with positive integer entries. Let $[\ui]\geq2$ be the length of the vector $\ui$ and $|\ui| = i_1 + i_2 + \cdots + i_{[\ui]}$ be the $\ell^1$ norm of the vector $\ui$. Define 
\begin{equation}
	L_h^{\ui} = L_h^{i_1}\Pip L_h^{i_2}\Pip\cdots\Pip L_h^{i_{[\ui]}}.
\end{equation}

\underline{\em Reformulation of the sdA-RKDG scheme.}  Subtracting \eqref{eq:RKDG} from \eqref{eq:sdA-RKDG} and replacing $u_h^n$ by $w_h$, we have
\begin{equation}\label{eq:RRR-R}
\RRR_h w_h = R_h w_h + \sum_{j=1}^s \alpha_j \tau^{j} L_h\left(\LL_h^{j-1}-L_h^{j-1}\right) w_h.
\end{equation}
We can expand $L_h(\LL_h^{i-1} - L_h^{i-1})$ to rewrite the sdA-RKDG scheme in the following form. 
\begin{equation}\label{eq:sds-general}
	\RRR_h w_h = R_h w_h + \sum_{\ui: |\ui|\leq s} \mu_{\ui}\tau^{|\ui|} L_h^{\ui} w_h.
\end{equation}
Here the real numbers $\{\mu_{\ui}\}$ are known and depend only on $\{\alpha_i\}$. As their specific values are not essential, we omit the characterizations of the coefficients $\mu_{\ui}$. \smallskip

\underline{\em Energy estimates.} By squaring both sides of \eqref{eq:sds-general} and using the elementary inequality $(\sum_{i}|a_i|)^2\leq C\sum_i|a_i|^2$, one can see that the sdA-RKDG scheme admits the following energy estimate. 
\begin{equation}\label{eq:stabg}
\begin{aligned}
	\nm{\RRR_hw_h}^2 =& \nm{R_h w_h}^2 +  \nm{\sum_{\ui: |\ui|\leq s}\mu_{\ui} \tau^{|\ui|}L_h^{\ui} w_h}^2 + 2\sum_{\ui: |\ui|\leq s} \mu_{\ui}\tau^{|\ui|} \ip{L_h^{\ui} w_h}{R_h w_h}\\
	\leq& \nm{R_h w_h}^2 +  C\sum_{\ui: |\ui|\leq s} \tau^{2|\ui|}\nm{L_h^{\ui} w_h}^2 + C\sum_{\ui: |\ui|\leq s} \tau^{|\ui|} \left|\ip{L_h^{\ui} w_h}{R_h w_h}\right|.	
\end{aligned}
\end{equation}

Similar to the second-order case, we need to estimate terms of the forms $\nm{L_h^{\ui} w_h}^2$ and $\ip{L_h^{\ui} w_h}{v_h}$. To this end, we need the following Lemmas \ref{lem:Luiwv} and \ref{lem:Luiw}.  The proof of Lemma \ref{lem:Luiwv} is technical and is  therefore postponed to Section \ref{sec:lem:Luiwv}. 
\begin{LEM}[Generalization of Lemma \ref{lem:LPL}]\label{lem:Luiwv}
	\begin{equation}\label{eq:Luiwv}
		\begin{aligned}
			\left|\ip{L_h^{\ui} w_h}{v_h}\right|
			\leq&C\sum_{j=0}^{i_1-1}h^{-|\ui|+i_{[\ui]}+j}\jpnm{L_h^{i_{[\ui]}-1} w_h}\jpnm{L_h^{j} v_h} .
		\end{aligned}
	\end{equation}	
\end{LEM}

\begin{LEM}[Generalization of Lemma \ref{lem:nmLPL}]\label{lem:Luiw}
	\begin{equation}\label{eq:Luiw}
		\begin{aligned}
			\nm{L_h^{\ui} w_h}
			\leq&C\sum_{j=0}^{i_1-1}h^{-|\ui|+i_{[\ui]}-\hf}\jpnm{L_h^{i_{[\ui]}-1} w_h}.
		\end{aligned}
	\end{equation}	
\end{LEM}
\begin{proof}
We use the inverse estimate \eqref{eq:inv1} and the fact $\nm{L_h}\leq Ch^{-1}$ to get 
\begin{equation}
\jpnm{L_h^{j} v_h}\leq Ch^{-\hf}\nm{L_h^j v_h} \leq Ch^{-\hf}\nm{L_h}^j \nm{v_h} \leq Ch^{-j-\hf}\nm{v_h}.
\end{equation} 
Hence Lemma \ref{lem:Luiwv} implies
\begin{equation}\label{eq:Luiwv}
\begin{aligned}
\left|\ip{L_h^{\ui} w_h}{v_h}\right|
\leq&C\sum_{j=0}^{i_1-1}h^{-|\ui|+i_{[\ui]}-\hf}\jpnm{L_h^{i_{[\ui]}-1} w_h}\nm{v_h} .
\end{aligned}
\end{equation}	
The proof of Lemma \ref{lem:Luiw} can be completed after taking $v_h = L_h^{\ui} w_h$ in \eqref{eq:Luiwv} and dividing by $\nm{L_h^{\ui} w_h}$ on both sides. 
\end{proof}

Now we continue to prove Lemma \ref{lem:keykey}. To estimate the second term in \eqref{eq:stabg}, we apply Lemma \ref{lem:Luiw} to get
\begin{equation}\label{eq:stabg-norm}
\begin{aligned}
\tau^{2|\ui|}\nm{L_h^{\ui} w_h}^2
\leq&C\sum_{j=0}^{i_1-1}\lambda^{2|\ui|-2i_{[\ui]}+1} \tau^{2i_{[\ui]}-1}\jpnm{L_h^{i_{[\ui]}-1} w_h}^2\leq C\lambda\tau \sum_{j=0}^{i_1-1}\jpnm{(\tau L_h)^{i_{[\ui]}-1} w_h}^2.
\end{aligned}
\end{equation}	
To estimate the last term in \eqref{eq:stabg}, we set $v_h = R_h w_h$ in Lemma \ref{lem:Luiwv}. Then
\begin{equation}
	\begin{aligned}
		\left|\ip{\tau^{|\ui|} L_h^{\ui} w_h}{R_h w_h}\right|
		\leq&C\tau \sum_{j=0}^{i_1-1}\lambda^{|\ui|-i_{[\ui]}-j}\jpnm{(\tau L_h)^{i_{[\ui]}-1} w_h}\jpnm{(\tau L_h)^{j} R_hw_h}.
	\end{aligned}
\end{equation}	
With the triangle inequality, we have $\jpnm{(\tau L_h)^{j} R_hw_h} \leq \sum_{\ell=0}^{s}|\alpha_\ell| \jpnm{(\tau L_h)^{\ell+j} w_h}$.	Therefore, 
\begin{equation}
	\begin{aligned}
		\left|\ip{\tau^{|\ui|} L_h^{\ui} w_h}{R_h w_h}\right|
		\leq&C\tau \sum_{j=0}^{i_1-1}\sum_{\ell = 0}^s\lambda^{|\ui|-i_{[\ui]}-j}\jpnm{(\tau L_h)^{i_{[\ui]}-1} w_h}\jpnm{(\tau L_h)^{\ell + j} w_h}.
	\end{aligned}
\end{equation}	
Note that $|\ui|-i_{[\ui]}-j \geq i_1 - j \geq 1$. Hence $\lambda^{|\ui|-i_{[\ui]}-j}\leq \lambda$. Applying Young's inequality $ab\leq C(a^2+b^2)$ gives
\begin{equation}\label{eq:stabg-inner}
\left|\ip{\tau^{|\ui|} L_h^{\ui} w_h}{R_h w_h}\right| \leq C\lambda \tau \sum_{j=0}^{s+i_1-1}\jpnm{(\tau L_h)^{j} w_h}^2.
\end{equation}
Substitute \eqref{eq:stabg-norm} and \eqref{eq:stabg-inner} into 
\eqref{eq:stabg}. Recalling that $1\leq i_1, i_{[\ui]} \leq s-1$, one can get
\begin{equation}
	\nm{u_h^{n+1}}^2 = \nm{\RRR_h w_h}^2 \leq  \nm{R_h w_h}^2 + C\lambda \tau \sum_{j = 0}^{2s-2} \jpnm{(\tau L_h)^j w_h}^2.
\end{equation} 
This completes the proof of Lemma \ref{lem:keykey}. 

\begin{REM}[Stability of other Class A sd-RKDG ] In addition to \eqref{eq:sdA-rkdg-Butcher}, we have also considered other sd-RKDG schemes of Class A in \cite{chen2024runge}. In such schemes, each DG operator in \eqref{eq:sdA-rkdg-Butcher1} does not have to be the $P^{k-1}$-DG operator $\LL_h$ and can instead be chosen as either a $P^k$-DG operator $L_h$ or a $P^{k-1}$-DG operator $\LL_h$. It is important to note that such schemes can also be expressed in the generic form of \eqref{eq:sds-general}, and all previous analyses can proceed accordingly. Hence, Lemma \ref{lem:keykey} holds for all Class A sd-RKDG schemes in \cite{chen2024runge}. Consequently, with the analysis in Section \ref{sec:stab-g-thm}, one can show that all Class A sd-RKDG schemes will inherit the same stability property as their parent RKDG schemes.
\end{REM}

	\begin{REM}[$Q^k$-DG]
		Here, we have been focusing on the $P^k$-DG method. The extension to the $Q^k$-DG method may not be very straightforward. The key point in the stability analysis with $P^k$ elements is that the strong derivative on $V_h^k$ returns a function in the reduced space $V_h^{k-1}$, namely $\beta^x (w_h)_x + \beta^y (w_h)_y \in V_h^{k-1}$ for all $w_h \in V_h^k$. This fact is used in the proof of Lemma \ref{lem:PiL}, the key to other preparatory lemmas such as Lemmas \ref{lem:LPL} and \ref{lem:Luiwv}, to ensure that the first term in \eqref{eq:stabkeylem-2D} vanishes. However, for the $Q^k$-DG method, $\beta^x (w_h)_x + \beta^y (w_h)_y$ is not a piecewise $Q^{k-1}$ polynomial, and the analysis would not hold if we naively set $V_h^{k-1}$ to be the $Q^{k-1}$ space. This technical difficulty can be circumvented by replacing $V_h^{k-1}$ with a richer space $\widetilde{V}_h^{k} \supseteq \{\beta^x (w_h)_x + \beta^y (w_h)_y: w_h \in V_h^k\}$. For example, in the second-order case, one can use the $(P^1,Q^1)$ space pair, rather than the $(Q^0,Q^1)$ space pair, for the inner and final stages. Then, the stability can be proved along the same lines as in Section \ref{sec:proof-RK2DG1}, and the accuracy order on quasi-uniform meshes can be proved by following similar arguments in Section \ref{sec:proof-thm:err2} and replacing $\Pig$ with the 2D Gauss--Radau projection. In our future work, we will further investigate the $(Q^{0},Q^1)$ space pair, along with its high-order extension. We will explore whether stability can be established using alternative proof techniques or enforced through appropriate post-processing methods.
	\end{REM}

\subsubsection{Sketch proof of Theorem \ref{thm:key}}\label{sec:stab-g-thm}

\;\smallskip

	From Lemma \ref{lem:energy}, one can see that the stability of the RKDG and sdA-RKDG schemes is determined by the sign of the leading term of $\beta_i$ and the size of the largest negative-definite leading principal minor of $(\gamma_{ij})$. Note that when $\lambda$ is sufficiently small, the additional term $C\lambda \sum_{i=0}^{2s-2} \tau^{2i+1}\jpnm{L_h^i w_h}^2$ in Lemma \ref{lem:keykey} neither affects $\beta_i$ nor changes the negative definiteness of the leading principal minors in $(\gamma_{ij})$. Hence, all the stability results of the original RKDG methods will be retained. 
	
	Below, as an example, we provide the proof of Case 3 on monotonicity stability in Theorem \ref{thm:stabr}. We will also briefly comment on the proof of Case 4 regarding strong($2$) stability. The proof of other cases is omitted. 
	
	\underline{\em Proof of Case 3.} Let $\zeta$ denote the index such that $\beta_i = 0$ for $i < \zeta$ and $\beta_\zeta \neq 0$. According to \cite[Lemma 4.1]{sun2019strong}, for $r \equiv 3 \Mod{4}$, it follows that $\zeta = (r + 1)/2$, $\beta_\zeta = 2(-1)^{\zeta+1}/(r+1)! < 0$, and $\Gamma_\zeta = (r_{ij})_{0\leq i,j\leq \zeta} < 0$ is negative definite. By following the proof outlined in \cite[Theorem 2.7]{sun2019strong}, one can derive \cite[(2.21)]{sun2019strong}. After simplifying the coefficients there, it can be seen that there exists a positive constant $\varepsilon>0$ such that
	\begin{equation}\label{eq:RR_Rh}
		\nm{R_h w_h}^2 \leq \nm{w_h}^2 + (\beta_\zeta + C\lambda) \tau^{2\zeta}\nm{L_h^\zeta w_h}^2 - \varepsilon \sum_{i = 0}^{\zeta-1}\tau^{2i+1}\jpnm{L_h^i w_h}^2.
	\end{equation} 
	Note that for $i\geq \zeta$, we can use the inverse estimate $\jpnm{v_h}\leq Ch^{-1/2}\nm{v_h}$  to get 
	\begin{equation}
	\tau^{2i+1}\jpnm{L_h^i w_h}^2 \leq C\lambda \tau^{2i}\nm{L_h^i w_h}^2 \leq C\lambda^{2i-2\zeta+1}\tau^{2\zeta}\nm{L_h^\zeta w_h}^2 \leq C\lambda \tau^{2\zeta}\nm{L_h^\zeta w_h}^2.
	\end{equation}
	Therefore, 
	\begin{equation}\label{eq:RR_jp}
		C\lambda \sum_{i = 0}^{2s-2}\tau^{2i+1}\jpnm{L_h^i w_h}^2 \leq C\lambda \sum_{i = 0}^{\zeta-1}\tau^{2i+1}\jpnm{L_h^i w_h}^2 + C\lambda \tau^{2\zeta}\nm{L_h^\zeta w_h}^2.
	\end{equation}
	Then substituting estimates \eqref{eq:RR_Rh} and \eqref{eq:RR_jp} into Lemma \ref{lem:keykey}, we have 
	\begin{equation}
		\nm{\RRR_h w_h}^2 \leq \nm{w_h}^2 + (\beta_\zeta + C\lambda) \tau^{2\zeta}\nm{L_h^\zeta w_h}^2 - (\varepsilon-C\lambda) \sum_{i = 0}^{\zeta-1}\tau^{2i+1}\jpnm{L_h^i w_h}^2.
	\end{equation}
	Since $\beta_\zeta < 0$ and $\varepsilon>0$, we have $\nm{\RRR_h w_h}^2 \leq \nm{w_h}^2$ when $\lambda \ll 1$, which completes the proof of monotonicity stability of $\RRR_h$ for $r \equiv 3 \Mod{4}$. 
	
	\underline{\em Comments on Case 4.} The main catch in the proof of Case 4 is that one cannot express multiple steps of $R_h$, namely $R_h^{m}$, in the form of \eqref{eq:sdA-RKDG}. However, it still yields the form of \eqref{eq:sds-general}. Therefore, one can use Lemma \ref{lem:keykey} for \eqref{eq:sds-general} to combine multiple steps and analyze the strong$(2)$ stability. To be more specific, we consider the fourth-order RK method as an example. The corresponding sdA-RK4DG$k$ scheme takes the form
	\begin{equation}
		\RRR_h = I + \tau L_h \left(I+ \frac{1}{2}(\tau\LL_h) + \frac{1}{6} (\tau \LL_h)^2 + \frac{1}{24} (\tau \LL_h)^3\right).
	\end{equation}
	When computing $\RRR_h^2$, we encounter the term $(\tau L_h)^2$, which is not included in \eqref{eq:sdA-RKDG}. However, direct computation  shows that
	\begin{equation} \label{eq:RK4-2step}
		\begin{aligned}
		\RRR_h^2 = &I + 2\tau L_h + \frac{1}{2}(2\tau L_h)^2 + \frac{1}{6} (2\tau L_h)^3 + \frac{1}{24}(2\tau L_h)^4 \\
		&+ \frac{1}{4}(\tau L_h)^5 + \frac{5}{72}(\tau L_h)^6 + \frac{1}{72}(\tau L_h)^7 + \frac{1}{576}(\tau L_h)^8 + \sum_{\ui: |\ui|\leq 8} \mu_{\ui}\tau^{|\ui|} L_h^{\ui} \\
		=&R_h^2 + \sum_{\ui: |\ui|\leq 8} \mu_{\ui}\tau^{|\ui|} L_h^{\ui},
		\end{aligned}
	\end{equation}
	where $R_h$ corresponds to the original RK4DG$k$ scheme and the detailed characterization of $\{\mu_{\ui}\}$ is omitted. Note that \eqref{eq:RK4-2step} is still in the form of \eqref{eq:sds-general}. Hence one can apply Lemma \ref{lem:keykey} and use the norm estimate of $R_h^2$ in \cite{sun2017rk4} to show that $\nm{\RRR_h^2 w_h}\leq \nm{w_h}$ if $\lambda \leq \lambda_2$. Similarly, one can expand $\RRR_h^3$ to show that  $\nm{\RRR_h^3 w_h}\leq \nm{w_h}$  if $\lambda \leq \lambda_3$. See \cite{sun2019strong} for the $R_h^3$ counterpart. Combining both the two-step and the three-step stability, one can use the norm inequality to show that $\nm{\RRR_h^m w_h}\leq \nm{w_h}$ for all $m\geq 2$, if $\lambda \leq \min(\lambda_2,\lambda_3)$. This completes the proof of strong$(2)$ stability of the fourth-order RK method. The extension to the general case where $r \equiv 0 \Mod{4}$ follows similar arguments and is omitted. 
\section{Optimal error estimates}\label{sec:err}

\subsection{Preliminaries: projection operators} 

\subsubsection{1D case} 

\;\smallskip

We denote by $H^\ell(\Omega_h) = \{v\in L^2(\Omega): v|_{I_i}\in H^\ell(I_i)\  \forall i = 1,\cdots,N\}$ and $\nm{v}_{H^\ell} = \left(\sum_{i=1}^N \nm{v}_{H^\ell(I_i)}^2\right)^{1/2}$.
The following Gauss--Radau projection \cite{castillo2002optimal,shu2009discontinuous} is used for the optimal error estimate of the DG method in 1D. $\Pig: H^1(\Omega_h) \to V_h^k$ is an operator such that for $w\in H^1(\Omega_h)$, the projection error $\eta = \Pig w - w$ satisfies
\begin{subequations}\label{eq:1dPig-def}
\begin{align}
	\ipi{\eta}{v_h} =&\; 0 \quad \forall v_h \in P^{k-1}(I_i) \quad \forall i = 1,\cdots, N,\\
	{\eta}_{i+\hf}^- =&\; 0\quad \forall i  = 1,\cdots, N.
\end{align}
\end{subequations}

The Gauss--Radau projection $\Pig$ satisfies the following approximation estimate and superconvergence property. 
\begin{PROP}[Gauss--Radau projection]\label{prop:1dPig}
	$\Pig: H^1(\Omega_h)\to V_h^k$ in \eqref{eq:1dPig-def} is well-defined. Furthermore, when $\flat \geq 1$, 
	\begin{enumerate}
		\item for $w\in  H^{\flat}(\Omega_h)$, we have the approximation estimate
		\begin{equation}
			\nm{(I- \Pig) w} \leq C \nm{w}_{H^{\min(\ell,\flat)}}h^{\min(\ell,\flat)} \qquad \forall 0\leq \ell \leq k+1;
		\end{equation}
		\item for $w\in C(\overline{\Omega})\cap H^{1}(\Omega_h)$,
		we have the superconvergence property
		\begin{equation}\label{eq:supc1d}
			\ip{(\Pik L - L_h\Pig)w}{v_h} = 0.
		\end{equation}
	\end{enumerate} 
	
\end{PROP}

\subsubsection{2D case} 

\;\smallskip

Let $H^\ell(\Omega_h) = \{v\in L^2(\Omega): v|_{K_{ij}}\in H^\ell(K_{ij})\  \forall i,j\}$ and  $\nm{v}_{H^\ell} = \left(\sum_{i,j} \nm{v}_{H^\ell(K_{ij})}^2\right)^{1/2}$. In addition, we assume uniform meshes for which $h_x^i \equiv h_x$ and $h_y^j \equiv h_y$ are constants. The following projection \cite{liu2020optimal} is used for the optimal error estimate of the $P^k$-DG method in 2D. $\Pig: H^1(\Omega_h) \to V_h^k$ is an operator such that for $w\in H^1(\Omega_h)$, the projection error $\eta = \Pig w - w$ satisfies
\begin{subequations}\label{eq:2dPig-def}
\begin{align}
	&\ip{\eta}{1}_{K_{ij}} =\; 0 \quad \forall  i,j,\\
	&\ip{\eta}{\beta^x(v_h)_x+\beta^y(v_h)_y}_{K_{ij}}
	-\int_{x_{i-\hf}}^{x_{i+\hf}}\beta^y \eta\left(x,y_{j+\hf}^-\right)\left( v_h\left(x,y_{j+\hf}^-\right)- v_h\left(x,y_{j-\hf}^+\right)\right)\dd x\nonumber\\
	&-\int_{y_{j-\hf}}^{y_{j+\hf}}  \beta^x \eta\left(x_{i+\hf}^-,y\right) \left( v_h\left(x_{i+\hf}^-,y\right)-v_h\left(x_{i-\hf}^+,y\right)\right)\dd y = 0 \quad \forall v_h \in P^{k}(K_{ij})\ \forall i,j.
\end{align}
\end{subequations}

The special projection $\Pig$ satisfies the following approximation and superconvergence properties, as explained in Lemma \ref{lem:2dPig}. The case $\ell = k+2$ for \eqref{eq:supc2d} is proved in \cite{liu2020optimal}. We have supplemented a sketched proof in Section \ref{sec:prop:2dPig} for the case $1\leq \ell \leq k+1$. 
\begin{LEM}[Liu--Shu--Zhang projection]\label{lem:2dPig}
	Let $\Omega_h$ be a uniform partition in 2D. $\Pig: H^1(\Omega_h) \to V_h^k$ in \eqref{eq:2dPig-def} is well-defined. Furthermore, when $\flat \geq 1$, 
	\begin{enumerate}
		\item for $w\in H^\flat(\Omega_h)$, we have the approximation estimate
		\begin{equation}\label{eq:lszproj-approx}
			\nm{(I- \Pig) w} \leq C \nm{w}_{H^{\min(\ell,\flat)}}h^{\min(\ell,\flat)}\qquad \forall 0\leq \ell \leq k+1;
		\end{equation}
		\item for $w\in C(\overline{\Omega})\cap H^\flat(\Omega_h)$, we have the superconvergence property
		\begin{equation}\label{eq:supc2d}
			\ip{(\Pik L - L_h\Pig)w}{v_h} \leq C \nm{w}_{H^{\min(\ell,\flat)}}h^{\min(\ell,\flat)-1}\nm{v_h}\qquad \forall 1\leq \ell \leq k+2.
		\end{equation}
	\end{enumerate}
\end{LEM}

\subsubsection{Notations for both 1D and 2D cases}

\;\smallskip

In order to obtain the optimal error estimates for both the 1D and 2D cases in one go, in this section, we unify the notations and write properties in Proposition \ref{prop:1dPig} and Lemma \ref{lem:2dPig} as a single proposition. Note that the regularity assumptions for the superconvergence properties in 1D and 2D are different, and we have adopted the stronger assumption below for \eqref{eq:piL-LhPig-g}. 

\begin{PROP}[Special projection $\Pig$]\label{prop:ndPig}
	On both 1D quasi-uniform meshes and 2D uniform meshes, there exists a projection operator $\Pig: H^{1}(\Omega_h)\to V_h^k$, such that when $\flat \geq 1$, 
	\begin{enumerate}
		\item for $w\in H^\flat(\Omega_h)$, we have the approximation estimate
		\begin{align}
		\nm{(I- \Pig) w} \leq C \nm{w}_{H^{\min(\ell,\flat)}}h^{\min(\ell,\flat)} \qquad \forall 0\leq \ell \leq k+1;\label{eq:uni-approx}
		\end{align}
		\item for $w\in C(\overline{\Omega})\cap H^\flat(\Omega_h)$, we have the superconvergence property
		\begin{equation}\label{eq:piL-LhPig-g}
		\ip{(\Pik L - L_h\Pig)w}{v_h} \leq C_d \nm{w}_{H^{\min(\ell,\flat)}}h^{\min(\ell,\flat)-1} \nm{v_h}\qquad  \forall 1\leq \ell \leq k+2.
		\end{equation}
	\end{enumerate}
\end{PROP}

\begin{REM}[Error estimates in other settings]
The optimal error estimates in later sections can be extended to other settings, as long as one can construct a projection operator with the properties outlined in Proposition \ref{prop:ndPig}. For example, using the Generalized Gauss--Radau projection in \cite{meng2016optimal,cheng2017application}, one can prove optimal error estimates for sdA-RKDG schemes with upwind-biased fluxes in 1D; with the projections from \cite{cockburn2008optimal,sun2023generalized}, one can prove optimal error estimates for the sdA-RKDG schemes with upwind or upwind-biased fluxes on special simplex meshes satisfying the so-called flow conditions in two and three dimensions.
\end{REM}

In our error estimates, we need to deal with functions involving both spatial and temporal variables. We define
\begin{equation}
	L^\infty(0,T;H^\flat(\Omega_h)) = \left\{u(\cdot,t):[0,T]\to H^\flat(\Omega_h)\bigg|\sup_{0\leq t\leq T}\nm{u(\cdot,t)}_{H^\flat} <+\infty\right\}
\end{equation} 
and denote by $\nm{u}_{L^\infty(H^\flat)} = \sup_{0\leq t\leq T}\nm{u(\cdot,t)}_{H^\flat}$. Here the ``$\sup$" function should be interpreted as the essential supreme. We will abuse the notation and formally identify $u(\pmb{x},t): \overline{\Omega} \times [0,T]\to \mathbb{R}$ with $u(\cdot,t): [0,T]\to H^\flat(\Omega_h)$, and use $u\in C^{\flat-1}(\overline{\Omega}\times [0,T])\cap L^\infty(0,T;H^\flat(\Omega_h))$ to represent $u(\pmb{x},t)\in C^{\flat-1}(\overline{\Omega}\times [0,T])$ and $u(\cdot,t)\in L^\infty(0,T;H^\flat(\Omega_h))$.

\subsubsection{A few useful estimates}

\;\smallskip

The following lemma, derived from Proposition \ref{prop:ndPig}, will be frequently used. Its proof is postponed to Section \ref{sec:proof-est2}.
\begin{LEM}\label{lem:i-PiL}
	On both 1D quasi-uniform meshes and 2D uniform meshes, we have the following results. 
	\begin{enumerate}
		\item If $w \in H^\flat(\Omega_h)$ and $\flat \geq 1$, then
		\begin{alignat}{2}
		\nm{(I-\PiPi) w} \leq &C \nm{w}_{H^{\min(\ell,\flat)}}h^{\min(\ell,\flat)}, \qquad &\PiPi = 
		\begin{cases}
		\Pi, \Pig,&\quad \forall 0\leq \ell \leq k+1,\\
		\Piz,	 &\quad \forall 0\leq \ell \leq k.
		\end{cases} \label{eq:I-Pi-0}
		\end{alignat}
		\item If $w \in H^\flat(\Omega_h)$ and $\flat \geq 2$, then
		\begin{alignat}{2}	
		\nm{(I-\PiPi) L w} \leq& C \nm{w}_{H^{\min(\ell,\flat)}}h^{\min(\ell,\flat)-1}, \qquad &\PiPi = 
		\begin{cases}
		\Pi, \Pig, &\quad  \forall1\leq \ell \leq k+2,\\
		\Piz, &\quad  \forall1\leq \ell \leq k+1.
		\end{cases}\label{eq:I-Pi-1}
		\end{alignat}
		\item If $w \in C^{\flat-1}(\overline{\Omega})\cap H^\flat(\Omega_h)$ and $\flat \geq 1$, then
		\begin{equation}\label{eq:PiL-LhPig-i}
		\nm{\left(\Pik L - L_h \Pig\right)L^{i}w}\leq C_d \nm{w}_{H^{\min(\ell,\flat)}}h^{\min(\ell,\flat)-i-1}\qquad \forall1\leq \ell \leq k+2, 
		\end{equation}
		where $0\leq i \leq \min(\ell,\flat)-1$.
		\item 	If $w \in C^{\flat-1}(\overline{\Omega})\cap H^{\flat}(\Omega_h)$ and $\flat \geq 1$, then 
		\begin{equation}
		\nm{(\Pik L- \LL_h \Pig) w}\leq C\nm{w}_{H^{\min(\ell,\flat)}}h^{\min(\ell,\flat)-1}\qquad \forall1\leq \ell \leq k+1.\label{eq:piL-Lpi-1-1}
		\end{equation}
	\end{enumerate}
\end{LEM}

\subsection{Optimal error estimates: the second-order case}

\subsubsection{Main results}

\;\smallskip

In this section, we prove the optimal error estimate of the sdA-RK2DG1 scheme \eqref{eq:sdArkdg2}, or equivalently, \eqref{eq:sds21line}. 
\begin{THM}[Optimal error estimate of sdA-RK2DG1]\label{thm:err2}
	Consider the sdA-RK2DG1 scheme \eqref{eq:sdArkdg2} on either 1D quasi-uniform meshes or 2D uniform meshes. Suppose the exact solution $u$ to \eqref{eq:adv} is sufficiently smooth, for example, with $u \in C^2(\overline{\Omega}\times[0,T])\cap L^\infty(0,T;H^3(\Omega_h))$, then the sdA-RK2DG1 scheme \eqref{eq:sdArkdg2} admits the optimal error estimate
	\begin{equation}
		\nm{u^n - u_h^n}\leq \nm{u^0 - u_h^0}+ C\left(t^n+1\right)\nm{u}_{L^\infty(H^{3})}\left(\tau^2 + h^2\right).
	\end{equation}
	Here $t^n = n\tau\leq T$, $u^n = u(\pmb{x},t^n)$, and $u^0 = u(\pmb{x},0)$. 
\end{THM}

The key to the proof of Theorem \ref{thm:err2} is to construct the special approximation operator in Lemma \ref{lem:newproj}, whose proof is postponed to Section \ref{sec:proof-proj}. This approximation operator serves a role similar to that of the projection operator in classical analysis for the DG method. However, we use a different name because it may not be an actual projection operator in the mathematical sense.

\begin{LEM}[Approximation operator for sdA-RK2DG1]\label{lem:newproj}
	Suppose $w\in C^{\flat-1}(\overline{\Omega})\cap H^\flat(\Omega_h)$ and $\flat \geq 2$. Then
	\begin{align}\label{eq:def-Pis}
		\Pis w = \left(I+\frac{\tau}{2}\LL_h\right)^{-1}\Pig \left(w+\frac{\tau}{2} L w\right)
	\end{align}
	is well-defined when ${\lambda}\ll 1$. Moreover, it admits the approximation estimate
	\begin{equation}\label{eq:Pis2}
	\nm{(I- \Pis) w} \leq C \nm{w}_{H^{\min(\ell,\flat)}}h^{\min(\ell,\flat)}\qquad \forall1\leq \ell \leq k+1.
	\end{equation}
\end{LEM}

\subsubsection{Proof of Theorem \ref{thm:err2}}\label{sec:proof-thm:err2}

\;\smallskip

\underline{\em Reference function.} Using the Taylor expansion and the equation $u_t = Lu$, it yields
\begin{equation}\label{eq:taylor2}
	u^{n+1} = u^n + \tau u_t^n + \frac{\tau^2}{2}u_{tt}^n + \tau \rho^n = u^n + {\tau} L u^{n} + \frac{\tau^2}{2} L^2 u^{n} +\tau \rho^n,
\end{equation}
where $\rho^n$ is the remainder term in the Taylor expansion defined as
\begin{equation}\label{eq:rho-2}
\rho^n(\pmb{x}) = \frac{1}{2\tau}\int_0^\tau u_{ttt}(\pmb{x},t^n+\kappa)(\tau-\kappa)^2\dd \kappa.
\end{equation}
Here one can show that $u_{ttt} = \partial_t L^2 u = L^2 u_t = L^3 u$ is well-defined. The exchange of $L^2$ and $\partial_t$ in the second equality can be justified by using the definition of the weak derivatives, integrating by parts, and exchanging the partial derivatives at the smooth test function.  

Applying the projection $\Pi$ to both sides of \eqref{eq:taylor2}, it gives
\begin{equation}
	\Pik u^{n+1} = \Pik u^n + {\tau} \Pik L u^{n} + \frac{\tau^2}{2} \Pik L^2 u^{n} +\tau \Pik \rho^n,
\end{equation}
which, after adding and subtracting terms, can be rewritten as 
\begin{equation}\label{eq:pisu}
	\Pis u^{n+1} = \Pis u^n + {\tau} L_h \Pis u^{n} + \frac{\tau^2}{2} L_h\LL_h \Pis u^{n} +\tau \left(\Pik \rho^n + \nu_h^n+ \zeta_h^n\right).
\end{equation}
Here
\begin{align}
	\nu_h^n =& \tau^{-1}\left(\Pis-\Pik\right)\left(u^{n+1}-u^n\right):=\left(\Pis-\Pik\right)\varphi^n,\\
	\zeta_h^n =& \Pik L\left(I+\frac{\tau}{2}L \right)u^n - L_h\left(I+\frac{\tau}{2}\LL_h\right)\Pis u^n\label{eq:zeta-2}.
\end{align}
We will use \eqref{eq:pisu} as a reference function in our error estimate.  \small

\underline{\em Estimate of $\Pik \rho^n$.} To estimate \eqref{eq:rho-2}, note that $u_{ttt} = L^3 u$. One can get
\begin{equation}\label{eq:nmrho2}
		\nm{\Pik \rho^n}\leq \nm{\rho^n} \leq C\sup_{0\leq t\leq t^n} \nm{u_{ttt}(\cdot,t)}\tau^2  = C\sup_{0\leq t\leq t^n} \nm{L^3u(\cdot,t)}\tau^2\leq C\nm{u}_{L^\infty(H^3)}\tau^2.
\end{equation} 

\underline{\em Estimate of $\nu_h^n$.} Note that 
\begin{equation}
\varphi^n = \tau^{-1}(u^{n+1}-u^n) = \tau^{-1}\int_0^\tau u_t(\pmb{x},t^n+\kappa)\dd \kappa = \tau^{-1}\int_0^\tau L u(\pmb{x},t^n+\kappa)\dd \kappa.
\end{equation}
One can apply Jensen's inequality to see that $\nm{\varphi^n}_{H^{\ell}} \leq \sup_{0\leq \kappa\leq \tau} \nm{L u(\cdot,t^n+\kappa)}_{H^{\ell}} \leq \nm{u}_{L^\infty(H^{\ell+1})}$. Therefore, $u\in L^\infty(0,T;H^3({\Omega_h}))$ implies $\varphi^n \in H^2({\Omega_h})$. Using the triangle inequality, \eqref{eq:I-Pi-0}, and \eqref{eq:Pis2}, with $\flat = 2$ and $\ell = 2$, we have
\begin{equation}\label{eq:nmnu2}
	\nm{\nu_h^n}\leq \nm{(I-\Pis)\varphi^n}+\nm{(I-\Pi)\varphi^n}\leq  C\nm{\varphi^n}_{H^2}h^2\leq C\nm{u}_{L^\infty(H^3)}h^2.
\end{equation}

\underline{\em Estimate of $\zeta_h^n$.} According to the definition of $\Pis$ in Lemma \ref{lem:newproj}, we have 
\begin{equation}\label{eq:zeta2}
	\begin{aligned}
		\zeta_h^n =& \Pik L\left(I+\frac{\tau}{2}L \right)u^n - L_h\Pig \left(I+\frac{\tau}{2} L\right)u^n
		= \left(\Pik L - L_h \Pig\right)\left(I+\frac{\tau}{2} L \right)u^n.
	\end{aligned}	
\end{equation}
Apply the triangle inequality and use \eqref{eq:PiL-LhPig-i} in Lemma \ref{lem:i-PiL} with $i = 0,1$, $\ell = k+2 = 3$, and $\flat = 3$. Then
\begin{equation}\label{eq:est-zeta-2}
	\begin{aligned}
		\nm{\zeta_h^n} \leq& \nm{\left(\Pik L - L_h \Pig\right) \left(u^n + \frac{\tau}{2}L u^n\right)}\\
		\leq & \nm{\left(\Pik L - L_h \Pig\right) u^n} + \frac{\tau}{2}\nm{\left(\Pik L - L_h \Pig\right) L u^n}\\
		\leq & C\nm{u}_{L^\infty(H^{3})}\left(h^{2}+\tau h\right) \leq  C \nm{u}_{L^\infty(H^3)}h^2.
	\end{aligned}
\end{equation}

\underline{\em Estimate of the numerical error.} Subtract \eqref{eq:sds21line} from \eqref{eq:pisu}, and define $\xi_h^n = \Pis u^n - u_h^n$. It gives
\begin{equation}
\begin{aligned}
	\xi_h^{n+1} =& \Pis \xi_h^n + {\tau} L_h \Pis \xi_h^n + \frac{\tau^2}{2} L_h\LL_h \Pis \xi_h^n  + \tau (\Pik\rho^n +\nu_h^n+ \zeta_h^n )\\
	=&\RRR_{h,2} \xi_h^n + \tau (\Pik \rho^n +\nu_h^n+ \zeta_h^n ).
	\end{aligned}
\end{equation}
Applying the triangle inequality  and the monotonicity stability $\nm{\RRR_{h,2}\xi_h^n}\leq \nm{\xi_h^n}$ in Theorem \ref{thm:main-stab2}, it gives
\begin{equation}
	\nm{\xi_h^{n+1}} \leq  \nm{\xi_h^n} +  \tau \left(\nm{\Pik\rho^n} +\nm{\nu_h^n}+ \nm{\zeta_h^n}\right).
\end{equation}
With the estimates for $\nm{\Pik \rho^n}$ in \eqref{eq:nmrho2}, for $\nm{\nu_h^n}$ in \eqref{eq:nmnu2}, and for $\nm{\zeta_h^n}$ in \eqref{eq:est-zeta-2}, we have
\begin{equation}\label{eq:errstab}
	\nm{\xi_h^{n+1}} \leq \nm{\xi_h^n} + C\tau \nm{u}_{L^\infty(H^3)}\left(\tau^2 + h^2\right).
\end{equation}
Note that $\tau n = t^n$. Repeatedly applying \eqref{eq:errstab} for $n$ times gives 
\begin{equation}
	\begin{aligned}
		\nm{\xi_h^n}\leq \nm{\xi_h^0} + Ct^n\nm{u}_{L^\infty(H^{3})}\left(\tau^2 + h^2\right).
	\end{aligned}
\end{equation}
Recalling the approximation estimate \eqref{eq:Pis2} in Lemma \ref{lem:newproj}, we have  $\nm{u^n-\Pis u^n}\leq C\nm{u^n}_{H^2}h^2$ and $\nm{u^0-\Pis u^0}\leq C\nm{u^0}_{H^2}h^2$. With the triangle inequality, it gives
\begin{equation}
	\begin{aligned}
			\nm{u^n - u_h^n}\leq& \nm{u^n-\Pis u^n}+\nm{\xi_h^n}\leq \nm{\xi_h^0}+ C\left(t^n+1\right)\nm{u}_{L^\infty(H^{3})}\left(\tau^2 + h^2\right)\\
			\leq& \nm{u^0 - \Pis u^0 } + \nm{u^0 - u_h^0}+ C\left(t^n+1\right)\nm{u}_{L^\infty(H^{3})}\left(\tau^2 + h^2\right)\\
			\leq& \nm{u^0 - u_h^0}+ C\left(t^n+1\right)\nm{u}_{L^\infty(H^{3})}\left(\tau^2 + h^2\right).
	\end{aligned}
\end{equation}

\begin{REM}[Idea behind the definition of $\Pis$]
	The relationship \eqref{eq:zeta-2} and the cancellation in \eqref{eq:zeta2} motivate the definition of $\Pis$ in \eqref{eq:def-Pis}. In order to get rid of $(I+\tau \LL_h/2)$ in $\zeta_h^n$ in \eqref{eq:zeta-2}, we introduce the term $(I+\tau \LL_h/2)^{-1}$ in the definition of $\Pis$. Then in order to pull out the operator $\Pik L - L_h\Pig$ in \eqref{eq:zeta2}, we impose the term $\Pig(I+\tau L/2)$.
\end{REM}

\begin{REM}[Alternative error estimates for RK2DG1]
	When we set $\LL_h = L_h$, the sdA-RK2DG1 scheme \eqref{eq:sdArkdg2} formally retrieves the RK2DG1 scheme. It is worth noting that in this case, all the analysis in this section still holds. Consequently, it provides an alternative proof of the optimal error estimates for the standard RK2DG1 scheme. Further discussions on this matter can be found in Appendix \ref{sec:multistage}.
\end{REM}

\subsection{Optimal error estimates: general cases}

\subsubsection{Main results}

\;\smallskip

Theorem \ref{thm:err-g} presents the optimal error estimates of the sdA-RK$r$DG$k$ scheme for generic values of $r$ and $k$. Its proof is given in Section \ref{sec:err-g}. 

\begin{THM}[Optimal error estimates of sdA-RK$r$DG$k$]\label{thm:err-g}
	Consider the sdA-RK$r$DG$k$ scheme \eqref{eq:sdA-rkdg-Butcher} (or \eqref{eq:sdA-RKDG}) on either 1D quasi-uniform meshes or 2D uniform meshes. Suppose the exact solution $u$ to \eqref{eq:adv} is sufficiently regular, for example, $u\in C^{\flat-1}([0,T]\times \overline{\Omega})\cap L^\infty(0,T;H^\flat(\Omega_h))$. Let $q = \min(r,k+1,\flat-1)$. Then we have
	\begin{equation}
		\nm{u^n - u_h^n}\leq C_1(t^n)\nm{u^0 - u_h^0}+ C_2(t^n)\nm{u}_{L^\infty(H^{q+1})}\left(\tau^q+ h^q\right),
	\end{equation}
	where $t^n = n\tau\leq T$, $u^n = u(\pmb{x},t^n)$, $u^0 = u(\pmb{x},0)$, and
	\begin{enumerate}
		\item $C_1(t^n) = C_2(t^n) = e^{Ct^n}$ if $\RRR_h$ is weakly$(\gamma)$ stable and $\lambda \leq \lambda_0 h^{1/(\gamma-1)}$;
		\item $C_1(t^n) = 1$ and $C_2(t^n) = C(t^n+1)$ if $\RRR_h$ is monotonically stable and $\lambda \ll 1$;
		\item $C_1(t^n) = C$ and $C_2(t^n) = C(t^n+1)$ if $\RRR_h$ is strongly$(m_\star)$ stable and $\lambda \ll 1$.
	\end{enumerate}
\end{THM}

The key step for proving Theorem \ref{thm:err-g} is to construct the special approximation operator in Lemma \ref{lem:newproj-general}, whose proof is postponed to Section \ref{newproj-general}.

\begin{LEM}[Approximation operator for sdA-RK$r$DG$k$]\label{lem:newproj-general}
	For $w\in C^{\flat-1}(\overline{\Omega})\cap H^\flat(\Omega_h)$ and $q \leq \min(r,k+1,\flat-1)$, 
\begin{equation}
	\Pis w = \left(\sum_{i = 1}^s \alpha_i (\tau \LL_h)^{i-1}\right)^{-1}\left(\Pig\left(\sum_{i = 1}^q (i!)^{-1} (\tau L)^{i-1}\right)w + \tau^q\sum_{i = q+1}^s \alpha_i (\tau \LL_h)^{i-1-q{}} \Pig L^q w\right)
\end{equation}
is well-defined. Moreover, it satisfies the approximation estimate
\begin{equation}\label{eq:ae-Pis}
\nm{(I-\Pis) w} \leq C\nm{w}_{H^{q}}h^{q}.
\end{equation}
\end{LEM}

\begin{REM}[Regularity assumption required by $\Pis$]
	Following Lemma \ref{lem:newproj} for the second-order case, one may naturally consider defining 
	\begin{equation}
		\Pis w = \left(\sum_{i = 1}^s \alpha_i (\tau \LL_h)^{i-1}\right)^{-1}\left(\Pig\left(\sum_{i = 1}^s \alpha_i (\tau L)^{i-1}\right)w\right).
	\end{equation}
	However, this definition imposes a restrictive regularity assumption on $w \in H^s(\Omega_h)$, where $s\geq r$. For example, for $r$-stage $r$th-order RK schemes that are strongly$(m_\star)$ stable, as we need to combine $m_\star$ steps into a single step in the analysis, we could have $s = m_\star r$. However, typically, we anticipate that $H^r$ regularity should be adequate to achieve the desired accuracy. To circumvent the need for this additional regularity requirement, we modify high-order derivatives $L^{i-1}w$ with $q+1\leq i\leq s$ as $\LL_h^{i-1-q}\Pig L^q w$. This modification shares a similar flavor to the cutting-off technique in the error analysis of the original RKDG method \cite{xu2020superconvergence,xu2020error}.
\end{REM}

\subsubsection{Proof of Theorem \ref{thm:err-g}}\label{sec:err-g}

\;\smallskip

\underline{\em Reference function.} Since $u\in C^{\flat-1}([0,T]\times \overline{\Omega})$ and $q\leq \flat-1$, we have $\partial_t^i u^n = L^i u^n$ for $1\leq i\leq q$. Thus, using the Taylor expansion, we have
\begin{equation}\label{eq:taylorg}
	u^{n+1} = \sum_{i = 0}^{q}(i!)^{-1}\tau^{i}\partial_t^i u^n + \tau \rho^n = \sum_{i = 0}^{q}(i!)^{-1}\tau^{i}L^i u^n + \tau \rho^n. 
\end{equation}
The remainder $\rho^n$ is defined as 
\begin{equation}\label{eq:rho-g}
\rho^n(\pmb{x}) = \frac{1}{\tau q!}\int_0^\tau \partial_t^{(q+1)}u(\pmb{x}, t^n+\kappa)(\tau-\kappa)^{q}\dd \kappa.
\end{equation}
Recalling that $q\leq \flat -1$, one can show that $\partial_t^{(q+1)}u = L^{q+1} u$ is well-defined under the given regularity assumption.

We apply $\Pi$ to both sides of \eqref{eq:taylorg} to get
\begin{equation}
	\Pik u^{n+1} = \Pik u^n + \sum_{i = 1}^q(i!)^{-1}\tau^{i}\Pik L^i u^n + \tau \Pik \rho^n. 
\end{equation}
By adding and subtracting terms, it can be rewritten as 
\begin{equation}\label{eq:pis-general}
	\Pis u^{n+1} = \Pis u^n + \sum_{i = 1}^s \alpha_i \tau^i L_h \LL_h^{i-1} \Pis u^n + \tau \left(\Pik \rho^n+ \nu_h^n + \zeta_h^n\right). 
\end{equation}
where 
\begin{align}
	\nu_h^n =& \tau^{-1}\left(\Pis-\Pik\right)\left(u^{n+1}-u^n\right):=(\Pis-\Pi)\varphi^n,\\
	\zeta_h^n =& \sum_{i=1}^q (i!)^{-1} \tau^{i-1} \Pik L^{i}u^n - L_h\left(\sum_{i=1}^s \alpha_i \left(\tau \LL_h\right)^{i-1}\right)\Pis u^n.
\end{align}
\eqref{eq:pis-general} will be used as the reference function for our error estimate. \smallskip

\underline{\em Estimates of $\Pik\rho^n$ and $\nu_h^n$.} The estimates of ${\Pik \rho^n}$ and ${\nu_h^n}$ are similar to those in the second-order case. Further details are omitted.
\begin{equation}\label{eq:Pirho-g}
\nm{\Pik \rho^n} \leq C\sup_{0\leq t\leq t^n}\nm{\partial_t^{(q+1)} u(\cdot,t)} \tau^{q} = C\sup_{0\leq t\leq t^n}\nm{L^{q+1} u(\cdot,t)} \tau^{q}\leq C\nm{u}_{L^\infty(H^{q+1})}\tau^{q}.
\end{equation}
\begin{equation}\label{eq:nu-g}
\nm{\nu_h^n}\leq \nm{(I-\Pis)\varphi^n}+\nm{(I-\Pi)\varphi^n}\leq  C\nm{\varphi^n}_{H^q}h^q\leq C\nm{u}_{L^\infty(H^{q+1})}h^q.
\end{equation}
\;

\underline{\em Estimate of $\zeta_h^n$.} We then estimate the term $\zeta_h^n$. 

\begin{LEM}\label{eq:est-zetag}
	When $\lambda \ll 1$, we have 
	\begin{equation}\label{eq:est-zetag}
		\nm{\zeta_h^n} \leq C \nm{u}_{L^\infty(H^{q+1})}\left(\tau^q + h^q\right).
	\end{equation}
\end{LEM}
\begin{proof}
For $\zeta_h^n$, we can use the definition of $\Pis$ to get
\begin{equation}
	\begin{aligned}
	\zeta_h^n 
	=& \left(\Pik L - L_h \Pig\right)\sum_{i=1}^q (i!)^{-1} (\tau L)^{i-1}u^n + \tau^q L_h\sum_{i = q+1}^s \alpha_i (\tau \LL_h)^{i-1-q{}} \Pig L^q u^n\\
	:=& \zeta_h^{n,1} + \zeta_h^{n,2}.
	\end{aligned}
\end{equation}

For $\zeta_h^{n,1}$, we can apply a similar estimate as that for the second-order case. One can use the triangle inequality and \eqref{eq:PiL-LhPig-i} in Lemma \ref{lem:i-PiL} with $\ell = q+1\leq \flat$ to get
\begin{equation}\label{eq:nmzeta1}
	\begin{aligned}
		\nm{\zeta_h^{n,1}} 
		\leq & \sum_{i=1}^q (i!)^{-1} \tau^{i-1} \nm{\left(\Pik L - L_h \Pig\right)L^{i-1}u^n}\\
		\leq & C_d\sum_{i=1}^q \tau^{i-1} \nm{u^n}_{H^{q+1}}h^{q+1-i}\leq  C_d\nm{u}_{L^\infty(H^{q+1})}h^{q}.
	\end{aligned}
\end{equation}
For $\zeta_h^{n,2}$, we separate the first term in the summation to get
\begin{equation}
		\zeta_h^{n,2} = \tau^q \left(\alpha_{q+1} L_h\Pig L^q u^n + (\tau L_h)\sum_{i = q+2}^s \alpha_i (\tau \LL_h)^{i-2-q} \LL_h\Pig L^q u^n\right).
\end{equation}
Then we apply the triangle inequality and the estimates $\tau\nm{\LL_h}\leq \tau \nm{L_h} \leq C\lambda$ to get
\begin{equation}\label{eq:nmzeta2}
	\begin{aligned}
	\nm{\zeta_h^{n,2}} \leq & C \tau^q \left(\nm{L_h\Pig L^q u^n }+ \nm{\LL_h\Pig L^q u^n}\right)\\
	\leq & C \tau^q \left(\nm{(L_h\Pig - \Pik L) L^q u^n} +  \nm{(\LL_h\Pig-\Pik L) L^q u^n}+2\nm{\Pik L^{q+1} u^n} \right).
	\end{aligned}
\end{equation}
We can apply \eqref{eq:PiL-LhPig-i} in Lemma \ref{lem:i-PiL} with $w = u^n$ and $\ell = q+1\leq \flat$ to get 
\begin{equation}
\nm{(L_h\Pig - \Pik L) L^q u^n} \leq C_d \nm{u^n}_{H^{q+1}} \leq  C_d\nm{u}_{L^\infty(H^{q+1})}.
\end{equation}
Apply \eqref{eq:piL-Lpi-1-1} in Lemma \ref{lem:i-PiL} with $\ell = 1$ and $w = L^qu^n \in C^{\flat -q-1}(\overline{\Omega})\cap H^{\flat-q}(\Omega_h)\subset C(\overline{\Omega})\cap H^1(\Omega_h)$. It gives
\begin{equation}
	\nm{(\LL_h\Pig-\Pik L) L^q u^n} \leq C\nm{L^q u^n}_{H^{1}}\leq C\nm{u}_{L^\infty(H^{q+1})}.
\end{equation}
Moreover, 
\begin{equation}
	\nm{\Pik L^{q+1} u^n}\leq C\nm{u^n}_{H^{q+1}}\leq C\nm{u}_{L^\infty(H^{q+1})}.
\end{equation} 
Substituting these estimates into \eqref{eq:nmzeta2} gives
\begin{equation}\label{eq:nmzeta2-simp}
	\nm{\zeta_h^{n,2}} \leq C \nm{u}_{L^\infty(H^{q+1})}\tau^q. 	
\end{equation}
Combining \eqref{eq:nmzeta1} and \eqref{eq:nmzeta2-simp}, we complete the proof of \eqref{eq:est-zetag}. 
\end{proof}

\underline{\em Estimate of the numerical error.} Let $\xi_h = \Pis u - u_h$. Subtracting \eqref{eq:sdA-RKDG} from \eqref{eq:pis-general} gives 
\begin{equation}\label{eq:err-xi-g}
	\xi_h^{n+1} = \RRR_h \xi_h^n + \tau (\Pis\rho^n + \nu_h^n+ \zeta_h^n).
\end{equation}
Apply the triangle inequality and the estimates in \eqref{eq:rho-g}, \eqref{eq:nu-g}, and \eqref{eq:est-zetag}. Then we get
\begin{equation}\label{eq:xi-g}
	\begin{aligned}
		\nm{\xi_h^{n+1}} 
		\leq& \nm{\RRR_h \xi_h^n} + \tau\left(\nm{\Pis \rho^n}+\nm{\nu_h^n} +\nm{\zeta_h^n}\right) \\
		\leq& \nm{\RRR_h \xi_h^n} + C\nm{u}_{L^\infty(H^{q+1})}\tau \left(\tau^q + h^q\right).
	\end{aligned}
\end{equation}

Below, we discuss the cases with different types of stability. 
\begin{itemize}
	\item Suppose $\RRR_h$ is weakly($\gamma$) stable. Then $\nm{\RRR_h \xi_h^n} \leq (1+C\lambda^\gamma)\nm{\xi_h^n}$, and we have
\begin{equation}
	\nm{\xi_h^{n+1}} 
	\leq  (1+ C\lambda^{\gamma})\nm{\xi_h^n} + C\nm{u}_{L^\infty(H^{q+1})}\tau \left(\tau^q + h^q\right).
\end{equation}
Assuming $\lambda \leq \lambda_0 h^{1/(\gamma-1)}$, we can apply discrete Gronwall's inequality to get
\begin{equation}\label{eq:est-xi-g-w}
	\begin{aligned}
		\nm{\xi_h^{n}} 
		\leq e^{Ct^n}\left(\nm{\xi_h^0} + \nm{u}_{L^\infty(H^{q+1})} \left(\tau^q + h^q\right)\right).
	\end{aligned}
\end{equation}

\item Suppose $\RRR_h$ is monotonically stable. Then $\nm{\RRR_h \xi_h^n} \leq \nm{\xi_h^n}$, and we have
\begin{equation}\label{eq:est-xi-g-m}
	\nm{\xi_h^{n}} 
	\leq \nm{\xi_h^0} + Ct^{n}\nm{u}_{L^\infty(H^{q+1})} \left(\tau^q + h^q\right).
\end{equation}

\item Suppose $\RRR_h$ is strongly($m_\star$) stable. First, noting that $\nm{\RRR_h} \leq C$ under the condition $\lambda \ll 1$, \eqref{eq:xi-g} implies
	\begin{equation}\label{eq:err-sm1}
		\nm{\xi_h^{n+1}}  \leq C\nm{\xi_h^n} +  C\nm{u}_{L^\infty(H^{q+1})}\tau \left(\tau^q + h^q\right).
	\end{equation}
	Invoking \eqref{eq:err-sm1} for $\ell$ times, with $\ell \leq m_\star$, we have 
	\begin{equation}\label{eq:err-sm11}
		\nm{\xi_h^{n+\ell}}  \leq C\nm{\xi_h^n} +  C\nm{u}_{L^\infty(H^{q+1})}\tau \left(\tau^q + h^q\right).
	\end{equation}
	Next, applying \eqref{eq:err-xi-g} for $m_\star$ times, we get
	\begin{equation}\label{eq:err-sm2}
		\begin{aligned}
			\xi_h^{n+m_\star} =& \RRR_h \xi_h^{n+m_\star-1} + \tau (\Pis\rho^{n+m_\star-1} + \nu_h^{n+m_\star-1}+ \zeta_h^{n+m_\star-1})\\
			=& \RRR_h^{m_\star} \xi_h^{n} + \tau \sum_{\ell = 1}^{m_\star} \RRR_h^{\ell-1} (\Pis\rho^{n+m_\star-\ell} + \nu_h^{n+m_\star-\ell}+ \zeta_h^{n+m_\star-\ell}).
		\end{aligned}
	\end{equation}
	By the strong$(m_\star)$ stability $\nm{\RRR_h^{m_\star}w_h}\leq \nm{w_h}$, the estimate $\nm{\RRR_h}\leq C$, and the estimates of $\Pis\rho^{m}$, $\nu_h^{m}$, $\zeta_h^{m}$, one can get 
\begin{equation}\label{eq:err-sm3}
	\begin{aligned}
	\nm{\xi_h^{n+m_\star}} 
	\leq& \nm{\RRR_h^{m_\star} \xi_h^n} + \tau \sum_{\ell = 1}^{m_\star} \nm{\RRR_h}^{\ell-1} (\nm{\Pis\rho^{n+m_\star-\ell}} + \nm{\nu_h^{n+m_\star-\ell}}+ \nm{\zeta_h^{n+m_\star-\ell}})\\
	\leq& \nm{\xi_h^n} + C\tau\nm{u}_{L^\infty(H^{q+1})}\left(\tau^q + h^q\right).
	\end{aligned}
\end{equation}
Finally, applying \eqref{eq:err-sm11} and then \eqref{eq:err-sm3} for $\lfloor n/m_\star\rfloor$ times, it yields
\begin{equation}\label{eq:err-sm4}
	\begin{aligned}
	\nm{\xi_h^{n}}  \leq& C\nm{\xi_h^{\lfloor n/m_\star\rfloor m_\star}} +  C\tau\nm{u}_{L^\infty(H^{q+1})} \left(\tau^q + h^q\right)\\
	 \leq& C\nm{\xi_h^0} +  C\left\lfloor {n}/{m_\star}\right\rfloor\tau\nm{u}_{L^\infty(H^{q+1})} \left(\tau^q + h^q\right)\\
	 \leq& C\nm{\xi_h^0} +  Ct^n\nm{u}_{L^\infty(H^{q+1})} \left(\tau^q + h^q\right).
	\end{aligned}
\end{equation}
\end{itemize}

The rest of the proof can be completed by applying the triangle inequality to \eqref{eq:est-xi-g-w}, \eqref{eq:est-xi-g-m}, and \eqref{eq:err-sm4}, respectively. The derivation is similar to the second-order case, and details are omitted. Different constants can be attributed to the different types of stability.

\begin{REM}[Alternative error estimates for RK$r$DG$k$]
	Similar to the second-order case, when $\LL_h = L_h$, the analysis above gives an alternative proof for the optimal error estimates of the RK$r$DG$k$ method.
\end{REM}

\section{Numerical verification}\label{sec:num}

In this section, we verify our stability and error analysis presented in the previous sections through numerical experiments.

For the 1D tests, we consider the linear advection equation $u_t + u_x = 0$ on the interval $x\in [0,1]$, discretized with $N$ mesh cells ($h = 1/N$). For the 2D tests, we consider the equation $u_t + u_x + u_y = 0$ on the domain $(x,y)\in [0,1]\times [0,1]$, discretized with $N\times N$ mesh cells ($h_x = h_y = 1/N$). Both cases are subject to periodic boundary conditions.

In this study, we focus solely on linear advection equations. Additional numerical tests involving problems such as the nonlinear Burgers equation can be found in \cite{chen2024runge}.

We employ $r$-stage $r$th-order RK schemes with $r = 2,3,4,5$ \cite{gottlieb2001strong}. For linear autonomous problems with a fixed $r$, all $r$-stage $r$th-order RK schemes are equivalent to the $r$th-order Taylor method, given by $R_h = \sum_{i = 0}^r (i!)^{-1}\tau^i L_h^i$. Correspondingly, $\RRR_h = I + \tau L_h\sum_{i = 1}^r(i!)^{-1}\tau^{i-1}\LL_h^{i-1}$.

\subsection{Stability tests}

\;\smallskip

	In this section, we verify different types of stability of the RKDG and sdA-RKDG schemes. To this end, we follow \cite{xu2019L2,xu2024sdf} to compute 
	\begin{equation}\label{eq:delta}
	\delta = \max\left(\nm{(K_h)^m}^2-1,10^{-16}\right)\quad\text{with}\quad
	K_h = 
	\begin{cases}
	R_h&\text{for RKDG schemes,}\\
	\RRR_h&\text{for sdA-RKDG schemes.}
	\end{cases}
	\end{equation} 
	Here the operator norm can be evaluated by computing the 2-norm of the (sdA-)RKDG matrices assembled under the orthonormal Legendre basis. With $\tau = \text{CFL}/(dN)$ for the $d$-dimensional test ($d =1,2$), we note that $\delta$ is a function of the CFL number, $N$, and $m$.
	
\begin{figure}[h!]
		\centering
	\subfigure[(sdA-)RK2DG1, $m = 1$]{\includegraphics[width=0.32\textwidth]{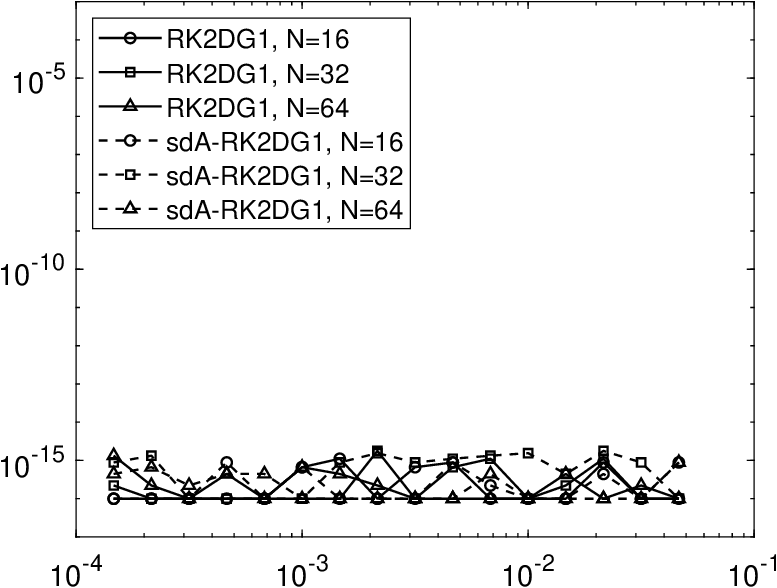}\label{fig:1d-RK2DG1}}
	\subfigure[(sdA-)RK3DG2, $m = 1$]{\includegraphics[width=0.32\textwidth]{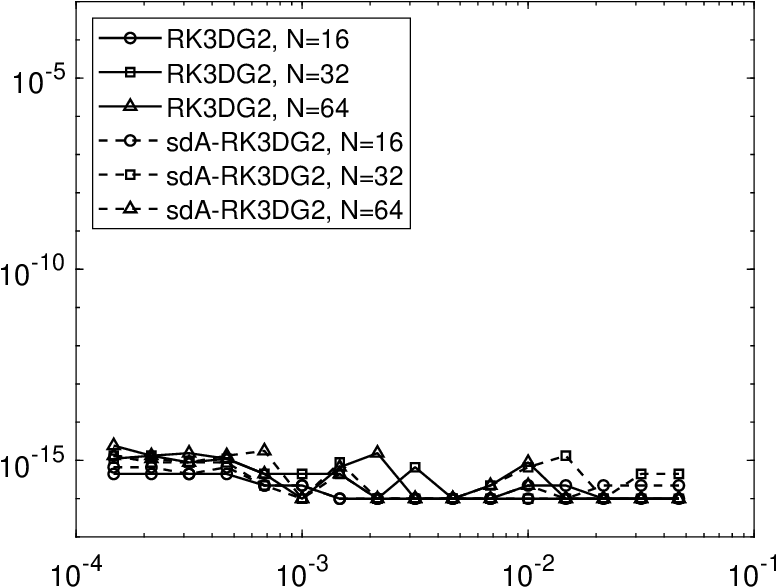}\label{fig:1d-RK3DG2}}
	\subfigure[(sdA-)RK4DG3, $m = 1$]{\includegraphics[width=0.32\textwidth]{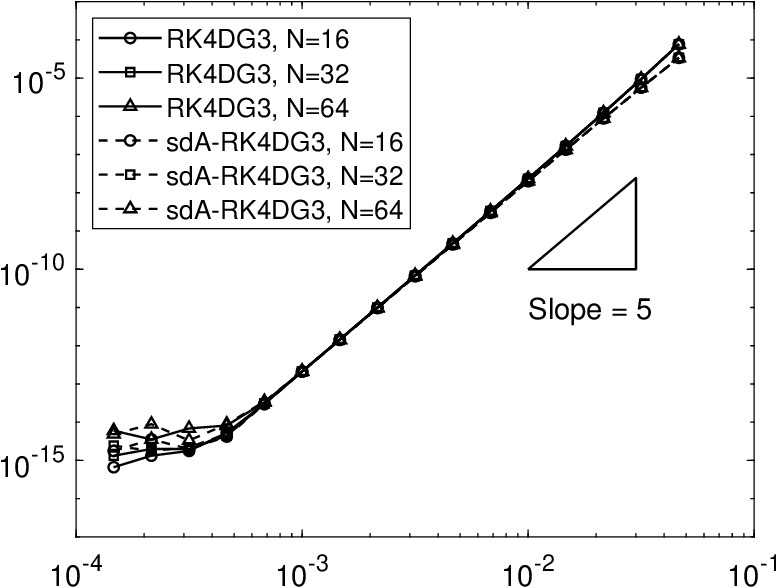}
	\label{fig:1d-RK4DG3}}
	\subfigure[(sdA-)RK2DG2, $m = 1$]{\includegraphics[width=0.32\textwidth]{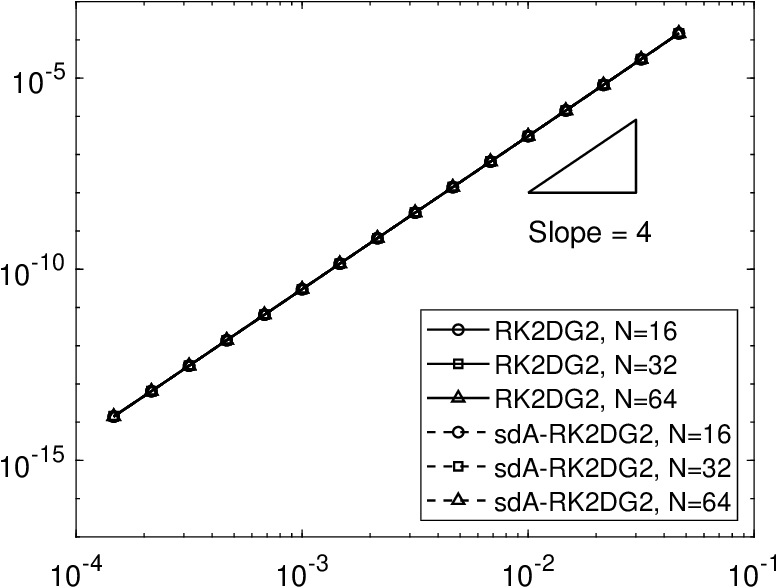}\label{fig:1d-RK2DG2}}
	\subfigure[(sdA-)RK3DG3, $m = 1$]{\includegraphics[width=0.32\textwidth]{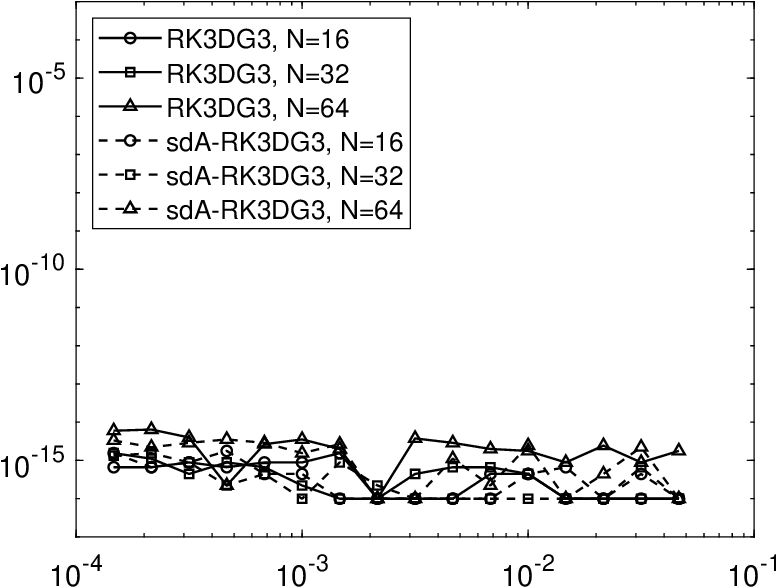}\label{fig:1d-RK3DG3}}
	\subfigure[(sdA-)RK4DG3, $m = 2$]{\includegraphics[width=0.32\textwidth]{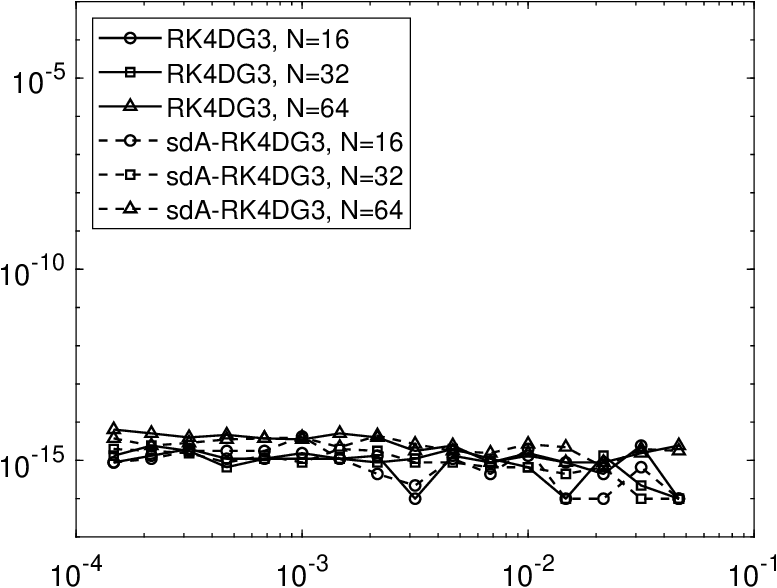}\label{fig:1d-RK4DG3m2}}
	\subfigure[(sdA-)RK5DG4, $m = 1$]{\includegraphics[width=0.32\textwidth]{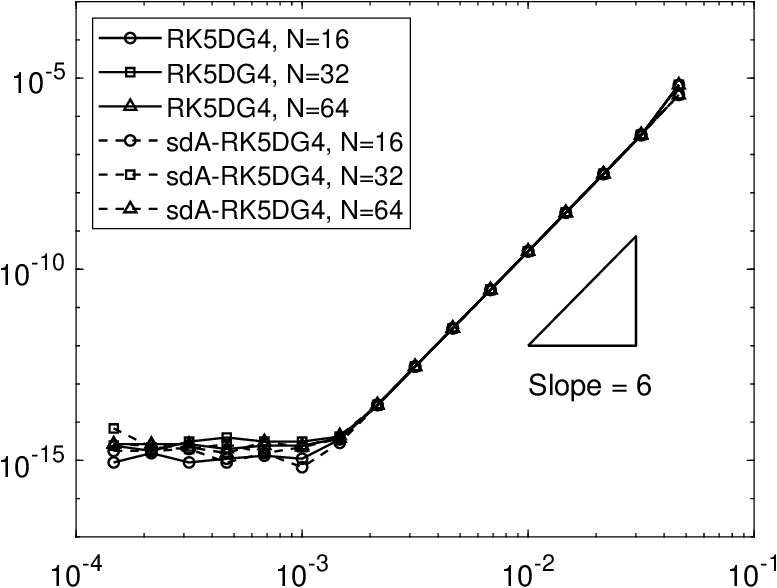}
	\label{fig:1d-RK5DG4}}
	\subfigure[(sdA-)RK5DG1, $m = 1$]{\includegraphics[width=0.32\textwidth]{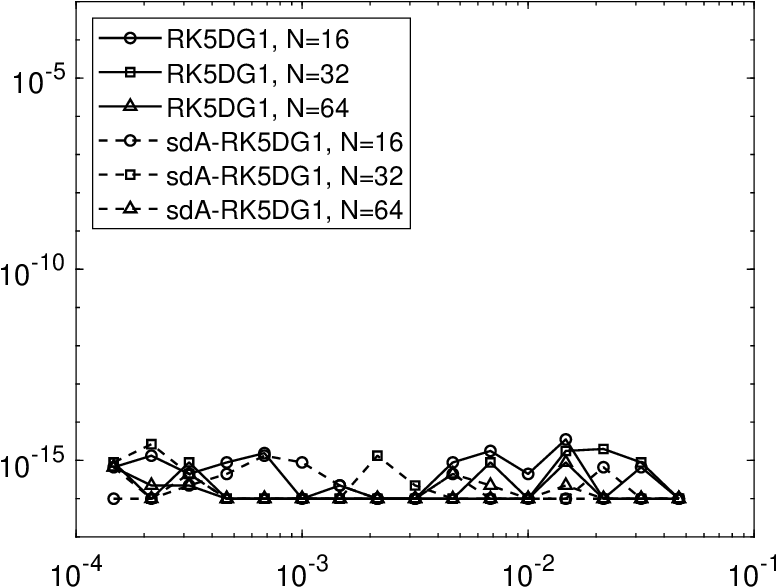}\label{fig:1d-RK5DG1}}
	\subfigure[(sdA-)RK5DG2, $m = 2$]{\includegraphics[width=0.32\textwidth]{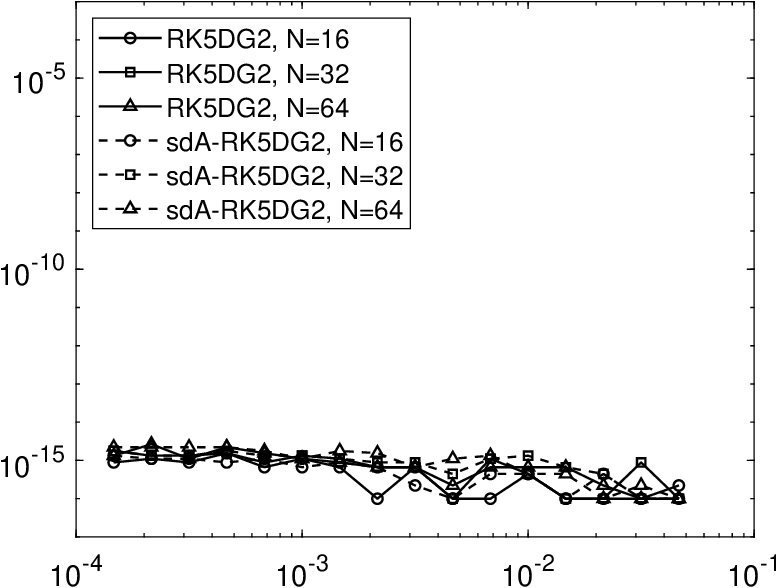}\label{fig:1d-RK5DG2}}
	\caption{Stability verification for 1D RKDG and sdA-RKDG schemes. The $x$ axis is the CFL number. The $y$ axis is the value of $\delta$ defined in \eqref{eq:delta}. }\label{fig:stab-1d}
\end{figure}

\begin{figure}[h!]{\setcounter{subfigure}{0}}
		\centering
		\subfigure[(sdA-)RK2DG1, $m = 1$]{\includegraphics[width=0.32\textwidth]{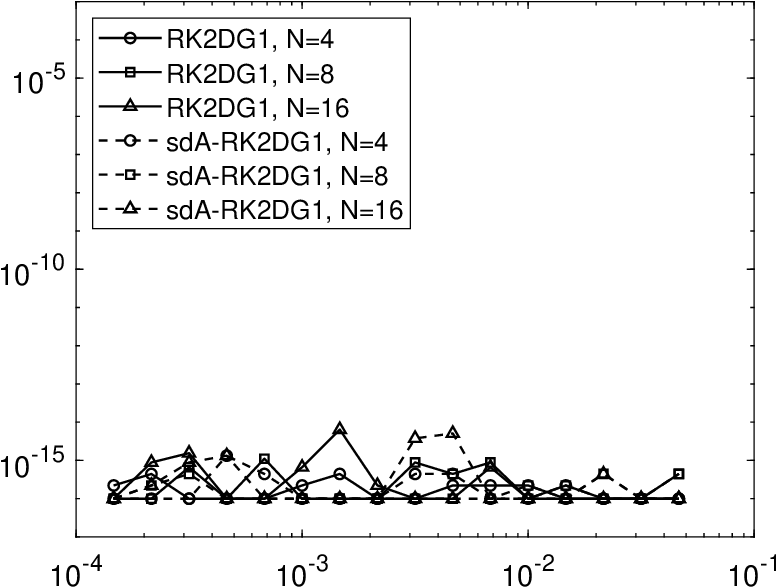}\label{fig:2d-RK2DG1}}
		\subfigure[(sdA-)RK3DG2, $m = 1$]{\includegraphics[width=0.32\textwidth]{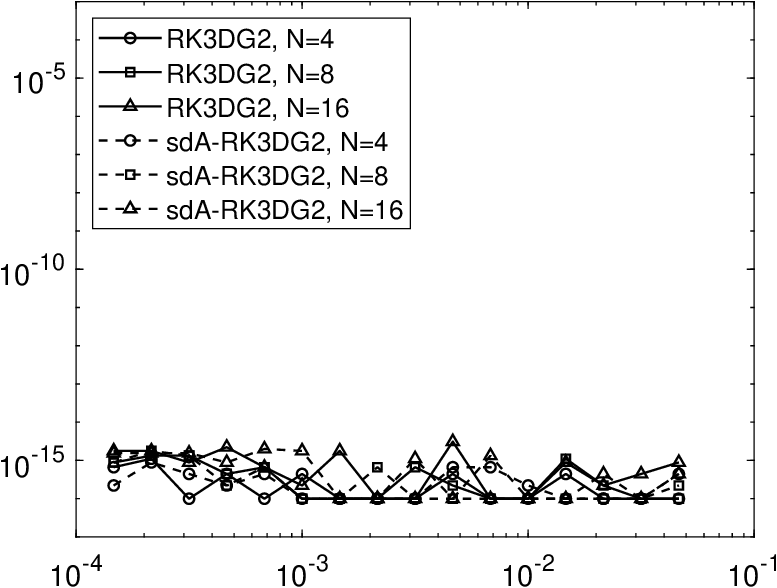}\label{fig:2d-RK3DG2}}
		\subfigure[(sdA-)RK4DG3, $m = 1$]{\includegraphics[width=0.32\textwidth]{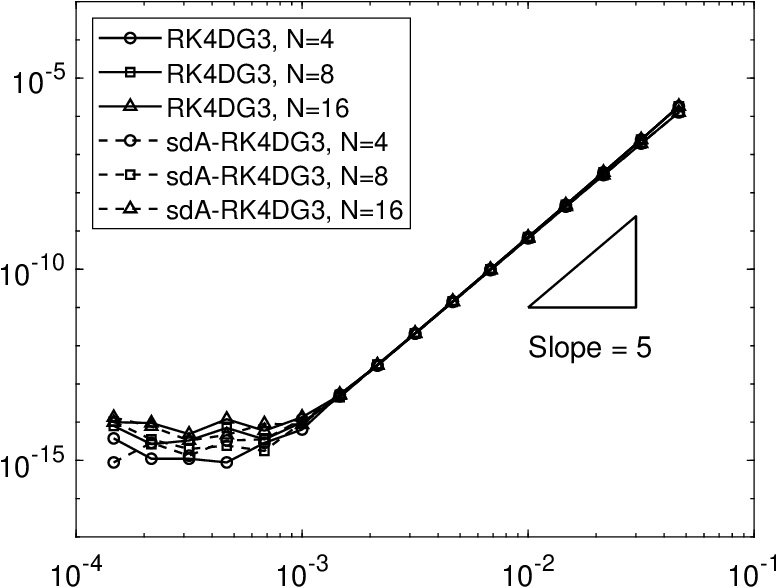}
			\label{fig:2d-RK4DG3}}
		\subfigure[(sdA-)RK2DG2, $m = 1$]{\includegraphics[width=0.32\textwidth]{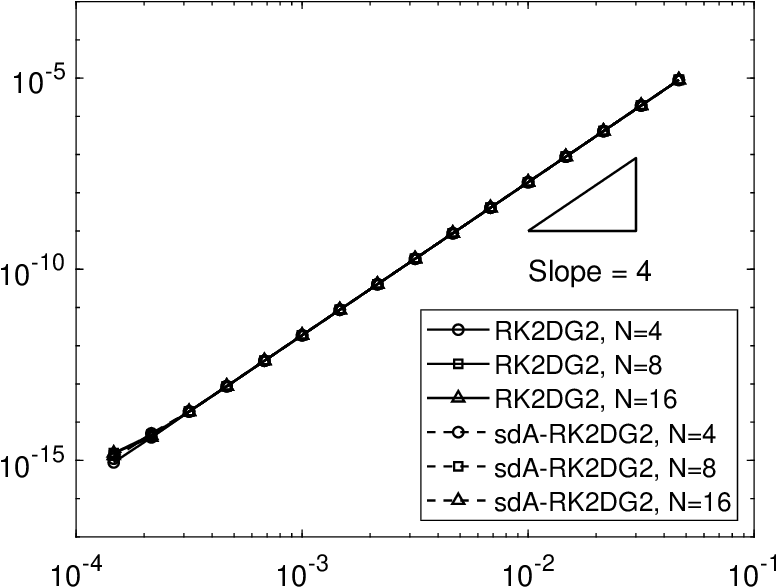}\label{fig:2d-RK2DG2}}
		\subfigure[(sdA-)RK3DG3, $m = 1$]{\includegraphics[width=0.32\textwidth]{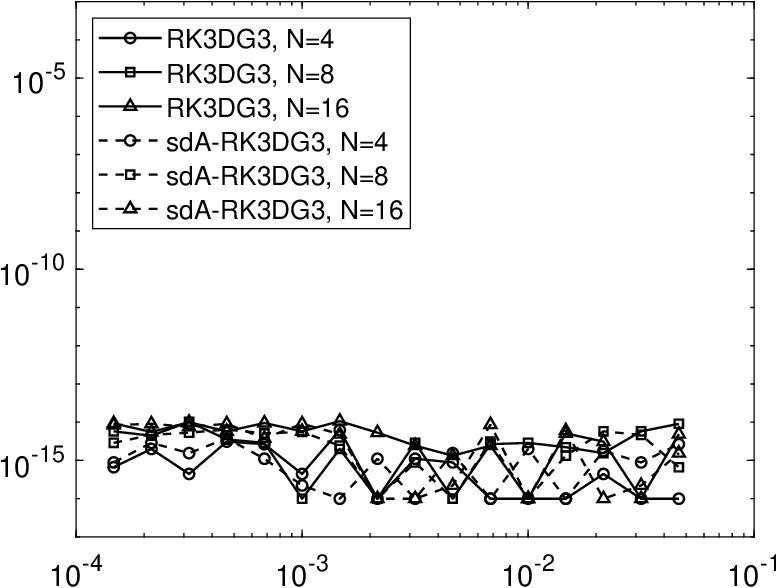}\label{fig:2d-RK3DG3}}
		\subfigure[(sdA-)RK4DG3, $m = 2$]{\includegraphics[width=0.32\textwidth]{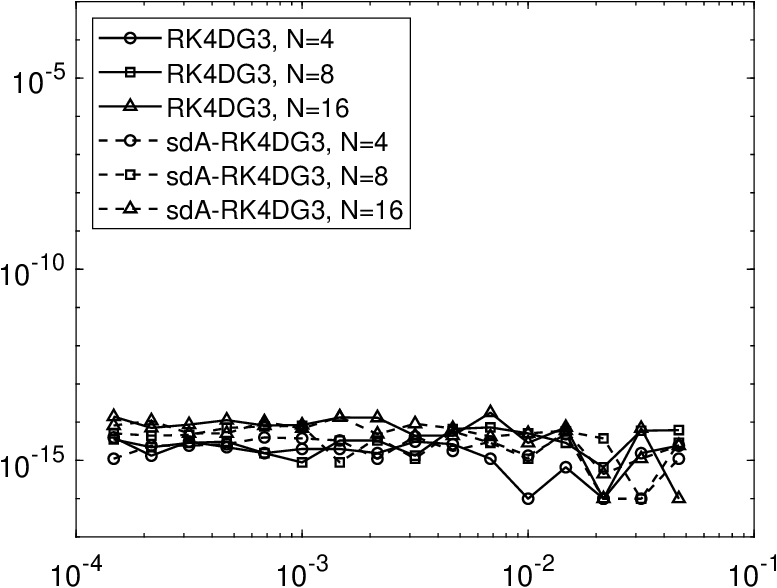}\label{fig:2d-RK4DG3m2}}
		\subfigure[(sdA-)RK5DG4, $m = 1$]{\includegraphics[width=0.32\textwidth]{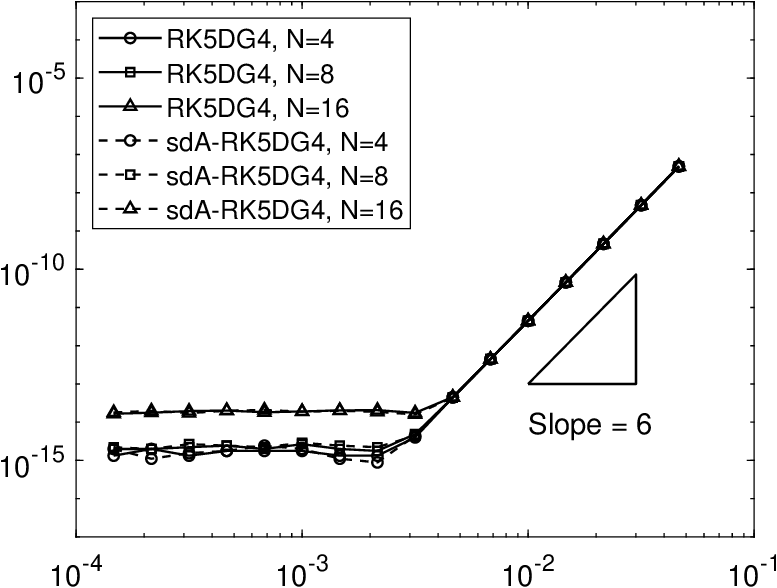}
			\label{fig:2d-RK5DG4}}
		\subfigure[(sdA-)RK5DG1, $m = 1$]{\includegraphics[width=0.32\textwidth]{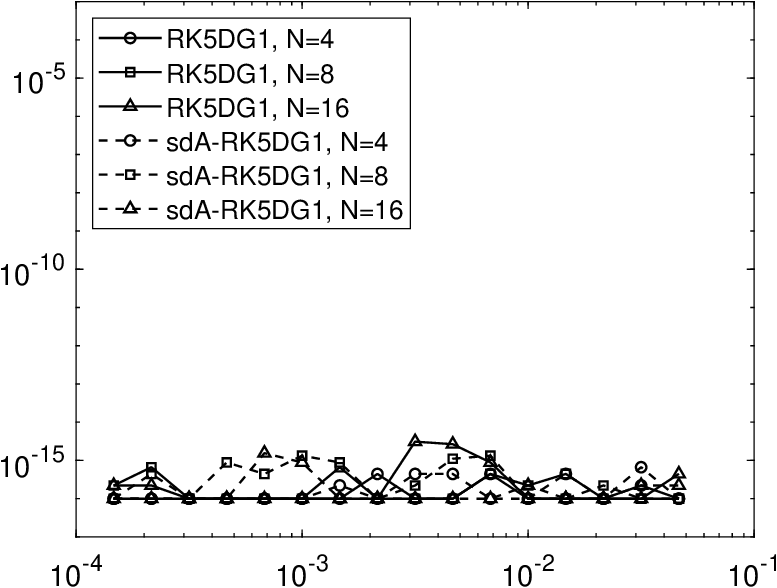}\label{fig:2d-RK5DG1}}
		\subfigure[(sdA-)RK5DG2, $m = 2$]{\includegraphics[width=0.32\textwidth]{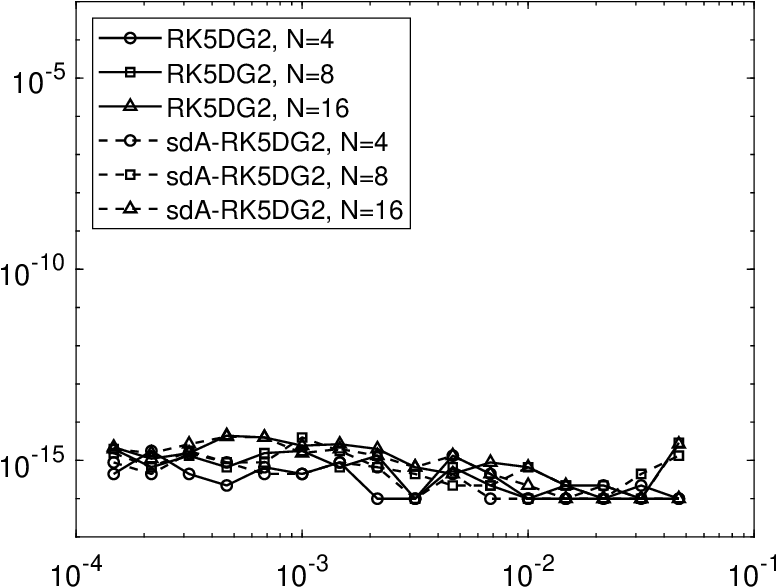}\label{fig:2d-RK5DG2}}
		\caption{Stability verification for 2D RKDG and sdA-RKDG schemes. The $x$ axis is the CFL number. The $y$ axis is value of $\delta$ defined in \eqref{eq:delta}.}\label{fig:stab-2d}
	\end{figure}

	Below, we present plots showing the dependence of $\delta$ on the CFL number for various $N$ values across different (sdA-)RK$r$DG$k$ schemes. The 1D results, with $N = 16, 32, 64$, are shown in Figure \ref{fig:stab-1d}. While the 2D results, with $N = 4, 8, 16$, are shown in Figure \ref{fig:stab-2d}.
	
	Consistent with our theoretical analysis, the stability behavior appears to be independent of dimensionality, with 1D and 2D results exhibiting similar trends. Additionally, when $r$ and $k$ are fixed, both RKDG and sdA-RKDG schemes demonstrate similar patterns of $\delta$ regardless of the number of mesh cells. These results verify different types of stability outlined in Theorems \ref{thm:stabr}, and \ref{thm:stabr-low} (and also in Theorem \ref{thm:key}). See Table \ref{tab:stab} for a summary.

	\begin{table}[h!]
		\centering
		\begin{tabular}{cccccc}
			\hline
			1D Figure&2D Figure&$r$&$k$&Stability&Verification of\\
			\hline
			\ref{fig:1d-RK2DG1}&\ref{fig:2d-RK2DG1}&$2$&$1$&Monotonicity&Theorem \ref{thm:stabr-low}, item 2\\
			\ref{fig:1d-RK2DG2}&\ref{fig:2d-RK2DG2}&$2$&$2$&Weak$(4)$&Theorem \ref{thm:stabr}, item 2\\
			\ref{fig:1d-RK3DG2}, \ref{fig:1d-RK3DG3}&\ref{fig:2d-RK3DG2}, \ref{fig:2d-RK3DG3}&$3$&$2,3$&Monotonicity&Theorem \ref{thm:stabr}, item 3\\
			\ref{fig:1d-RK4DG3}, \ref{fig:1d-RK4DG3m2}&\ref{fig:2d-RK4DG3}, \ref{fig:2d-RK4DG3m2}&$4$&$3$&Strong$(2)$&Theorem \ref{thm:stabr}, item 4\\
			\ref{fig:1d-RK5DG4}&\ref{fig:2d-RK5DG4}&$5$&$4$&Weak$(6)$&Theorem \ref{thm:stabr}, item 1\\
			\ref{fig:1d-RK5DG1}&\ref{fig:2d-RK5DG1}&$5$&$1$&Monotonicity&Theorem \ref{thm:stabr-low}, item 1\\
			\ref{fig:1d-RK5DG2}&\ref{fig:2d-RK5DG2}&$5$&$2$&Strong$(2)$&Theorem \ref{thm:stabr-low}, item 2\\
			\hline
		\end{tabular}
		\caption{Verification of Theorems \ref{thm:key}, \ref{thm:stabr}, and \ref{thm:stabr-low} with Figures \ref{fig:stab-1d} and \ref{fig:stab-2d}.}\label{tab:stab}
	\end{table}

\subsection{Accuracy tests with smooth solutions}

\;\smallskip

	In this test, we consider the linear advection equation with sufficiently smooth solutions. For the 1D problem, we set $u(x,0) = \sin(2\pi x)$, and the exact solution is $u(x,t) = \sin(2\pi(x-t))$. In the 2D problem, we set $u(x,y,0) = \sin(2\pi(x+y))$, and the exact solution is $u(x,y,t) = \sin(2\pi(x+y-2t))$. The final time is set as $t^n = 1$. We solve the problem using RK$r$DG$k$ and sdA-RK$r$DG$k$ methods with $r = k+1 = 2,3,4,5$. The time step size is set as $\tau = 0.1/(dN)$ for $r = 2,3,4$, and $\tau = 0.1/(dN^{6/5})$ for $r = 5$.
	
	The $L^2$ errors of the numerical solutions on 1D uniform meshes, 1D randomly perturbed nonuniform meshes, and 2D uniform meshes are provided in Tables \ref{tab:acc-sin-1d-uni}, \ref{tab:acc-sin-1d-nonuni}, and \ref{tab:acc-sin-2d-uni}, respectively. As observed, in all cases, the sdA-RKDG schemes achieve the optimal convergence rate, with numerical errors comparable to the corresponding RKDG schemes.
\begin{table}[htbp]
	\centering
	\begin{tabular}{c|cc|cc|cc|cc}
		\hline
		RKDG&\multicolumn{2}{c|}{RK2DG1}&\multicolumn{2}{c|}{RK3DG2}&\multicolumn{2}{c|}{RK4DG3}&\multicolumn{2}{c}{RK5DG4}\\
		\hline
		$N$& $L^2$ error&Order&$L^2$ error&Order&$L^2$ error&Order&$L^2$ error&Order\\
		\hline
		20 & 6.90e-03 & - & 5.67e-04 & - & 3.46e-05 & - & 1.71e-06 & - \\  
		40 & 1.73e-03 & 2.00 & 7.12e-05 & 2.99 & 2.17e-06 & 4.00 & 5.67e-08 & 4.91 \\ 
		80 & 4.37e-04 & 1.98 & 8.91e-06 & 3.00 & 1.35e-07 & 4.00 & 1.62e-09 & 5.13 \\ 
		160& 1.10e-04 & 1.99 & 1.11e-06 & 3.00 & 8.46e-09 & 4.00 & 5.07e-11 & 5.00 \\ 	
		320& 2.77e-05 & 1.99 & 1.39e-07 & 3.00 & 5.29e-10 & 4.00 & 1.58e-12 & 5.00 \\ 
		\hline
		sdA-RKDG&\multicolumn{2}{c|}{sdA-RK2DG1}&\multicolumn{2}{c|}{sdA-RK3DG2}&\multicolumn{2}{c|}{sdA-RK4DG3}&\multicolumn{2}{c}{sdA-RK5DG4}\\
		\hline
		$N$& $L^2$ error&Order&$L^2$ error&Order&$L^2$ error&Order&$L^2$ error&Order\\
		\hline
		20 & 8.23e-03 & - & 7.67e-04 & - & 4.98e-05 & - & 2.10e-06 & - \\  
		40 & 2.10e-03 & 1.97 & 9.62e-05 & 3.00 & 3.12e-06 & 4.00 & 7.03e-08 & 4.90 \\ 
		80 & 5.34e-04 & 1.98 & 1.20e-05 & 3.00 & 1.95e-07 & 4.00 & 2.12e-09 & 5.05 \\ 
		160& 1.35e-04 & 1.99 & 1.51e-06 & 3.00 & 1.22e-08 & 4.00 & 6.03e-11 & 5.13 \\ 
		320& 3.38e-05 & 1.99 & 1.88e-07 & 3.00 & 7.63e-10 & 4.00 & 1.85e-12 & 5.03 \\ 
		\hline 
	\end{tabular}
	\caption{Error table of the RKDG and sdA-RKDG methods for 1D linear advection equation on uniform meshes. $u(x,0) = \sin(2\pi x)$.  $t^n = 1$. }
	\label{tab:acc-sin-1d-uni}
\end{table}

\begin{table}[htbp]
	\centering
	\begin{tabular}{c|cc|cc|cc|cc}
		\hline
		RKDG&\multicolumn{2}{c|}{RK2DG1}&\multicolumn{2}{c|}{RK3DG2}&\multicolumn{2}{c|}{RK4DG3}&\multicolumn{2}{c}{RK5DG4}\\
		\hline
		$N$& $L^2$ error&Order&$L^2$ error&Order&$L^2$ error&Order&$L^2$ error&Order\\
		\hline
		20 & 7.12e-03 & - & 6.54e-04 & - & 4.16e-05 & - & 2.43e-06 & - \\ 
		40 & 1.84e-03 & 1.95 & 8.06e-05 & 3.02 & 2.65e-06 & 3.97 & 7.51e-08 & 5.02 \\ 
		80 & 4.68e-04 & 1.97 & 1.03e-05 & 2.97 & 1.70e-07 & 3.96 & 1.99e-09 & 5.24 \\ 
		160& 1.16e-04 & 2.01 & 1.32e-06 & 2.97 & 1.05e-08 & 4.01 & 7.72e-11 & 4.69 \\ 
		320& 2.96e-05 & 1.97 & 1.60e-07 & 3.04 & 6.55e-10 & 4.01 & 2.23e-12 & 5.11 \\ 
		\hline
		sdA-RKDG&\multicolumn{2}{c|}{sdA-RK2DG1}&\multicolumn{2}{c|}{sdA-RK3DG2}&\multicolumn{2}{c|}{sdA-RK4DG3}&\multicolumn{2}{c}{sdA-RK5DG4}\\
		\hline
		$N$& $L^2$ error&Order&$L^2$ error&Order&$L^2$ error&Order&$L^2$ error&Order\\
		\hline
		20 & 8.12e-03 & - & 7.91e-04 & - & 5.64e-05 & - & 2.97e-06 & - \\ 
		40 & 2.09e-03 & 1.96 & 1.01e-04 & 2.97 & 3.67e-06 & 3.94 & 8.72e-08 & 5.09 \\ 
		80 & 5.33e-04 & 1.97 & 1.28e-05 & 2.98 & 2.11e-07 & 4.12 & 2.75e-09 & 4.99 \\ 
		160& 1.34e-04 & 1.99 & 1.58e-06 & 3.02 & 1.31e-08 & 4.01 & 7.92e-11 & 5.12 \\ 
		320& 3.33e-05 & 2.01 & 2.01e-07 & 2.98 & 8.77e-10 & 3.90 & 2.48e-12 & 5.00 \\ 
		\hline 
	\end{tabular}
	\caption{Error table of the RKDG and sdA-RKDG methods for 1D linear advection equation on randomly perturbed nonuniform meshes. Each node is randomly perturbed by $15\%$ of the mesh size. $u(x,0) = \sin(2\pi x)$. $t^n = 1$.}
	\label{tab:acc-sin-1d-nonuni}
\end{table}

\begin{table}[htbp]
	\centering
	\begin{tabular}{c|cc|cc|cc|cc}
		\hline
		RKDG&\multicolumn{2}{c|}{RK2DG1}&\multicolumn{2}{c|}{RK3DG2}&\multicolumn{2}{c|}{RK4DG3}&\multicolumn{2}{c}{RK5DG4}\\
		\hline
		$N$& $L^2$ error&Order&$L^2$ error&Order&$L^2$ error&Order&$L^2$ error&Order\\
		\hline
		20 & 2.54e-02 & - & 4.06e-03 & - & 4.36e-04 & - & 3.88e-05 & - \\  
		40 & 6.98e-03 & 1.87 & 5.14e-04 & 2.98 & 2.74e-05 & 3.99 & 1.23e-06 & 4.98 \\   
		80 & 1.87e-03 & 1.90 & 6.45e-05 & 3.00 & 1.72e-06 & 4.00 & 3.82e-08 & 5.00 \\   
		160& 4.86e-04 & 1.94 & 8.06e-06 & 3.00 & 1.07e-07 & 4.00 & 1.19e-09 & 5.00 \\   
		320& 1.24e-04 & 1.97 & 1.01e-06 & 3.00 & 6.72e-09 & 4.00 & 3.73e-11 & 5.00 \\   
		\hline
		sdA-RKDG&\multicolumn{2}{c|}{sdA-RK2DG1}&\multicolumn{2}{c|}{sdA-RK3DG2}&\multicolumn{2}{c|}{sdA-RK4DG3}&\multicolumn{2}{c}{sdA-RK5DG4}\\
		\hline
		$N$& $L^2$ error&Order&$L^2$ error&Order&$L^2$ error&Order&$L^2$ error&Order\\
		\hline
		20 & 2.83e-02 & - & 4.72e-03 & - & 5.28e-04 & - & 4.37e-05 & - \\    
		40 & 7.88e-03 & 1.85 & 5.97e-04 & 2.98 & 3.32e-05 & 3.99 & 1.37e-06 & 4.99 \\   
		80 & 2.10e-03 & 1.91 & 7.48e-05 & 3.00 & 2.08e-06 & 4.00 & 4.20e-08 & 5.03 \\   
		160& 5.41e-04 & 1.95 & 9.36e-06 & 3.00 & 1.30e-07 & 4.00 & 1.30e-09 & 5.02 \\ 
		320& 1.38e-04 & 1.98 & 1.17e-06 & 3.00 & 8.13e-09 & 4.00 & 4.02e-11 & 5.02 \\  
		\hline 
	\end{tabular}
	\caption{Error table of the RKDG and sdA-RKDG methods for 2D linear advection equation on uniform meshes. $u(x,y,0) = \sin(2\pi(x+y))$. $t^n = 1$.}
	\label{tab:acc-sin-2d-uni}
\end{table}

\subsection{Accuracy tests with solutions of limited regularity}

\;\smallskip

	When $r = k+1$, Theorem \ref{thm:err-g} indicates that the optimal $r$th-order convergence rate for RK$r$DG$k$ schemes can only be attained if $\flat \geq r+1 = k+2$. To verify the sharpness of this regularity number, we consider the linear advection equation with the initial condition
	\begin{equation}
		u(x,0) = \left(\sin(2\pi x)\right)^{\flat -\frac{1}{3}}
	\end{equation}
	for the 1D problem.  Note that the solution is in $C^{\flat-1}(\overline{\Omega})\cap H^\flat(\Omega_h)$ but not in $C^{\flat}(\overline{\Omega})\cap H^{\flat+1}(\Omega_h)$ \cite{xu2020error,xu2024sdf}. When we set $\flat = r = k+1$, observations from Tables \ref{tab:1dacc-reg23} and \ref{tab:1dacc-reg45} show that both the RK$r$DG$k$ and the sdA-RK$r$DG$k$ methods exhibit suboptimal convergence for $r = k+1 = 2,3,4,5$.  Here we have chosen a large final time to make the order degeneracy observable for $r = k+1 = 4,5$. Conversely, under the same settings with $\flat = r + 1 = k+2$, all schemes achieve optimal convergence rates. A similar test has been conducted for the 2D case with the initial condition $u(x,y,0) = \left(\sin(2\pi (x+y))\right)^{\flat -\frac{1}{3}}$, where $r = k + 1 = 2,3$. See Table \ref{tab:2dacc-reg23}.

\begin{table}[htbp]
	\centering
	\begin{tabular}{c|cc|cc|cc|cc}
		\hline
		$\flat$&\multicolumn{2}{c|}{$\flat = r$}&\multicolumn{2}{c|}{$\flat = r+1$}&\multicolumn{2}{c|}{$\flat = r$}&\multicolumn{2}{c}{$\flat = r+1$}\\
		\hline
		RKDG&\multicolumn{4}{c|}{RK2DG1}&\multicolumn{4}{c}{RK3DG2}\\
		\hline
		$N$& $L^2$ error&Order&$L^2$ error&Order&$L^2$ error&Order&$L^2$ error&Order\\
		\hline
		80  & 2.25e-03 & - & 1.16e-03 & -& 5.93e-05 & - & 7.22e-05 & -\\ 
		160 & 7.14e-04 & 1.66 & 2.70e-04 & 2.10& 8.24e-06 & 2.85 & 9.07e-06 & 2.99\\ 
		320 & 2.30e-04 & 1.63 & 6.59e-05 & 2.04& 1.18e-06 & 2.81 & 1.13e-06 & 3.00\\ 
		640 & 7.49e-05 & 1.62 & 1.63e-05 & 2.02& 1.75e-07 & 2.76 & 1.42e-07 & 3.00 \\ 
		1280& 2.45e-05 & 1.61 & 4.84e-06 & 2.00 & 2.68e-08 & 2.70 & 1.77e-08 & 3.00 \\ 		
		\hline
		sdA-RKDG&\multicolumn{4}{c|}{sdA-RK2DG1}&\multicolumn{4}{c}{sdA-RK3DG2}\\
		\hline
		$N$& $L^2$ error&Order&$L^2$ error&Order&$L^2$ error&Order&$L^2$ error&Order\\
		\hline
		80  & 2.00e-03 & - & 1.24e-03 & -& 7.36e-05 & - & 9.76e-05 & -\\ 
		160 & 6.35e-04 & 1.66 & 3.07e-04 & 2.01 & 9.73e-06 & 2.92 & 1.22e-05 & 2.99\\ 
		320 & 2.05e-04 & 1.63 & 7.71e-05 & 2.00 & 1.30e-06 & 2.90 & 1.53e-06 & 3.00\\ 
		640 & 6.66e-05 & 1.62 & 1.93e-05 & 2.00 & 1.79e-07 & 2.86 & 1.92e-07 & 3.00\\  
		1280& 2.19e-05 & 1.61 & 4.84e-06 & 2.00 & 2.56e-08 & 2.81 & 2.40e-08 & 3.00\\    
		\hline 
	\end{tabular}
	\caption{Error table of the RK$r$DG$k$ and the sdA-RK$r$DG$k$ methods for 1D linear advection equation on uniform meshes. $u(x,0) = (\sin(2\pi x))^{\flat-1/3}$. $r = k+1 = 2,3$. $t^n = 1$. }
	\label{tab:1dacc-reg23}
\end{table}

\begin{table}[htbp]
	\centering
	\begin{tabular}{c|cc|cc|cc|cc}
		\hline
		$\flat$&\multicolumn{2}{c|}{$\flat = r$}&\multicolumn{2}{c|}{$\flat = r+1$}&\multicolumn{2}{c|}{$\flat = r$}&\multicolumn{2}{c}{$\flat = r+1$}\\
		\hline
		RKDG&\multicolumn{4}{c|}{RK4DG3}&\multicolumn{4}{c}{RK5DG4}\\
		\hline
		$N$& $L^2$ error&Order&$L^2$ error&Order&$L^2$ error&Order&$L^2$ error&Order\\
		\hline
		20 & 4.47e-03 & - & 7.23e-03 & - & 3.69e-04 & - & 5.43e-04 & - \\ 
		40 & 2.45e-04 & 4.19 & 2.72e-04 & 4.73& 1.10e-05 & 5.07 & 1.43e-05 & 5.25 \\   
		80 & 1.75e-05 & 3.81 & 1.50e-05 & 4.18& 3.65e-07 & 4.91 & 4.43e-07 & 5.01 \\  
		160& 1.38e-06 & 3.66 & 9.03e-07 & 4.05& 1.26e-08 & 4.86 & 1.39e-08 & 5.00 \\   
		320& 1.17e-07 & 3.56 & 5.57e-08 & 4.02& 4.47e-10 & 4.82 & 4.33e-10 & 5.00 \\   
		\hline
		sdA-RKDG&\multicolumn{4}{c|}{sdA-RK4DG3}&\multicolumn{4}{c}{sdA-RK5DG4}\\
		\hline
		$N$& $L^2$ error&Order&$L^2$ error&Order&$L^2$ error&Order&$L^2$ error&Order\\
		\hline
		20 & 3.55e-03 & - & 5.68e-03 & -& 3.97e-04 & - & 5.94e-04 & - \\  
		40 & 2.19e-04 & 4.02 & 2.80e-04 & 4.34& 1.24e-05 & 5.00 & 1.76e-05 & 5.08 \\   
		80 & 1.58e-05 & 3.79 & 1.62e-05 & 4.11& 3.97e-07 & 4.96 & 5.37e-07 & 5.04 \\  
		160& 1.27e-06 & 3.64 & 9.80e-07 & 4.05& 1.31e-08 & 4.92 & 1.65e-08 & 5.03 \\ 
		320& 1.10e-07 & 3.52 & 6.04e-08 & 4.02& 4.49e-10 & 4.87 & 5.07e-10 & 5.02 \\   
		\hline 
	\end{tabular}
	\caption{Error table of the RK$r$DG$k$ and the sdA-RK$r$DG$k$ methods for 1D linear advection equation on uniform meshes. $u(x,0) = (\sin(2\pi x))^{\flat-1/3}$. $r = k+1 = 4,5$. $t^n = 500$. }
	\label{tab:1dacc-reg45}
\end{table}

\begin{table}[htbp]
	\centering
	\begin{tabular}{c|cc|cc|cc|cc}
		\hline
		$\flat$&\multicolumn{2}{c|}{$\flat = r$}&\multicolumn{2}{c|}{$\flat = r+1$}&\multicolumn{2}{c|}{$\flat = r$}&\multicolumn{2}{c}{$\flat = r+1$}\\
		\hline
		RKDG&\multicolumn{4}{c|}{RK2DG1}&\multicolumn{4}{c}{RK3DG2}\\
		\hline
		$N$& $L^2$ error&Order&$L^2$ error&Order&$L^2$ error&Order&$L^2$ error&Order\\
		\hline
		20 & 7.28e-02 & - & 9.18e-02 & - & 1.79e-02 & - & 2.44e-02 & - \\ 
		40 & 2.20e-02 & 1.73 & 2.12e-02 & 2.11 & 2.76e-03 & 2.70 & 3.88e-03 & 2.65 \\   
		80 & 6.52e-03 & 1.75 & 4.55e-03 & 2.22 & 3.88e-04 & 2.83 & 5.17e-04 & 2.91 \\  
		160& 2.06e-03 & 1.66 & 1.07e-03 & 2.08 & 5.34e-05 & 2.86 & 6.55e-05 & 2.98 \\ 
		320& 6.66e-04 & 1.63 & 2.72e-04 & 1.98 & 7.41e-06 & 2.85 & 8.21e-06 & 3.00 \\ 
		\hline
		sdA-RKDG&\multicolumn{4}{c|}{sdA-RK2DG1}&\multicolumn{4}{c}{sdA-RK3DG2}\\
		\hline
		$N$& $L^2$ error&Order&$L^2$ error&Order&$L^2$ error&Order&$L^2$ error&Order\\
		\hline
		 20 & 6.89e-02 & - & 8.15e-02 & -& 2.05e-02 & - & 2.88e-02 & - \\   
		 40 & 2.03e-02 & 1.77 & 1.98e-02 & 2.04& 3.15e-03 & 2.70 & 4.54e-03 & 2.67 \\  
		 80 & 6.09e-03 & 1.74 & 4.56e-03 & 2.12& 4.36e-04 & 2.85 & 6.00e-04 & 2.92 \\ 
		 160& 1.93e-03 & 1.66 & 1.15e-03 & 1.99& 5.88e-05 & 2.89 & 7.60e-05 & 2.98 \\  
		 320& 6.22e-04 & 1.63 & 2.97e-04 & 1.95& 7.94e-06 & 2.89 & 9.53e-06 & 3.00 \\ 
		\hline 
	\end{tabular}
	\caption{Error table of the RK$r$DG$k$ and the sdA-RK$r$DG$k$ methods for 2D linear advection equation on uniform meshes. $u(x,y,0) = (\sin(2\pi (x+y)))^{\flat-1/3}$. $r = k+1 = 2,3$. $t^n = 1$.}
	\label{tab:2dacc-reg23}
\end{table}



\section{Proofs of technical lemmas}\label{sec:proofs}

\subsection{Proof of Lemma \ref{lem:PiL}}\label{sec:proof-lem:PiL}

\;\smallskip

We will consider the 1D and 2D cases separately. \smallskip

\underline{\em 
	1D case.} Using the self-adjoint property of $\Pip$, the definition of $L_h$, and the integration by parts, one can get
\begin{equation}
	\begin{aligned}
		\ip{\Pip L_h w_h}{v_h}
		=& \ip{L_h w_h}{\Pip v_h} \\
		=& \beta\ip{w_h}{(\Pip v_h)_x} -\beta \sum_{i=1}^N \left(({w}_{h}^-(\Pip v_{h})^-)_{i+\hf} - ({w}_{h}^-(\Pip v_{h})^+)_{i-\hf}\right)\\
		=& -\beta\ip{(w_h)_x}{\Pip v_h} - \beta \sum_{i=1}^N  ([w_{h}](\Pip v_{h})^+)_{i-\hf}.
	\end{aligned}
\end{equation}
Note that $(w_h)_x \in V_h^{k-1}$ and $\Pip v_h \perp V_h^{k-1}$. Thus, $\ip{(w_h)_x}{\Pip v_h} = 0$ and 
\begin{equation}\label{eq:PipLwv}
	\ip{\Pip L_h w_h}{v_h} = - \beta \sum_{i=1}^N  ([w_{h}](\Pip v_{h})^+)_{i-\hf} \leq C\jpnm{w_h}\nmG{\Pip v_h}.
\end{equation}
Here we have used the Cauchy--Schwarz inequality. With the inverse estimate \eqref{eq:inv1}, we get 
\begin{equation}\label{eq:est-nmG-Pipv}
	\nmG{\Pip v_h}\leq Ch^{-1/2}\nm{\Pip v_h} \leq Ch^{-1/2}\nm{v_h}.
\end{equation} 
Substituting \eqref{eq:est-nmG-Pipv} into \eqref{eq:PipLwv} completes the proof of \eqref{eq:key1} in the 1D case. 

For the proof of \eqref{eq:key2}, note we have \cite{zhang2010stability,sun2017stability}
\begin{equation}\label{eq:Ltop}
	\ip{L_h^\top w_h}{v_h} = -\beta\ip{w_h}{(v_h)_x} +\beta \sum_{i=1}^N \left((w_h^+ v_h^-)_{i+\hf} - (w_h^+ v_h^+)_{i-\hf}\right). 
\end{equation} 
Comparing \eqref{eq:Ltop} to \eqref{eq:Lh1d}, $L_h^\top$ can be defined by replacing the numerical flux $w_h^-$ in \eqref{eq:Lh1d} by $w_h^+$ and changing the sign of each term. One can take advantage of this similar structure to prove \eqref{eq:key2} by following similar lines as those in the proof of \eqref{eq:key1}. 
\begin{equation}
\begin{aligned}
&\ip{\Pip L_h^\top w_h}{v_h}
= \ip{L_h^\top w_h}{\Pip v_h} \\
=& -\beta\ip{w_h}{(\Pip v_h)_x} +\beta \sum_{i=1}^N \left(({w}_{h}^+(\Pip v_{h})^-)_{i+\hf} - ({w}_{h}^+(\Pip v_{h})^+)_{i-\hf}\right)\\
=& \beta\ip{(w_h)_x}{\Pip v_h} + \beta \sum_{i=1}^N  ([w_{h}](\Pip v_{h})^-)_{i+\hf}
= \beta \sum_{i=1}^N  ([w_{h}](\Pip v_{h})^-)_{i+\hf}\\
\leq& C \jpnm{w_h}\nmG{\Pip v_h}\leq Ch^{-\hf} \jpnm{w_h}\nm{\Pip v_h}\leq Ch^{-\hf} \jpnm{w_h}\nm{v_h}.
\end{aligned}
\end{equation}
\smallskip

\underline{\em 2D case.} Similar to the 1D case, we use the self-adjoint property of $\Pip$, the definition of $L_h$, and the integration by parts. It yields

\begin{equation}\label{eq:stabkeylem-2D}
	\begin{aligned}
		&\ip{\Pip L_h w_h}{v_h}
		= \ip{L_h w_h}{\Pip v_h}\\
		=& \ip{w_h}{\beta^x(\Pip v_h)_x+\beta^y(\Pip v_h)_y}\\
		&-\sum_{i=1}^{N_x}\sum_{j=1}^{N_y}\int_{x_{i-\hf}}^{x_{i+\hf}} \beta^y \left(w_h\left(x,y_{j+\hf}^-\right) (\Pip v_h)\left(x,y_{j+\hf}^-\right)-w_h\left(x,y_{j-\hf}^-\right) (\Pip  v_h)\left(x,y_{j-\hf}^+\right)\right)\dd x\\
		&-\sum_{i=1}^{N_x}\sum_{j=1}^{N_y}\int_{y_{j-\hf}}^{y_{j+\hf}}  \beta^x \left(w_h\left(x_{i+\hf}^-,y\right) (\Pip v_h)\left(x_{i+\hf}^-,y\right)-w_h\left(x_{i-\hf}^-,y\right) (\Pip v_h)\left(x_{i-\hf}^+,y\right)\right)\dd y\\
		=& - \ip{\beta^x(w_h)_x + \beta^y(w_h)_y}{\Pip v_h} - \sum_{i=1}^{N_x}\sum_{j=1}^{N_y}\int_{x_{i-\hf}}^{x_{i+\hf}} \beta^y [w_h]_{i,j-\hf}(\Pip  v_h)\left(x,y_{j-\hf}^+\right)\dd x\\
		&-\sum_{i=1}^{N_x}\sum_{j=1}^{N_y}\int_{y_{j-\hf}}^{y_{j+\hf}}  \beta^x [w_h]_{i-\hf,j} (\Pip v_h)\left(x_{i-\hf}^+,y\right)\dd y.
	\end{aligned}
\end{equation}
Note that the first term is $0$ since $(\beta^x\partial_x + \beta^y\partial_y)w_h\in V_h^{k-1}$ and $\Pip v_h \perp V_h^{k-1}$. As we have done for the 1D case, we can then prove the lemma by applying the Cauchy--Schwarz inequality, the inverse estimate, and the fact $\nm{\Pip v_h}\leq \nm{v_h}$.  The 2D case of \eqref{eq:key2} can be proved similarly by rewriting $L_h^\top$ in a similar form as \eqref{eq:Ltop}. Details are omitted.

\subsection{Proof of Lemma \ref{lem:Luiwv}}\label{sec:lem:Luiwv}

\;\smallskip

To prove Lemma \ref{lem:Luiwv}, we need a few preparatory lemmas.

\begin{LEM}[Discrete integration by parts]\label{lem:ibp}
	\begin{equation}
		\ip{L_h^i w_h}{v_h} = \sum_{j=0}^{i-1} (-1)^{j+1}\jpip{L_h^{i-j-1} w_h}{L_h^{j} v_h} +(-1)^i\ip{w_h}{L_h^iv_h}.
	\end{equation}	
\end{LEM}
\begin{proof}
	Recall Proposition \ref{prop:sn}. The lemma can be proved by inductively applying $\ip{L_h w_h}{v_h} = -\jpip{w_h}{v_h} -\ip{w_h}{L_h v_h}$. Details are omitted. 	
\end{proof}

\begin{LEM}\label{lem:LPLjp-general}
	\begin{equation}\label{eq:LPLjp-general}
		\jpnm{L_h^i\Pip L_h w_h} \leq C h^{-i-1} \jpnm{w_h}\qquad \forall i\geq 0.
	\end{equation}
\end{LEM}
\begin{proof}
	With the inverse estimate $\jpnm{v_h}\leq Ch^{-\hf}\nm{v_h}$ and the property $\nm{L_h}\leq C h^{-1}$, we have
	\begin{equation}\label{eq:LPLjp-general-p}
		\jpnm{L_h^i\Pip L_h w_h} \leq Ch^{-\hf}\nm{L_h^i\Pip L_h w_h} \leq Ch^{-\hf}\nm{L_h}^i\nm{\Pip L_h w_h}\leq Ch^{-i-\hf} \nm{\Pip L_hw_h}.
	\end{equation}
	Note from \eqref{eq:key1}, we have $\nm{\Pip L_hw_h} \leq Ch^{-1/2}\jpnm{w_h}$. Applying this estimate to \eqref{eq:LPLjp-general-p} completes the proof of \eqref{eq:LPLjp-general}.  
\end{proof}

Now we prove Lemma \ref{lem:Luiwv}.

\begin{proof}[Proof of Lemma \ref{lem:Luiwv}]
	We define $z_h$ such that $L_h^{\ui}w_h = L_h^{i_1}\Pip L_h z_h$. In other words, 
	\begin{equation}\label{eq:def-zh}
		\begin{aligned}
			z_h =& L_h^{i_2-1}\Pip L_h^{i_3}\cdots \Pip L_h^{i_{|\ui|}}w_h \\
			=& \left(L_h^{i_2-1}\Pip L_h\right)
			\left(L_h^{i_3-1}\Pip L_h\right)\cdots 
			\left(L_h^{i_{|\ui|-1}-1}\Pip L_h \right) L_h^{i_{|\ui|}-1}w_h.	
		\end{aligned}
	\end{equation}
	First, with Lemma \ref{lem:ibp}, the triangle inequality, and the Cauchy--Schwarz inequality, we have
	\begin{equation}
		\begin{aligned}
			&\left|\ip{L_h^{\ui} w_h}{v_h}\right|\\
			 =& 	\left|\ip{L_h^{i_1-1}(L_h\Pip L_h z_h)}{v_h}\right|\\
			=& \left|\sum_{j=0}^{i_1-2} (-1)^{j+1}\jpip{L_h^{i_1-j-1}\Pip L_h z_h}{L_h^{j} v_h} +(-1)^{i_1-1}\ip{L_h\Pip L_h z_h}{L_h^{i_1-1}v_h}\right|\\
			\leq &\sum_{j=0}^{i_1-2}\jpnm{L_h^{i_1-j-1}\Pip L_h z_h}\jpnm{L_h^{j} v_h}+|\ip{L_h\Pip L_h z_h}{L_h^{i_1-1}v_h}|.
		\end{aligned}
	\end{equation}
	Then we apply Lemmas \ref{lem:LPLjp-general} and \ref{lem:LPL} to get
	\begin{equation}\label{eq:Luiwv}
		\begin{aligned}
			\left|\ip{L_h^{\ui} w_h}{v_h}\right|
			\leq &C\sum_{j=0}^{i_1-2}h^{-(i_1-j)}\jpnm{z_h}\jpnm{L_h^{j} v_h} +Ch^{-1}\jpnm{ z_h}\jpnm{L_h^{i_1-1}v_h}\\
			=&C\sum_{j=0}^{i_1-1}h^{-(i_1-j)}\jpnm{z_h}\jpnm{L_h^{j} v_h} .
		\end{aligned}
	\end{equation}	
	Next, recall the definition of $z_h$ in \eqref{eq:def-zh} and apply Lemma \ref{lem:LPLjp-general} for $[\ui]-2$ times. We have
	\begin{equation}\label{eq:est-zh}
		\jpnm{z_h} \leq C h^{-i_2-i_3-\cdots - i_{[\ui]-1}}\jpnm{L_h^{i_{[\ui]}-1}w_h}\leq Ch^{-|\ui|+i_1+i_{[\ui]}} \jpnm{L_h^{i_{[\ui]}-1}w_h}.
	\end{equation}
	Finally, one can complete the proof of Lemma \ref{lem:Luiwv} by substituting \eqref{eq:est-zh} into \eqref{eq:Luiwv}. 
\end{proof}

\subsection{Proof of Lemma \ref{lem:2dPig}}\label{sec:prop:2dPig}

\;\smallskip

	First, the well-definedness of $\Pig$ was prove in \cite[Lemma 2.4]{liu2020optimal}. Then, \eqref{eq:lszproj-approx} follows from the standard approximation theorem \cite{brenner2008mathematical}. See also \cite[(2.19)]{liu2020optimal}. 
	Next, for \eqref{eq:supc2d}, the case $\ell = k+2$ is a Cauchy--Schwarz version of \cite[(2.30)]{liu2020optimal}. Finally, for the case $1\leq \ell \leq k+1$, one can follow a similar argument as the standard suboptimal error estimate of the DG scheme.  To be more specific, with $\eta = \Pig w - w$, we have
	\begin{equation}\label{eq:supc2d-split}
		\begin{aligned}
			&\ip{(\Pik L - L_h\Pig)w}{v_h} = \ip{Lw}{v_h} - \ip{ L_h\Pig w}{v_h} = \ip{\eta}{\beta^x(v_h)_x+\beta^y(v_h)_y}\\
			&-\sum_{i=1}^{N_x}\sum_{j=1}^{N_y}\int_{x_{i-\hf}}^{x_{i+\hf}} \beta^y \left(\eta\left(x,y_{j+\hf}^-\right) v_h\left(x,y_{j+\hf}^-\right)-\eta\left(x,y_{j-\hf}^-\right) v_h\left(x,y_{j-\hf}^+\right)\right)\dd x\\
			&-\sum_{i=1}^{N_x}\sum_{j=1}^{N_y} \int_{y_{j-\hf}}^{y_{j+\hf}}  \beta^x \left(\eta\left(x_{i+\hf}^-,y\right) v_h\left(x_{i+\hf}^-,y\right)-\eta\left(x_{i-\hf}^-,y\right) v_h\left(x_{i-\hf}^+,y\right)\right)\dd y\\
			&\leq C\nm{\eta}\nm{v_h}_{H^1} + C\nmG{\eta}\jpnm{v_h} \leq C\left(h^{-1}\nm{\eta}+h^{-\hf}\nmG{\eta}\right)\nm{v_h}.
		\end{aligned}
	\end{equation}
	Note with the approximation property of $\eta$ and the trace inequality, it can be seen that $\nm{\eta}\leq C\nm{w}_{H^{\min(\ell,\flat)}}h^{\min(\ell,\flat)}$ and $\nmG{\eta} \leq C\nm{\eta}^{1/2}\nm{\eta}_{H^1}^{1/2}\leq C\nm{w}_{H^{\min(\ell,\flat)}}h^{\min(\ell,\flat)-1/2}$. Substituting these estimates into \eqref{eq:supc2d-split} completes the proof of \eqref{eq:supc2d} for the case $1\leq \ell\leq k+1$. 

\subsection{Proof of Lemma \ref{lem:i-PiL}}\label{sec:proof-est2}

\;\smallskip

1. \eqref{eq:I-Pi-0} follows from \eqref{eq:uni-approx} and the standard approximation theory \cite{brenner2008mathematical}. 
	
2. To prove \eqref{eq:I-Pi-1} for $\PiPi = \Pi, \Pig$, we consider 
	\eqref{eq:I-Pi-0} with $w = L \tilde{w}$, $\ell = \tilde{\ell}-1$, and $\flat = \tilde{\flat} -1$. Note that, in the case 
	$\tilde{w} \in H^{\tilde{\flat}}(\Omega_h)$, $1\leq \tilde{\ell} \leq k+2$ and $\tilde{\flat}\geq 2$, we have $w\in H^{\flat}(\Omega_h)$, $0\leq \ell \leq k+1$ and $\flat\geq 1$. Hence \eqref{eq:I-Pi-0} can be applied, which gives
	\begin{equation}
	\nm{(I-\PiPi) L \tilde{w}} \leq C \nm{L \tilde{w}}_{H^{\min(\tilde{\ell}-1,\tilde{\flat} -1)}}h^{\min(\tilde{\ell}-1,\tilde{\flat}-1)}\leq C \nm{ \tilde{w}}_{H^{\min(\tilde{\ell},\tilde{\flat})}}h^{\min(\tilde{\ell},\tilde{\flat})-1}.
	\end{equation}
	We can complete the proof of \eqref{eq:I-Pi-1} after removing the tilde symbol. The proof for $\PiPi = \Piz$ is similar. 
	
3. Similarly, for \eqref{eq:PiL-LhPig-i}, we consider \eqref{eq:piL-LhPig-g} with $w = L^i \tilde{w}$, $\ell =  \tilde{\ell}-i$, $\flat = \tilde{\flat} -i$, and $v_h = \left(\Pik L - L_h \Pig\right)L^{i}\tilde{w}^n$. Note that, in the case $\tilde{w} \in C^{\tilde{\flat}-1}(\overline\Omega) \cap H^{\tilde{\flat}}(\Omega_h)$, $1\leq \tilde{\ell} \leq k+2$, and $0\leq i \leq \min(\tilde{\ell},\tilde{\flat})-1$, we have $w\in C(\overline{\Omega})\cap H^{\flat}(\Omega_h)$ and $1\leq \ell\leq k+2$. Thus, \eqref{eq:piL-LhPig-g} holds and 
	\begin{equation}
	\begin{aligned}
	\nm{\left(\Pik L - L_h \Pig\right)L^{i}\tilde{w}}\leq& C_d \nm{L^i\tilde{w}}_{H^{\min(\tilde{\ell}-i,\tilde{\flat}-i)}}h^{\min(\tilde{\ell}-i,\tilde{\flat}-i)-1}\\
	\leq& C_d \nm{\tilde{w}}_{H^{\min(\tilde{\ell},\tilde{\flat})}}h^{\min(\tilde{\ell},\tilde{\flat})-i-1}.
	\end{aligned}
	\end{equation}
	The proof is completed after removing the tilde symbol.

4. Since $\LL_h = (I-\Pip)L_h$, adding and subtracting terms gives
\begin{equation}\label{eq:pigL-LLpig}
\begin{aligned}
(\Pik L- \LL_h \Pig) w
=& 
\left(\Pik L - L_h \Pig\right) w + \Pip L_h \Pig w  \\
=& 
\left(\Pik L - L_h \Pig\right) w + \Pip\left( L_h \Pig -\Pik L\right) w + \Pip \Pik L w  \\
=& (I-\Pip)\left(\Pik L - L_h \Pig\right) w - \Pip(I-\Pi) Lw + \Pip L w.
\end{aligned}
\end{equation}
For $1\leq \ell \leq k+1$, with \eqref{eq:I-Pi-1} in Lemma \ref{lem:i-PiL}, it can be seen that
\begin{align}
& \nm{\Pip(I-\Pi) Lw } \leq \nm{(I-\Pi)Lw } \leq C\nm{w}_{H^{\min(\ell,\flat)}}h^{\min(\ell,\flat)-1},\\
&\nm{\Pip L w} = \nm{(I-\Piz)Lw}  \leq C\nm{w}_{H^{\min(\ell,\flat)}}h^{\min(\ell,\flat)-1}.
\end{align} 
With \eqref{eq:PiL-LhPig-i} in Lemma \ref{lem:i-PiL}, we have 
\begin{equation}
\nm{(I-\Pip)\left(\Pik L - L_h \Pig\right) w}\leq \nm{\left(\Pik L - L_h \Pig\right) w}\leq C_d \nm{w}_{H^{\min(\ell,\flat)}}h^{\min(\ell,\flat)-1}.
\end{equation}
Applying the triangle inequality to \eqref{eq:pigL-LLpig}, and then substituting in the estimates above, we obtain \eqref{eq:piL-Lpi-1-1}. 

\subsection{Proof of Lemma \ref{lem:newproj}}\label{sec:proof-proj}

\;\smallskip

\begin{LEM}[Invertibility of $I+\tau \LL_h/2$]\label{lem:invert2}
	$(I+\tau \LL_h/2)^{-1}$ is well-defined when $\lambda \ll1$. Moreover, $\nm{(I+\tau \LL_h/2)^{-1}w_h}\leq C\nm{w_h}$.
\end{LEM}
\begin{proof}
	Consider 
	\begin{equation}
		\left(I+\frac{\tau}{2} \LL_h\right)v_h = w_h.
	\end{equation}
	We apply the triangle inequality to get
	\begin{equation}
		\nm{v_h} - \frac{\tau}{2}\nm{\LL_h}\nm{v_h} \leq \nm{w_h}.
	\end{equation}
	Recall that $\tau\nm{ \LL_h} \leq C \lambda$. When $\lambda\ll 1$, we have $1-\tau\nm{ \LL_h}/2\geq 1-C\lambda/2 > 0$. Hence
	\begin{equation}\label{est:I+tauLLinv}
		\nm{v_h}\leq\left(1-\frac{\tau}{2}\nm{ \LL_h}\right)^{-1}\nm{w_h}\leq  \left(1-\frac{C\lambda}{2}\right)^{-1}\nm{w_h}\leq C\nm{w_h}.
	\end{equation} 
	When $w_h = 0$, \eqref{est:I+tauLLinv} forces $\nm{v_h} = 0$, and thus $v_h = 0$. Therefore $I+\tau\LL_h/2$ is an injection. 	Note that $I+\tau \LL_h/2: V_h^k \to V_h^k$ is an operator on a finite-dimensional space. $I + \tau \LL_h/2$ is an injection implies that $I + \tau \LL_h/2$ is invertible. By rewriting $v_h = (I+\tau \LL_h/2)^{-1}w_h$ in \eqref{est:I+tauLLinv}, we obtain the desired estimate. 
\end{proof}

\begin{LEM}[Commutation error of $\LL_h$]\label{lem:est2-2}
	If $w \in C^{\flat-1}(\overline{\Omega})\cap H^\flat(\Omega_h)$ and $\flat \geq 2$, then 
	\begin{equation}
		\nm{(\Pig L- \LL_h \Pig) w}\leq C\nm{w}_{H^{\min(\ell,\flat)}}h^{\min(\ell,\flat)-1}\qquad 1\leq \ell \leq k+1.\label{eq:piL-Lpi-1}
	\end{equation}
\end{LEM}
\begin{proof}
	By triangle inequality, we have
	\begin{equation}\label{eq:pigLLhpig}
		\nm{(\Pig L- \LL_h \Pig) w}\leq \nm{(\Pik L- \LL_h \Pig) w}+ \nm{(I -\Pik ) Lw}+\nm{(I -\Pig ) Lw}
	\end{equation}
	Since $w\in C^{\flat-1}(\overline{\Omega})\cap H^\flat(\Omega_h)$ and $\flat \geq 2$, the regularity assumptions in both \eqref{eq:piL-Lpi-1-1} and \eqref{eq:I-Pi-1} of Lemma \ref{lem:i-PiL} hold. Applying these estimates to \eqref{eq:pigLLhpig} would complete the proof of Lemma \ref{lem:est2-2}. 
\end{proof}

Now we prove Lemma \ref{lem:newproj}. 
\begin{proof}[Proof of Lemma \ref{lem:newproj}]
	
	Since for $w \in C^{\flat-1}(\overline{\Omega})\cap H^\flat(\Omega_h)$ with $\flat \geq 2$, we have $w + \tau L w/2 \in H^1(\Omega_h)$. Thus $\Pig (w + \tau L w/2)$ is well-defined. By Lemma \ref{lem:invert2}, $\Pis w = (I+\tau \LL_h/2)^{-1}\Pig (w + \tau L w/2)$ is also well-defined.  
	
	To prove the approximation estimate \eqref{eq:uni-approx}, note that
	\begin{equation}\label{eq:diff-Pis-Pig}
		\begin{aligned}
			\left(\Pis  - \Pig\right) w =& \left(I+\frac{\tau}{2} \LL_h\right)^{-1}\Pig\left(I+\frac{\tau}{2} L\right)w - \Pig w\\
			=& \left(I+\frac{\tau}{2} \LL_h\right)^{-1}\Pig\left(I+\frac{\tau}{2} L\right)w - \left(I+\frac{\tau}{2} \LL_h\right)^{-1}\left(\Pig w+\frac{\tau}{2} \LL_h \Pig w\right) \\
			=& \frac{\tau}{2} \left(I+\frac{\tau}{2} \LL_h\right)^{-1}\left(\Pig L- \LL_h \Pig\right)w.
		\end{aligned}
	\end{equation}
	Applying Lemmas \ref{lem:invert2} and \ref{lem:est2-2}, for all $1\leq \ell \leq k+1$, one can get
	\begin{equation}
	\begin{aligned}
		\nm{(\Pis - \Pig) w}\leq  C\tau \nm{(\Pig L- \LL_h \Pig)w}\leq & C \tau \nm{w}_{H^{\min(\ell,\flat)}}h^{\min(\ell,\flat)-1}\\
		\leq & C \nm{w}_{H^{\min(\ell,\flat)}}h^{\min(\ell,\flat)}.
	\end{aligned}
	\end{equation}
	Using the triangle inequality and \eqref{eq:I-Pi-0}, we have
	\begin{equation}\label{eq:I-Pis-2}
		\nm{\left(I -\Pis\right)w}\leq \nm{\left(I -\Pig\right)w} + 	\nm{\left(\Pig -\Pis\right)w}\leq C \nm{w}_{H^{\min(\ell,\flat)}}h^{\min(\ell,\flat)}.
	\end{equation}
\end{proof}

\subsection{Proof of Lemma \ref{lem:newproj-general}}\label{newproj-general}

\;\smallskip

\begin{LEM}[Generalization of Lemma \ref{lem:invert2}]\label{lem:invertg}
	Let $\alpha_i = (i!)^{-1}$ for $i = 0,1, \cdots,r$. The operator inversion $\left(\sum_{i = 1}^s \alpha_i (\tau \LL_h)^{i-1} \right)^{-1}$ is well-defined if $\lambda \ll 1$. Moreover, 
	\begin{equation}
		\nm{\left(\sum_{i = 1}^s \alpha_i (\tau \LL_h)^{i-1} \right)^{-1} w_h}\leq C\nm{w_h}.
	\end{equation}
\end{LEM}

\begin{proof}
	Consider $\left(\sum_{i = 1}^s \alpha_i (\tau \LL_h)^{i-1} \right)v_h = w_h$. Note that $\sum_{i = 1}^s \alpha_i (\tau \LL_h)^{i-1} = I+\sum_{i = 2}^s \alpha_i (\tau \LL_h)^{i-1}$. Then by the triangle inequality, we have 
	\begin{equation}\label{eq:inj-general}
		\left(1-\sum_{i = 2}^s |\alpha_i|\left( \tau\nm{ \LL_h}\right)^{i-1}\right)\nm{v_h}\leq \nm{w_h}.
	\end{equation} 
	Recall that $\tau \nm{\LL_h} \leq C\lambda$. By taking $\lambda \ll1$, we can make $1-\sum_{i = 2}^s |\alpha_i| ( \tau\nm{ \LL_h})^{i-1}>0$ and hence $\nm{v_h}\leq C\nm{w_h}$. As a result, $w_h = 0$ implies $v_h = 0$. Therefore, $\sum_{i = 1}^s \alpha_i (\tau \LL_h)^{i-1}$ is an injection on a finite-dimensional space, and hence invertible. The norm estimate also follows from \eqref{eq:inj-general} by taking $v_h = \left(\sum_{i = 1}^s \alpha_i (\tau \LL_h)^{i-1} \right)^{-1} w_h$. 
\end{proof}

\begin{LEM}[Generalization of Lemma \ref{lem:est2-2}]\label{lem:estg}
	If $w \in C^{\flat-1}(\overline{\Omega})\cap H^\flat(\Omega_h)$ and $\flat\geq 2$, then 
	\begin{equation}\label{eq:piL-Lpi}
		\nm{(\Pig L^i- \LL_h^i \Pig) w}\leq C\nm{w}_{H^{\min(\ell,\flat)}}h^{\min(\ell ,\flat)-i} \qquad \forall 1\leq \ell \leq k+1, 
	\end{equation}
	where $0\leq i\leq \min(\ell,\flat-1)$. 
\end{LEM}
\begin{proof}
	We prove by induction with respect to $i$. The case $i = 0$ is trivial, and the case $i = 1$ is proved in Lemma \ref{lem:est2-2}. Now we assume \eqref{eq:piL-Lpi} holds for $i-1$, and we prove the inequality holds for $i$ with $i\geq 2$. Indeed, 
	\begin{equation}\label{eq:ind0}
		\begin{aligned}
			(\Pig L^{i}- \LL_h^{i} \Pig) w = (\Pig L^{i-1} - \LL_h^{i-1} \Pig) Lw + \LL_h^{i-1} (\Pig Lw - \LL_h \Pig w).
		\end{aligned}
	\end{equation}
	Noting that $2\leq i\leq  \flat -1$, we have $Lw\in C^{\flat-2}(\overline{\Omega})\cap H^{\flat-1}(\Omega_h)\subset C^{i-1}(\overline{\Omega})\cap H^i(\Omega_h)\subset C^1(\overline{\Omega})\cap H^2(\Omega_h)$. Using the induction hypothesis with $\ell$ replaced by $\ell-1$, we have 
	\begin{equation}\label{eq:ind1}
		\begin{aligned}
			\nm{(\Pig L^{i-1} - \LL_h^{i-1} \Pig) Lw} \leq& C\nm{L w}_{H^{\min(\ell-1,\flat-1)}}h^{\min(\ell-1,\flat-1)-(i-1)}\\
			\leq&  C\nm{w}_{H^{\min(\ell,\flat)}}h^{\min(\ell,\flat)-i} .
		\end{aligned}
	\end{equation}
	With the fact $\nm{\LL_h}\leq \nm{L_h} \leq Ch^{-1}$ and Lemma \ref{lem:est2-2}, one can get
	\begin{equation}\label{eq:ind2}
	\begin{aligned}
		\nm{\LL_h^{i-1} (\Pig Lw - \LL_h \Pig w)}\leq& Ch^{-(i-1)}\nm{(\Pig L - \LL_h \Pig) w}\\
		 \leq&  C\nm{w}_{H^{\min(\ell,\flat)}}h^{\min(\ell,\flat)-i}.
	\end{aligned}
	\end{equation}
	Finally, by applying the triangle inequality to \eqref{eq:ind0}, and invoking \eqref{eq:ind1} and \eqref{eq:ind2}, we prove \eqref{eq:piL-Lpi} for $i$. By induction, \eqref{eq:piL-Lpi} holds for all $0\leq i\leq \min(\ell,\flat-1)$. 
\end{proof}

Now we prove the approximation estimate of $\Pis$. 
\begin{proof}[Proof of Lemma \ref{lem:newproj-general}]
	\begin{equation}
		\begin{aligned}
			&\Pis w - \Pig w \\
			=& \left(\sum_{i = 1}^s  \alpha_i (\tau \LL_h)^{i-1}\right)^{-1}\left(\Pig\left(\sum_{i = 1}^q (i!)^{-1} (\tau L)^{i-1}\right)w + \tau^q\sum_{i = q+1}^s \alpha_i (\tau \LL_h)^{i-1-q{}} \Pig L^q w\right) \\
			&- \Pig w\\
			=& \left(\sum_{i = 1}^s  \alpha_i (\tau \LL_h)^{i-1}\right)^{-1}\left(\Pig\left(\sum_{i = 1}^q (i!)^{-1} (\tau L)^{i-1}w\right) + \tau^q\sum_{i = q+1}^s \alpha_i (\tau \LL_h)^{i-1-q{}} \Pig L^q w\right)\\
			&-\left(\sum_{i = 1}^s  \alpha_i (\tau \LL_h)^{i-1}\right)^{-1}\left( \left(\sum_{i = 1}^s  \alpha_i (\tau \LL_h)^{i-1}\right)\Pig w\right).
			\end{aligned}
			\end{equation}
			Combining the like terms, it yields
				\begin{equation}
			\begin{aligned}\Pis w - \Pig w 
			=& \left(\sum_{i = 1}^s  \alpha_i (\tau \LL_h)^{i-1}\right)^{-1}\left(\sum_{i = 1}^q (i!)^{-1} \tau^{i-1} \left(\Pig L^{i-1}-\LL_h^{i-1}\Pig\right) w\right)\\
			&+ \left(\sum_{i = 1}^s  \alpha_i (\tau \LL_h)^{i-1}\right)^{-1}\left(\tau^q\sum_{i = q+1}^s \alpha_i (\tau \LL_h)^{i-1-q}\left( \Pig  L^q - \LL_h^q\Pig \right)w \right).
		\end{aligned}
	\end{equation}
	Thus applying the estimate $\nm{ \left(\sum_{i = 1}^s  \alpha_i (\tau \LL_h)^{i-1}\right)^{-1}}\leq C$ in Lemma \ref{lem:invertg}, the triangle inequality, and $\tau\nm{ \LL_h}\leq \lambda$, we have 
	\begin{equation}\label{eq:est-Pis-Pig-g}
		\begin{aligned}
			&\nm{(\Pis- \Pig) w}\\
			\leq &
			C\sum_{i = 1}^q  \tau^{i-1} \nm{\left(\Pig L^{i-1}-\LL_h^{i-1}\Pig\right) w}
			+C\tau^q\sum_{i = q+1}^s \nm{\Pig  L^q w- \LL_h^q\Pig w} \\
			\leq& 
			C\sum_{i = 0}^{q}  \tau^{i} \nm{\left(\Pig L^{i}-\LL_h^{i}\Pig\right) w}.
		\end{aligned}
	\end{equation}
	We apply Lemma \ref{lem:estg} with $\ell = q\leq \min(r,k+1,\flat-1)$. With the fact $\tau = \lambda h$, one can get
	\begin{align}\label{eq:est-tauPigLi-MiPig}
		\tau^i \nm{\left(\Pig L^{i}-\LL_h^{i}\Pig\right) w}\leq & C \nm{w}_{H^{q}}\tau^ih^{q-i}\leq C\nm{w}_{H^q}h^q\quad \forall 0\leq i\leq q.
	\end{align}
	Therefore, substituting \eqref{eq:est-tauPigLi-MiPig} into \eqref{eq:est-Pis-Pig-g}, it gives
	\begin{equation}
		\nm{(\Pis - \Pig) w} \leq C\nm{w}_{H^{q}}h^{q}. 
	\end{equation}
	Finally, recall that $\nm{(I-\Pig) w} \leq C\nm{w}_{H^{q}}h^{q}$ and apply the triangle inequality. We have
	\begin{equation}
		\nm{(I - \Pis) w} \leq \nm{(I - \Pig) w} + \nm{(\Pis - \Pig) w} \leq C\nm{w}_{H^{q}}h^{q}. 
	\end{equation}
\end{proof}

\section{Conclusions}\label{sec:conclusion}

In this paper, we perform stability analysis and optimal error estimates of a variant version of the widely used RKDG method. We focus on the linear advection equation in both one and two spatial dimensions and consider RKDG schemes of arbitrarily high order. We prove that even if we drop the highest-order polynomial modes at all inner RK stages, the resulting method will still maintain the same type of stability and optimal order of accuracy as those of the original RKDG method.

Through this paper, we establish theoretical groundwork to demonstrate the following fact: Within a multi-stage RK solver for time-dependent PDEs, it is feasible to employ cost-effective and low-order spatial discretizations at specific stages. Doing so allows for the creation of a more efficient fully discrete numerical solver that maintains the same stability type and accuracy order as the original method. The mathematical tools developed in this paper could be useful for the analysis of other numerical schemes beyond the classical method-of-lines framework and benefit the future development of numerical solvers for time-dependent PDEs.

\section*{Acknowledgment}

The work of the author was partially supported by the NSF grant DMS-2208391. This work was also partially supported by the NSF grant DMS-1929284 while the author was in residence at the ICERM in Providence, RI, during the semester program ``Numerical PDEs: Analysis, Algorithms, and Data Challenges". The author would like to thank the organizers, participants, and local faculty and staff for the wonderful workshops and many inspiring discussions. The author would also like to thank College of Arts and Sciences at The University of Alabama for the teaching release through the ASPIRE program, which made the visit to ICERM possible.

\appendix

\section{Remarks on error estimates in multi-stage form}\label{sec:multistage}

In the proof of Theorems \ref{thm:err2} and \ref{thm:err-g}, we formally rewrite the sdA-RKDG scheme into a one-step one-stage scheme and construct the desirable approximation operator accordingly. However, in previous works for proving error estimates of the RKDG method, when constructing the reference function for deriving the error equation, one usually adheres to the multi-stage form of the RK method. Here, we provide remarks on the connections between the two approaches. To avoid unnecessary technicality, we will focus on the second-order case. We assume the exact solution is smooth and consider the 1D case only, for which we have $\Pik L = L_h \Pig$. We assume $\LL_h = L_h$ in this section so that the sdA-RKDG scheme formally coincides with the original RKDG scheme. \smallskip

\underline{\em Analysis in multi-stage form.} 
For the error estimate of RK2DG1 scheme, one starts with writing the exact solution $u$ as
\begin{equation}
u^{(2)} = u^n + \frac{\tau}{2} L u^n\quand
u^{n+1} = u^n + \tau Lu^{(2)} + \tau \rho^n.
\end{equation}
Next, we multiply a test function in $V_h^k$ and then integrate by parts, which yields a similar equation as that of the RK2DG1 scheme. While in the strong form, it can be written as
\begin{equation}
\Pik u^{(2)} = \Pik u^n + \frac{\tau}{2} \Pik L u^n\quand
\Pik u^{n+1} = \Pik u^n + \tau \Pik Lu^{(2)} + \tau \Pik \rho^n.
\end{equation}
After applying the relationship $\Pik L= L_h\Pig$, we have
\begin{equation}
\Pik u^{(2)} = \Pik u^n + \frac{\tau}{2} L_h\Pig u^n\quand
\Pik u^{n+1} = \Pik u^n + \tau L_h \Pig u^{(2)} + \tau \Pik \rho^n.
\end{equation}
However, if one directly works on $\Pik u^n - u_h^n$, it will lead to a suboptimal error estimate. The way to get around is to  consider $\xi_h^n = \Pig u^n - u_h^n$. To this end, we add and subtract terms on both sides to get
\begin{subequations}\label{eq:sch2-ref}
	\begin{align}
	\Pig u^{(2)} =& \Pig u^n + \frac{\tau}{2} L_h\Pig u^n + \tau \nu_h^{n,1} ,\\
	\Pig u^{n+1} =& \Pig u^n + \tau L_h \Pig u^{(2)} + \tau \nu_h^{n,2} + \tau \Pik \rho^n.
	\end{align}
\end{subequations}
Here $\nu_h^{n,1}=\tau^{-1}(\Pig-\Pi)(u^{(2)}-u^n)$ and $\nu_h^{n,2}=\tau^{-1}(\Pig-\Pi)(u^{n+1}-u^n)$. The $\nu_h^{n,\ell}$ terms here are similar to the $\tilde{\eta}$ terms in \cite{xu2020error}.

Then, one can subtract \eqref{eq:sdArkdg2} from \eqref{eq:sch2-ref} to get the error equation
\begin{subequations}\label{eq:xi2}
	\begin{align}
	\xi_h^{(2)} =& \xi_h^n + \frac{\tau}{2} L_h \xi_h + \tau \nu_h^{n,1} ,\\
	\xi_h^{n+1} =& \xi_h^n + \tau L_h \xi_h^{(2)} + \tau \nu_h^{n,2} + \tau \Pik \rho^n.
	\end{align}
\end{subequations}
Note that the error term $\xi_h^{n+1}$ satisfies an RK scheme for the linear advection equation together with a source term. Hence to derive optimal error estimates, one typically first derives a stability estimate of the RK scheme for the non-homogeneous equation $u_t + \beta u_x = f({x},t)$. Then, one applies such stability estimate to \eqref{eq:xi2}. Together with the fact $\nu_h^{n,1}$, $\nu_h^{n,2}$, and $\Pik \rho^n$ are $\mathcal{O}(\tau^2+h^2)$ terms, one can complete the optimal error estimate. \smallskip

\underline{\em Difference between $\Pig$ and $\Pis$.} The first difference we notice is that in \eqref{eq:sch2-ref}, we use $\Pig$ to define the reference function, and in \eqref{eq:pisu}, we use $\Pis$ to define the reference function. When $\LL_h = L_h$, it yields
\begin{equation}
\Pis w= \left(I+\frac{\tau}{2}L_h\right)^{-1} \Pig \left(I + \frac{\tau}{2} L\right)w. 
\end{equation}
By simple algebraic manipulation and the relationship $L_h \Pig = \Pik L$, we have
\begin{equation}
\begin{aligned}
\left(\Pis - \Pig\right)w =&\left(I+\frac{\tau}{2}L_h\right)^{-1}\left( \Pig \left(I + \frac{\tau}{2} L\right)w -  \left(I + \frac{\tau}{2} L_h\right)\Pig w\right)\\
=& \frac{\tau}{2} \left(I+\frac{\tau}{2} L_h\right)^{-1}\left(\Pig L- L_h \Pig\right)w\\
=& \frac{\tau}{2} \left(I+\frac{\tau}{2} L_h\right)^{-1}\left(\Pig - \Pik \right)Lw.
\end{aligned}
\end{equation}
One can show that, when $\tau \nm{L_h} \leq C \lambda \ll 1$, $\nm{\left(I+{\tau}L_h/2\right)^{-1}}\leq C$ and thus
\begin{equation}
\nm{\left(\Pis - \Pig\right)w} \leq C\tau \nm{\left(\Pig - \Pik \right)Lw}\leq C\nm{Lw}_{H^{k+1}}\tau h^{k+1}\leq C\nm{w}_{H^{k+2}}h^{k+2}.
\end{equation}
Moreover, we have
\begin{equation}
\nm{(\Pik L - L_h\Pis)w} = \nm{\left(L_h\Pig-L_h\Pis\right)w} \leq C\nm{L_h} \nm{\left(\Pig - \Pik \right)Lw}\leq C\nm{w}_{H^{k+2}}h^{k+1}.
\end{equation}
It can be seen that $\Pis w$ is a $\mathcal{O}(h^{k+2})$ order approximation of $\Pik w$ and the operator $\Pis$ also satisfies the superconvergence property.  Although $\Pig$ and $\Pis$ are mathematically different, they are indeed very close and both satisfy the essential properties in Proposition \ref{prop:ndPig}. \smallskip

\underline{\em Compact form of \eqref{eq:sch2-ref}.} For further comparison, we rewrite \eqref{eq:sch2-ref} into a compact form
\begin{equation}\label{eq:ref-rkdg2}
\Pig u^{n+1} = \Pig u^n + \tau L_h \Pig u^{n} + \frac{\tau^2}{2}L_h^{2} \Pig u^{n} + \tau\left(\tau L_h \nu_h^{n,1} + \nu_h^{n,2} + \Pik \rho^n\right). 
\end{equation}
Note, if one replaces $\Pig$ in \eqref{eq:ref-rkdg2} by $\Pis$, then this is of a similar form as that in \eqref{eq:pisu}. In \eqref{eq:ref-rkdg2}, one needs to show that $\tau L_h\nu_h^{n,1}$, $\nu_h^{n,2}$, and $\Pik \rho_h^n$ are $\mathcal{O}(\tau^2+h^2)$ terms. While in \eqref{eq:ref-rkdg2}, one needs to show that $\nu_h^{n}$, $\zeta_h^{n}$, and $\Pik \rho^n$ are $\mathcal{O}(\tau^2+h^2)$ terms. In both cases, we write the reference function as $\PiPi u^{n+1} = R_{h,2} \PiPi u^n + \tau \mathcal{O}(\tau^2+h^2)$, where $\PiPi = \Pig, \Pis$.  \smallskip

\underline{\em Multi-stage form of \eqref{eq:pisu}.} When $\LL_h = L_h$, one can rewrite \eqref{eq:pisu} into a multi-stage form.
\begin{subequations}\label{eq:ref2-rkdg2}
	\begin{align}
	\Pis u^{(2)} =& \Pis u^n + \frac{\tau}{2}L_h\Pis u^n,\\
	\Pis u^{n+1} =& \Pis u^n + \tau L_h \Pis u^{(2)} + \tau\left(\Pik \rho^n+\nu_h^n + \zeta_h^n\right).
	\end{align}
\end{subequations}
Note, if one replaces $\Pis$ in \eqref{eq:ref2-rkdg2} by $\Pig$, this is of a similar form as that in \eqref{eq:pisu}. Only the source terms are distributed differently across stages. 

\bibliography{refs}
\bibliographystyle{abbrv}

\end{document}